# Librationism & its classical and extraclassical set theories

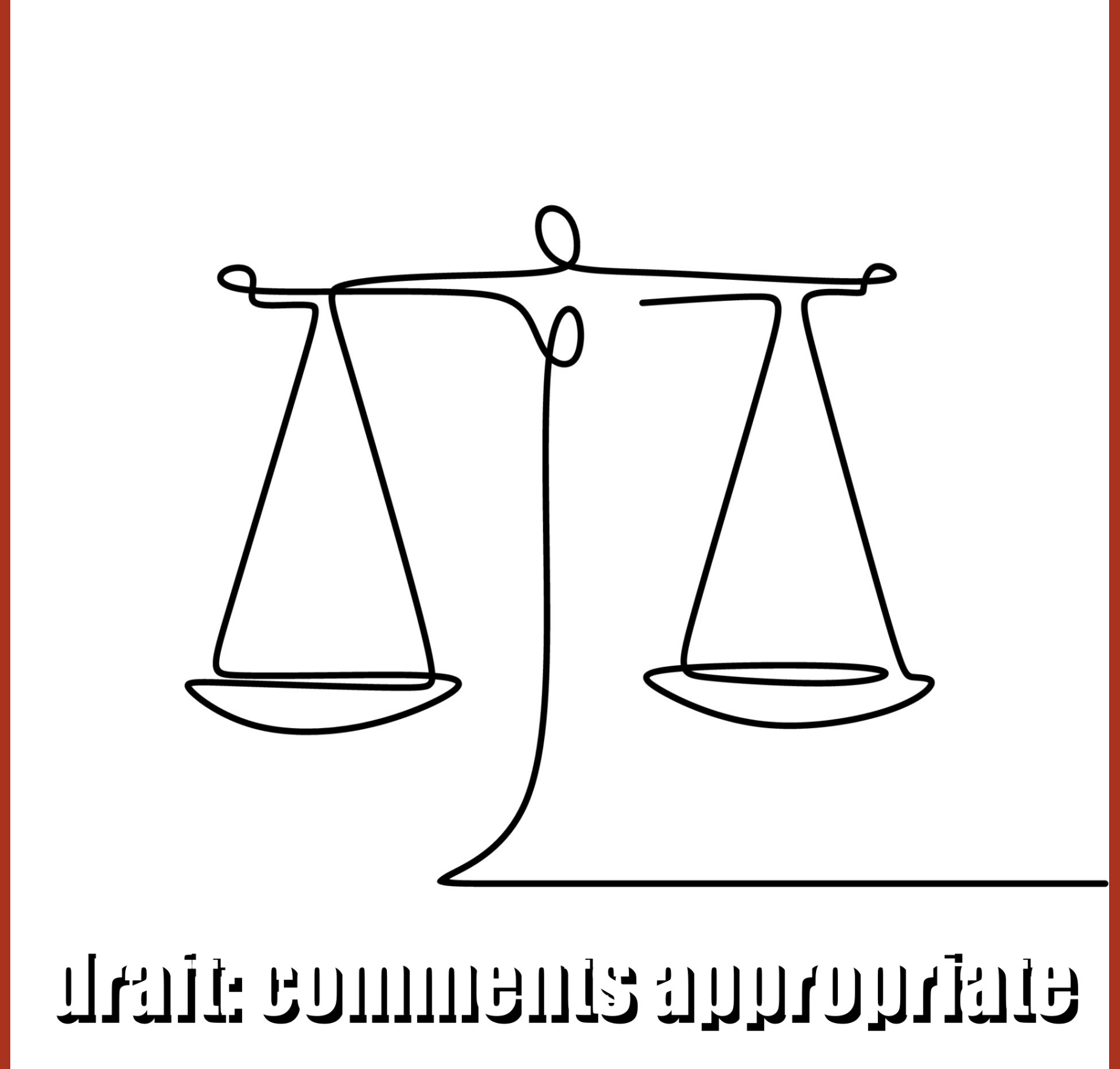

Responsible Extraclassicism, Paradox, Truth and Reality

Frode Alfson Bjørdal

...............................

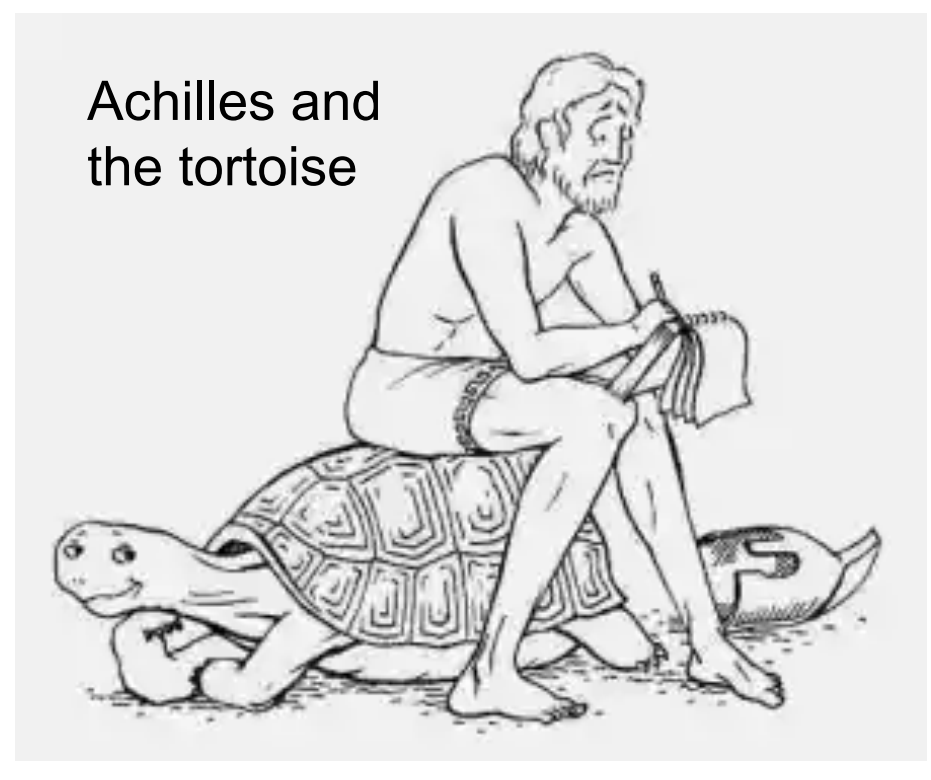

deliberating on modus ponens.

# LIBRATIONISM &

## its classical and extraclassical set theories

RESPONSIBLE EXTRACLASSICISM – PARADOX, TRUTH AND REALITY

Frode Alfson Bjørdal

*To my godson Jon Lervåg Walderhaug, 2/11/2012,*

*& my dear friend, nephew Jon Vegard Lervåg, 2/17/1979 --- 7/22/2011.*

# 1 Preface

*One cannot see it all from one point of view.*

---

1.0.1. Theorem Title

text

*NB! As indicated on the front page, this is still a draft. So there are still many things in the manuscript which need to be tidied up.*

The investigations, which led to *librationism* began in the spring of 1993. It was a struggle to build upon imprecise thoughts to express beliefs precise enough to be useful. But the author pressed on, as the investigations were addictive, and began giving programmatic talks already with (Bjørdal 1996b), and publishing unfinished ideas with (Bjørdal 1998). Motivation to persist was always found, so the hard work continued over all those years, in between other rather taxing tasks.

The problems with the paradoxes in set theory are of great importance not only when one attempts to find useful and philosophically reasonable foundational theories for the formal sciences, but also when one seeks to account for a variety of important problems in metaphysics, epistemology and other areas of philosophy.

Paradoxes, are important in metaphysics. As stressed in (Grim 1991), thinking according to the contemporary book has the awkward consequence that one must think that the world is not complete. It will be seen that librationism is not according to the current book, for it has, as shown in § 24.6, the result that there are only denumerably many objects in the world; importantly, the *validity* of Cantor's arguments for uncountability is of course not challenged.

Oddly enough, adherents of paraconsistency did not, as the author, even since before (Bjørdal 2005), explore the point of view that Cantor's Theorem in reality camouflages a paradoxicality. On the contrary, one saw that e.g. (Z. Weber 2012) attempted to retain the force of Cantor's arguments for uncountable infinities, and that others who championed paraconsistency also accepted the uncountability paradigm. Nevertheless, (Petersen 2024) has brought out that Cantor's Theorem indeed *is* paradoxical according to the set theory suggested by (Z. Weber 2012); I do not know that there are other paraconsistent set theories which may address the question differently.

We can see the more comprehensive relevance of the paradoxes rather directly from the fact that modal logics are very important philosophical tools for reasoning about ethics, knowledge, and other central philosophical concepts. But threats of paradox undermine the use of reasonable modal logics, with more than just a minimum of linguistic resources and plausible closure principles. This is on account of such limitative results as were discovered by (Montague 1963), and investigated further by others, like (Friedman and Sheard 1987), and (Cantini 1996). A takeaway is for example that if a modal logic is expressive enough to license the inference from *Smith ought to sell his house* to *there is something Smith ought to do*, then Russell like paradoxes arise.

So paradoxes are pervasive in philosophy. Standard evasions of Montague like limits put restrictions upon the linguistic resources, and the librationist resolutions are preferable as such restrictions are not needed.

The focus in this essay will be upon the paradoxes in the context of mathematics. The extension ℔ of **£**, characterized by Definition 25.0.1, succeeds with defining very strong theories. One of the central results of this treatise is that the librationist set theoretic extension ℔$\mathcal{HR}$(**G**) of ℔ accounts for **Neumann-Bernays-Gödel** set theory with the **Axiom of Choice** and **Tarski's Axiom**. Moreover, ℔ succeeds with defining a manifestation set **W**, *die Welt*, which is such that the extension ℔$\mathcal{H}$(**W**) of ℔ accounts for Quine's (Quine 1937) set theory NF, probably best known as **New Foundations**.

The points of view developed in this treatise entail that the truth-paradoxes and the set-paradoxes, *pace* Ramsey, are often as two sides of the same coin, via *Alethic comprehension*, so that the librationist resolutions of the set theoretic paradoxes are at the same time resolutions of the truth theoretic paradoxes. Both the librationist resolutions of the set theoretic paradoxes and the truth theoretic paradoxes have non-trivial philosophical implications: librationist set theories have the consequence that there are no absolutely uncountable sets, and librationist truth theories allow the use of syntactical modalities in ways which circumvent limitations as those of (Montague 1963), and a truth predicate which is useful for more precise philosophical discourse.

# 2 Others, times and their places

> “However, not everything that can be counted counts, and not everything that counts can be counted.”
>
> William Bruce Cameron in (Cameron 1963, p. 13)

.

Let me before being more specific extend my thanks to all those who taught me, and to all those who communicated with me, about some of these matters which interested me so greatly! Thanks as well to all those who helped me with bibliographical and technical problems!

The investigations, which led to **£**, and its strengthenings, began in the spring of 1993. Many thanks to colleague, professor Ingemund Gullvåg, for interesting conversations about these problems when I worked at the University of Trondheim from 1992 to 1994. It was a struggle to build upon imprecise thoughts to express beliefs precise enough to be useful. But the author pressed on, as the investigations were addictive. I began giving programmatic talks already with (Bjørdal 1996b). Thanks to Dag Prawitz for friendly comments on the failed project I presented in Uppsala in 1996. I continued publishing unfinished ideas with (Bjørdal 1998), and am thankful for the comments I received from Krister Segerberg, who attended my presentation in Villa Lanna, in Prague, September 1997.

Motivation to persist was always found, so the hard work continued over all those years, in between other rather taxing tasks.

The problems with the paradoxes in set theory are of great importance not only when one attempts to find useful and philosophically reasonable foundational theories for the formal sciences, but also when one seeks to account for a variety of important problems in metaphysics, epistemology and other areas of philosophy.

Paradoxes, are important in metaphysics. As argued, successfully, I think, in (Grim 1991), thinking according to the contemporary book has the awkward consequence that one must think that the world is not complete. It will be seen that librationism is not according to the current book, for it has, as shown in § 24.6, the result that there are only denumerably many objects in the world; importantly, the *validity* of Cantor’s arguments for uncountability is of course not challenged. Moreover: librationism has a universal set in a framework which respects classical logic. Instead, it on the other hand turns out that the standard power set operation creates paradoxes when it is not of an unparadoxical universal set! For more on this, see §19.

Oddly enough, adherents of paraconsistency did not, as the author, even since before the third millennium, explore the point of view that Cantor’s Theorem in reality camouflages a paradoxicality. On the contrary, one saw that e.g. (Z. Weber 2012) attempted to retain the force of Cantor’s arguments for uncountable infinities, and that others who championed for paraconsistency also accepted the uncountability paradigm. Nevertheless, (Petersen 2024) has recently brought out that Cantor’s Theorem indeed *is* paradoxical according to the set theory suggested by (Z. Weber 2012); I do not know that there are paraconsistent set theories which may address the question differently. Uwe Petersen remarks, in a private email of Februar 25, 2025, that he is not aware of any replies to (Petersen 2024). Many thanks to Uwe Petersen, for several useful conversations which we have had since the end of the 1990s. Many thanks as well to Zack Weber, for some discussions we had at *The Fourth World Congress of Paraconsistency*, Melbourne, Australia, 13-18 July 2008.

We can see the more comprehensive relevance of the paradoxes rather directly from the fact that modal logics are very important philosophical tools for reasoning about ethics, knowledge, and other central philosophical concepts. But threats of paradox undermine the use of reasonable modal logics, with more than just a minimum of linguistic resources and plausible closure principles. This is on account of such limitative results as were discovered by (Montague 1963), and investigated further by others, like (Friedman and Sheard 1987), and (Cantini 1996). A takeaway is for example that if a modal logic is expressive enough to license the inference from *Smith ought to sell his house* to *there is something Smith ought to do*, then Russell like paradoxes arise.

Thanks to Andrea Cantini for conversations on email, and at the Universal Logic conference in Lisbon, April 18-25, 2010, and for authoring the important book (Cantini 1996). Thanks also to Harvey Friedman for conversations on email, and at the Workshop in Ghent, the Netherlands, upon his receiving an Honorary doctorate at the University of Ghent.

The focus in this essay will be upon the paradoxes in the context of mathematics.

The extension ℔ of **£**, determined by Definition 25.0.1, succeeds as a useful framework for defining further extensions which turn out to be strong set theories, as compared with ZF. One of the central results of this treatise is that the librationist set theoretic extension ℔$\mathcal{H}\mathcal{R}(\mathbf{G})$ of ℔ accounts for **Neumann-Bernays-Gödel** set theory with the **Axiom of Choice** and **Tarski's Axiom**; which is useful as it accounts naturally for category theory. Moreover, ℔ succeeds with defining a manifestation set **W**, *die Welt*, which is such that the extension ℔$\mathcal{H}(\mathbf{W})$ of ℔ accounts for Quine's (Quine 1937) set theory, best known as **New Foundations**.

Thanks to my logic teacher, professor Kåre Johnsen, at the University of Bergen, in the period 1984 to 1985; and to my philosophy teacher, professor Hans Skjervheim, at the same university, in the period from 1983 to 1984.

The academic year 1985/86 I studied, with a DAAD-scholarship, at the Goethe University in Frankfurt am Main. I profited from seminars by professors Karl Otto Apel and Jūrgen Habermas on theories of truth, and from conversations on a range of logic topics with professor Wilhelm K. Essler. One thing that I was preoccupied with then was with the so-called *Reflexive argument* of Karl Otto Apel, which he suggested as a remedy to obtain epistemical certainty on issues such as the principle of non-contradiction, and seemingly other statements such as the *sum*—conclusion of the Cartesian *cogito-ergo-sum*—argument. I had a lot of contact with Privatdozent Wolfgang Kuhlmann, who had just published his Habilitationsschrift (Kuhlmann 1985), and wrote an essay, *Eine Kritik der Tranzendentalpragmatik*. He kindly supervised its progress, as I presented it for an oral exam upon which I was accepted as a *Doktorand*, in the same way as German students with Magister-degrees. I then discovered the book (Rescher and Brandom 1980), in the library, and it contributed to me becoming rather skeptical concerning the viability of Apel's transcendental pragmatism. But the simple realization that we may state that something elsewhere is contradictory, without committing to the principle of non-contradiction by stating that, was in the end more important. At any rate, I became very much interested in the principle of non-contradiction, then, and, as well, very much interested in the paradoxes.

The points of view developed in this treatise entail that the truth-paradoxes and the set-paradoxes are sometimes as two sides of the same coin, via *Alethic comprehension*. Both the librationist resolutions of the set theoretic paradoxes and the truth theoretic paradoxes have non-trivial philosophical implications: librationist set theories have the consequence that there are no absolutely uncountable sets, and librationist truth theories in adequate contexts allow the use of syntactical modalities in ways which may circumvent limitations as those of (Montague 1963), so librationism has a truth predicate which is useful for precise philosophical discourse.

I did not end up studying for my doctorate in Frankfurt am Main, but instead attended the PhD program in Philosophy at the University of California, Santa Barbara, from 1988 to 1992. I also took classes at University of California, Los Angeles, as that was allowed. Many thanks to my advisor Anthony Brueckner for his help with my dissertation *A Theory of Knowledge*, which was written to deal with the Gettier problem. Lots of thanks to my teacher Nathan U. Salmon, for many conversations and discussions, and for convincing me, through his book (Salmon 1986), that his Millian theory of names is correct for belief contexts - as testified by my publication (Bjørdal 1996a); unfortunately the *Nordic Journal of Philosophical Logic* did not at the time offer me the opportunity to edit the article before it was sent to press, so it unfortunately lacked the acknowledgement which it should have had, so I hope this somewhat makes up for it! Thanks also to my teacher William Forgie, for his seminar in the philosophy of religion, which was an important reason why I ended up writing (Bjørdal 2018). Thanks as well to David Kaplan, for his seminar on demonstratives at UCLA. I am also very grateful to professor Alonzo Church for teaching me in his last doctoral seminar on logic, with a concluding oral exam, in the spring of 1989, as he turned 86 years old. Many thanks also to mathematics professor John Doner at UCSB, for many conversations on the paradoxes, and for his seminar on Gentzen's consistency proof for arithmetic, in the spring of 1987.

In some semesters as from 2013 I held some seminaries on a variety of topics, for master and doctorate students of Philosophy, as *Professor Colaborador Voluntário*, at Programa de Pós-Graduação em Filosofia na Universidade Federal do Rio Grande do Norte, in Natal, Brasil: In the second semester of 2013 the seminary *The librationist foundation for reasoning* was offered, and in the second semester of 2014 the seminary *A teoria libracionista das coisas*.

Talks were held for the LOGICA congresses of the Czech Academy of Sciences, in 1998, 2004, 2005 and 2010. Associated papers (Bjørdal 1999; Bjørdal 2005; Bjørdal 2006; Bjørdal 2011) were published in the LOGICA Yearbook series. I extend my heartfelt thanks for stimulating interaction to the Czech participants Jaroslav Peregrin, Pavel Materna, Vít Punčochář, Andrej Majer, Vladimír Svoboda, Vítězslav Švejdar, Petr Cintula, Marta Bílková, Carles Noguera, Marie Duží, Karel Berka, Libor Běhounek, Ivan Kramosil, Petr Hájek & to the international participants Igor Sedlár, Juraj Podroužek, Göran Sundholm, Hartley Slater, James A. Woodbridge, Melvin Fitting, Philip Kremer, Franz Huber, Jordan Howard Sobel, Gerhard Schurz, Colin Cheyne, Catarina Dutilh Novaes, Gary Kemp, Michelle Friend, Jan Woleński, Yuri Gurevich, Kosta Došen, Massimiliano Carrara, Shawn Standefer, Theodora Achourioti, Gil Sagi, Andrea Iacona, Hannes Leitgeb, Albert (Tyke) Nunez, David Ripley, Martin Pleitz, Sebastian Sequoiah-Grayson, Peter Verdee, Peter Verdee, Giuseppina Ronzitti, Tero Tulenheimo, Raymond Turner, Arianna Betti, Maria van der Schaar, Wlodek Rabinowitz, Scott A. Shalkowski, Timothy Childers, Tomis Kapitan, Fred Johnson, Otávio Bueno, Robert Brandom, Giuseppe Primiero, John T. Kearns, Emil Bădici, Mircea Dumitru, Douglas Patterson, Greg Restall, Bjørn Jespersen, Ole Thomassen Hjortland, Franz Berto, Stephen Read, J. C.

Beall, Andreas (Pietz) Kapsner, David Makinson

Lectures were held for the Logic Colloquia, under the auspices of the Association of Symbolic Logic: in Barcelona, 2011 – Helsinki, 2014 – Stockholm, 2017 and Prague, 2019; for World Congresses on Paraconsistent Logic: Melbourne, 2008 and Kolkata, 2014; for International Conferences on Non-classical Mathematics: Guangzhou, 2011 and Vienna, 2014; for Sociedade Brasileira de Lógica: Petrópolis, 2014 – Pirenópolis, 2017 and Salvador, 2022; for World Congresses on Logic, Methodology, and Philosophy of Science: Helsinki 2014 and Prague, 2019; and for the World Congresses on Universal Logic in Lisbon, 2010, Rio de Janeiro, 2013 and Istanbul, 2015. The Kolkata lecture was reworked and published as (Bjørdal 2015). A lecture was delivered for Logic, Reasoning and Rationality at Centre for Logic & Philosophy of Science, Ghent University, Belgium, in 2010. Many thanks to Arnon Avron, Andrea Cantini, Alasdair Urquhart, Philip Welch, Jan von Plato, João Marcos, Pierluigi Minari, Sara Negri, Graham Priest, Melvin Fitting, J. Michael Dunn, Greg Restall, Ed Mares, Benedict Löwe, Heinrich Wansing, Carlos Caleiro, Alexandre Costa-Leite, Jean-Yves Béziau, Jaime Ramos, Hartry Field, Giovanni Sambin, Gerhard Jäger, Marcus Kracht, Sergei Artemov, Elena Y. Nogina, Denis Igorevich Saveliev, Lev Beklemishev, Max Urchs, Samuel Gomes, Reinhard Kahle, John T. Baldwin and others for positive interactions, and for interesting discussions.

An invited lecture was given for the mathematics department at the Federal University of Bahia, Salvador, Brazil, on October 9th 2014. Many thanks to Samuel Gomes, Andrey Bovykin and Andreas Brunner for friendly and useful discussions then.

Here comes acknowledgements in conection with *Encontros Brasileiros de Lógica*.

On October 16th, 2017, an invited lecture was held, on the First Online publication of my article (Bjørdal 2018). Also, an invited lecture was held at the Steklov Mathematical Institute with the Moscow division of the Russian Academy of Sciences, on the 4th of September 2014. Moreover, a talk was accepted for the Third St. Petersburg Days of Logic and Computability, August 24-26 at the Euler International Mathematical Institute at the Saint Petersburg division of Steklov Mathematical Institute, in 2015. Many thanks to Yuri Matiyasevich, Lev Beklemishev, Fedor Pakhomov, Denis Igorevich Saveliev, Sergei Soloviev and others for friendly discussions on those occasions.

During my years as acting professor at the University of Oslo many lectures on Logic and mathematics were offered by guests, and I am glad for the interaction on such occasions with Dana Scott, Patrick Suppes, Per-Martin Löf, Krister Segerberg, Göran Sundholm and Saul Kripke.

I have benefitted from exchanges with many helpful research mathematicians at the scientific discussion site MathOverflow for more than twelve years. So great thanks to Alec Rhea, Alex Kruckman, Andreas Blass, Andrej Bauer, Andrés E. Caicedo, Asaf Karagila, Branden Fitelson, Carl Mummert, Carlo Beenakker, Colin McLarty, David Madore, Denis Roman Hirschfeldt, Emil Jeřábek, Farmer Schlutzenberg, Fedor Pakhomov, François G. Dorais, Gabe Goldberg, Hanul Jeon, James E. Hanson, Jason Rute, Jim Hefferon, Joel David Hamkins, Julia Kameryn Williams, Martin Sleziak, Mauro Allegranza, Michael Greinecker, Monroe Eskew, Nate Eldredge, Nik Weaver, Noah Schweber, Noel Vaillant, Paul Taylor, Peter LeFanu Lumsdaine, Rodrigo Freire, Sam Hopkins, Srinivasan Keshav, Terence Tao, Timothy Chow, Todd Trimble, Trevor Wilson, Vladimir Dotsenko, William Russell Easterly, Zhen Lin Low, Zuhair Al-Johar and Zuzana Haniková! If one studies the more than 100 questions I posted over the years on Mathoverflow, one will find that all the persons

I just thanked – and I hope that I did not forget anyone – answered and commented on questions I posted there, and the communications have been good learning experiences.

I received important bibliographical help by persons who responded to my queries on the History of Science and Mathematics site at Stack Exchange: Georg Essl, user6530, njuffa, BakerStreet, Mauro Allegranza and Pugs. Many thanks for your contributions!

Many thanks to Karl Yngve Lervåg, for lots of help with writing in LaTeX.

I also read, and sometimes took part in discussion on the Foundations of Mathematics email-list for some years. I have have not yet reconstructed my contributions.

Most importantly, I delivered many lectures for the Seminary in Logic of the Mathematics Department at the University of Oslo, where I have had the honor to talk about these matters virtually every semester since I started working at the University of Oslo, in 1996. Those seminaries were, and are, very useful, and I learned a lot from participating there, with colleagues and advanced students. Let me extend my thanks to Jens Erik Fenstad, Roger Antonsen, Stål Aanderaa, Dag Normann, Herman Ruge Jervell, Lars Kristiansen, Eivind Briseid, Jon Henrik Forssell, Bjørn Kjos-Hanssen, Øystein Linnebo, Laura Crosilla, André Rognes, Juvenal Murwanashyaka, David Ross, Davide Sutto, Arild Torolv Søetorp Waaler, Mathias Barra and Tom Louis Lindstrøm for insightful discussions over all those years!

All encounters with the mentioned persons, were very important for the mathematical and philosophical maturation of the author. Resistance was usually useful. Needless to say, but no one shares responsibility for undetected errors.

# 3 Introduction

*"Thus mathematics may be defined as the subject in which we never know what we are talking about, nor whether what we are saying is true."*

---

`Bertrand Russell`, in (Russell 1901, p. 84).

It is presupposed that $A \land \neg A$ is a contradiction, and that a theory is inconsistent just if it has contradictory theses. As per § 13.3, **£** is consistent and not contradictory. It follows that the librationist points of view are not dialetheic, for *dialetheism* is canonically characterized, in (Priest, Berto, et al. 2024), as a view which takes some contradictions to be true. Moreover, **£** is not a paraconsistent point of view, as the latter are not serene in the sense of Definition 13.3.14. Librationism is instead, as per § 13.3, taken to offer classical and *extraclassical* points of view. To complete the terminological distinction, take librationism to also offer a *bialethic* point of view, and not a dialetheic one.

It will be showns in § 16.7 that Librationism meets a challenge which it is difficult to see can be met if one presupposes that contradictions, as $p \land \neg p$, are true, viz. to offer an account of what a true sentence $p$ says, in a paradoxical situation, which its true negation $\neg p$ does not say in that situation.

A remark on designator is called for. One might hold that a theory is not a set theory if it does not presuppose exactly the same linguistic resources as the language of set theory according to what has been the *practice* after Zermelo, viz. the language of first order logic plus the symbol $\in$. This austere tenet is not abided by here, and it is instead presupposed that set theoretic reality should be investigated with such resources which best reveal it. As will become clear, *set abstracts* are used, and these are not eliminable, due to the fact that **£** is a highly non-extensional theory. The symbol $\in$, however, *is* eliminable, by means of apposition.

As the slight **£** extension ₶, pronounced (”libra”), with some more natural strengthening assumptions interprets classical set theory **NBG** plus the rather strong Tarski's Axiom, it would seem rather disingenious, to hold that **£** and ₶ are not themselves set theories.

# 4 Librationism and its formal language

*A true paradoxicality has a true complement.*

## 4.1 Nomenclature

**£** has been used to denote the librationist foundational system after the publication of (Bjørdal 2012). It will be indicated, as in §§ 25–26, when additional assumptions are made, which result in the use of ℔ for an extension of **£**.

The sign **£** is most prevalently used for the currency of the United Kingdom. It derives from Latin *libra pondo*. The last word is an adverb which means *by weight*. *Libra* was used for the Roman pound - which was about 327 g, and also for scales and balances. Such scales were as well an attribute to the Greek Goddess for Divine Justice, Themis, and for her daughter Dike, who was the Goddess for Human justice. The roles of the attributes were thought to be the weighing of the consequences of acts to find balance, and, therefore, justice. According to the *interpretatio Romana*, the Goddess corresponding to Themis and Dike in the Roman religion was the blinded *Justitia* who also had a scale as attribute. Justitia is also referred to as *Lady Justice* in English.

In the context of librationism, **£** may be taken to symbolize the weighing and gauging of balances between sentences, and perhaps most interestingly, from the librationst points of view, in the case of sentences which are *incompatible* or *complementary*, in the sense of Definitions 13.4.1 and 13.4.3.

## 4.2 Numeralism - the metanumber standpoint

The *metanumbers* are the *numbers-of-the-meta-language*. The *ordinal* metanumbers are defined à la von Neumann by means of the meta mathematical *metasettheory*, which is the set-theory-of-the-meta-language. One must distinguish ordinal metanumbers from ordinal numbers of the *set theories* expressed, and accounted for, by the object language. The *metanumerals* are *numerals-of-the-meta-language*, denoting finite ordinal metanumbers. The *natural* metanumbers are the finite ordinal metanumbers, and the metaset of *counting* metanumbers is the metaset of the natural metanumbers minus 0. The *integer* metanumbers are the natural metanumbers extended with the negative counterparts of the counting metanumbers.

The metanumber standpoint presupposed here is stronger than the view presupposed by (Gödel 1931), which was that formulas, and expressions akin, may be *correlated* via a coding with numerals denoting natural numbers. For the symbols and expressions of **£** are taken to *be* counting metanumbers, and their syntactical manipulations are accounted for by the metasettheory presupposed.

## 4.3 The inclusion of set forming abstracts

The inclusion of term, or set forming abstracts is a trait shared with (Gandy 1959), and with contributions to the literature on non-classical set theories, including some which were called *property* *theories*.

The term "property theory" seems to have an origin with Gödel, as the opening sentence of Roger Myhill's article *Paradoxes*, in Synthese 60 (1984), 129-143, is: "Gödel said to me more than once "There never were any set-theoretic paradoxes, but the property-theoretic paradoxes are still unresolved"; and he may well have said the same thing in print."

(Cantini 1988) is following up on the suggestion by Myhill, and introduced a property theory, conjointly with classical set theory, in (Cantini 1996), where (Gilmore 1974) is as well related as contribution of , as a property theory where abstracts were used because the principle of extensionality fails.

Was Gödel aware of the contribution in (Scott 1962), or did he study (Friedman 1973). (Shapiro 1985) is another witness to modern research into set theories without extensionality.

Nevertheless, there are now so many non-extensional set theories in the literature, beyond attempts to deal with the paradoxes, that it seems unreasonable to consider them *property theoretic*, as opposed to *set theoretic*.

## 4.4 "=" and "∈" are not primitive in £

The formal language of **£** is Polish, and without symbols for identity or membership. A Polish ↓-connective is used. The membership relation can be defined with apposition of terms, because there as a consequence of the Polish policy are no parentheses in the formal language of **£**.

Moreover, § 9 shows that the identity of a and b can be delineated adequately by the statement that b is an element of all sets that have a as an element, as in Definition 9.1.2.

## 4.5 Metalinguistic conventions

4.5.1. Definition of some important meta language symbols:

(1) **Σ** is the existential quantifier.

(2) **Π** is the universal quantifier.

(3) $\sim$ is for negation.

(4) & is for conjunction.

(5) ⓡ is for disjunction.

(6) $\Rightarrow$ is for implication.

(7) $\Leftrightarrow$ is for bi-implication.

(8) [x: ... ] is the set notation for use in the metalanguage.

(9) ε is the metalanguage symbol for membership.

(10) $==$ is for metamathematical identification and definition.

(11) $\alpha, \beta, \gamma, \delta, \ldots$ are for ordinal numbers of the metalanguage.

(12) $\prec, \preceq, \succeq$, and $\succ$ are the orderings on the ordinal numbers of the meta language.

(13) $\mu$ is for the least operator of the metalanguage.

4.5.2. Definition of notation for the finite metanumbers and the integers:

(1) $\Omega \Longleftarrow 0,1,2\ldots$ is the term for finite ordinal metanumbers, i.e. natural metanumbers.

(2) $\Omega_{+} \Longleftarrow 1,2,3...$ is the term for positive natural numbers, i.e. *counting metanumbers*.

(3) $\Omega^{\pm} \Longleftarrow 0,1,-1,2,-2\ldots$ is the term for the metaintegers.

(4) $\Omega_{-} \Longleftarrow -1,-2,-3\ldots$ is the term for the negative metaintegers.

(5) $\Omega^{-} \Longleftarrow 0,-1,-2\ldots$ is the term for the metaintegers which are not positive.

4.5.3. Definition of the primitive symbols, their metanumerals and metanumbers

(1) •

(2) v̈

(3) ↓

(4) ∀

(5) ς

(6) #

(7) $\mathcal{T}$

(8) $\mathcal{F}$

are the symbols, which stand for the metanumbers denoted by the bijective base-8 digit metanumerals 1, 2, 3, 4, 5, 6, 7, and 8, respectively.

4.5.4. Definition of the bijective base-8 metanumeral strings:

(1) The base-8 digit metanumerals are the metanumeral strings of length 1.

(2) Let $n_i$ and $n_j$ be base-8 metanumeral strings.

(3) $\ell(n_i) = \!= \lfloor \log_8(n_i+1) \rfloor$ invokes the floor function $\lfloor\ \rfloor$, and defines the length of the bijective base-8 metanumeral needed to express the number denoted by metanumeral $n_i$.

(4) Concatenation $\frown$ is the function given by $n_i \frown n_j = \!= n_i \cdot 2^{\ell(n_j)} + n_j$.

(5) We know that $\frown$, so defined, is associative.

(6) If $\sigma_k$ is a metanumeral string of length $m$ and metanumeral string $\sigma_l$ has length $n$, then metanumeral string $\sigma_k\,\sigma_l = \sigma_i \frown \sigma_j$ has length $m + n$

(7) $n_i \frown n_j$ is taken to be the denotatum of the apposition $n_i\, n_j$.

4.5.5. Definition of *expression*:

(1) A string of symbols from Definition 4.5.3 is an *expression*, in *symbolic form*, just if it is formed according to the formation rules.

(2) A bijective base-8 metanumeral is an expression in *canonical form* just if it corresponds with it's symbolic form via the coding of metanumerals into symbols, as given in Definition 4.5.3.

4.5.6. Definition of predicates for terms, formulas, sentences and expressions:

**Ve**(v), **Ct**(c), **Tm**(a), **Fa**(A), **Se**(B), **En**(X) are written to state that v is a variable, c is a constant, a is a term, A is a formula, B is a sentence and X is an expression.

4.5.7. Definition of the underlines:

To remind that expressions in the last analysis are metanumbers, denoted by metanumerals, we in the remainder of this section underline, and write $\underline{\text{expression}}$, $\underline{\text{variable}}$, $\underline{\text{constant}}$, $\underline{\text{term}}$, formula, sentence, and expression. To ease the reading, the underlines will not be used as from the next section.

4.5.8. Definition of variables:

(1) $\ddot{\mathrm{v}}$ is a variable.

(2) A variable succeeded by $\bullet$ is a variable.

(3) $\mathrm{v}_0$ is variable $\ddot{\mathrm{v}}$, and $\mathrm{v}_{n+1}$ is variable $\ddot{\mathrm{v}}_n \frown \bullet$.

(4) Nothing else is a variable.

(5) Variables are terms.

4.5.9. Definition of primitive constants

(1) $\mathcal{T}$ and $\mathcal{F}$ are the primitive *set* constants.

(2) Primitive constants are terms.

4.5.10. Definition for arbitrary term $\mathrm{a}_i$ and arbitrary formula $\mathrm{A}_i$:

(1) If $\mathrm{a}_i$ is a variable or a primitive set constant then $\mathrm{a}_i$ is a term.

(2) If $\mathrm{a}_0$ and $\mathrm{a}_1$ are terms, $\mathrm{a}_1\mathrm{a}_0$ is a formula.

(3) If $\mathrm{A}_0$ and $\mathrm{A}_1$ are formulas, $\downarrow \mathrm{A}_0\mathrm{A}_1$ is a formula.

(4) If $\mathrm{A}_0$ is a formula and $\mathrm{v}_0$ is a variable, $\forall \mathrm{v}_0\mathrm{A}_0$ is a formula.

(5) If $\mathrm{A}_0$ is a formula and $\mathrm{v}_0$ is a variable, $\varsigma \mathrm{v}_0\mathrm{A}_0$ is a term.

(6) Nothing else is a term or a formula.

(7) Just terms and formulas are *expressions.*

4.5.11. Definition concerning suppression of subscripts

When possible, $\mathrm{a}, \mathrm{b}, \ldots$ is written for $\mathrm{a}_0, \mathrm{a}_1, \ldots$, while $\mathrm{v}, \mathrm{w}, \ldots$ are for $\mathrm{v}_0, \mathrm{v}_1, \ldots$, and $\mathrm{m}, \mathrm{n}, \ldots$ for $\mathrm{n}_0, \mathrm{n}_1, \ldots$, along with $\mathrm{A}, \mathrm{B}, \ldots$ instead of A with subscripts. Other letters, or letter-like symbols, may be used for special terms, or formulas.

4.5.12. Definition of binders, binds, ties and scopes:

(1) In $\forall \mathrm{v}\mathrm{A}$, $\forall$ is the *binder*. v is the *bind* of A and the *tie* of $\forall$. A is the *scope* of $\forall$.

(2) In $\varsigma \mathrm{v}\mathrm{A}$, $\varsigma$ is the *binder*. v is the *bind* of A and the *tie* of $\varsigma$. A is the *scope* of $\varsigma$.

4.5.13. Definition of free and bound variables

(1) A variable occurrence in a formula, or term, is bound, just if it is a bind, or it is in the scope of a binder with another occurrence as tie.

(2) Variable occurrences in a formula, or term, are *free* if not bound.

(3) A variable is free in a formula, or term, just if an occurrence is.

(4) A variable is bound in a formula, or term, just if an occurrence is.

4.5.14. Definition of sentences and constants:

(1) A term without free variables is a *constant*.

(2) A formula without free variables is a sentence.

4.5.15. Definition of substitution

If **En**(X), a is a term and v is a variable, $X^{a}_{v}$ is the expression obtained by substituting all free occurrences of v in **X** with term a.

4.5.16. Definition of substitutability

Term a is substitutable for variable v in A just if A is atomic, or A is $\uparrow$ BC and a is substituable for v in both B and C, or A is $\forall$wB and v is not free in B, or, w does not occur in a and a is substitutable for v in B.

4.5.17. Definition of postfixed variable vector notation

$$\mathcal{E}(v, w, x) \text{ signifies that variables } v, w \text{ and } x \text{ are free in } \mathcal{E}.$$

4.5.18. Presentation resolve

$$A(v, w, a) \text{ may be written for } A(v, w, x)^{a}_{x}.$$

4.5.19. Definition of prefixed variable vector notation:

Occasionally $\forall\vec{\mathbb{v}}A$ is used for a sentence which either is A, or for some $n > 0$ and variables $\mathbb{v}_0 \ldots \mathbb{v}_{n-1}$, $\forall\vec{\mathbb{v}}A$ is $\forall\mathbb{v}_0 \ldots \forall\mathbb{v}_{n-1}A$.

4.5.20. Definition on parentheses, and defined operators for the object language:

(1) Delimiters for punctuation: (, ), [, ], …

(2) $\neg A == \downarrow AA$

(3) $(A \wedge B) == \downarrow \neg A \neg B$

(4) $(A \vee B) == \neg \downarrow AB$

(5) $(A \to B) == (\neg A \vee B)$

(6) $(A \leftrightarrow B) == (A \to B) \wedge (B \to A)$

(7) $\exists vA == \neg \forall v \neg A$

(8) $a \in b == ba$

(9) $\{v|A\} == \varsigma vA$

4.5.21. Definition on notation for binders restricted to set $b$:

(1) $A^b$ and $a^b$ signifiy that all variables bound in A and a are restricted to b.

(2) $v^b$ is v.

(3) $(c \in d)^b$ is $c^b \in d^b$.

(4) $\neg A^b$ is $\neg(A^b)$, $(A \wedge B)^b$ is $(A^b \wedge B^b)$, and so on for other connectives.

(5) $\{v|A\}^b \mathrel{=\!=\!=} \{v|v \in b \wedge A^b\}$.

(6) $(\forall v)A^b \mathrel{=\!=\!=} (\forall v)(v \in b \to A^b)$.

(7) $(\forall \vec{v})A^b \mathrel{=\!=\!=}$ the *sentence* given by the least $n \geq 0$ such that

$$\big(n > 0 \;\&\; (\forall v_0 \ldots \forall v_{n-1})(v_0 \in b \wedge \ldots \wedge v_{n-1} \in b \to A^b)\big)$$

ⓡ

$$\big(n = 0 \;\&\; A^b\big).$$

## 4.6 The expression codes and their expressions

We define a coding scheme which is such that the internal set theory **£** considers the code of an expression to be one of the sets in an omega ordered progression of sets, while the metatheory considers the same code to be a metanatural number.

4.6.1. Definition of expression codes

Presuppose bijective base-8 numerals:

(1) $\ulcorner n \urcorner \mathrel{=\!=\!=} \#n$

(2) $\ulcorner 1 \urcorner \mathrel{=\!=\!=}$ the symbolical *austerity* $\varsigma\, \ddot{v}\, \forall\, \ddot{v} \bullet {\downarrow}{\downarrow}{\downarrow}\, \ddot{v} \bullet \ddot{v}\ddot{v} \bullet \ddot{v}\ddot{v} \bullet \ddot{v} \downarrow \ddot{v} \bullet \ddot{v}\ddot{v} \bullet \ddot{v}\ddot{v} \bullet \ddot{v}$ of the set term $\{x|\forall y(x \in y \to x \in y)\}$, and so by Definitions 4.5.3 and 4.5.4 numerically identical to the bijective base-8 number 524213332122122123212212212.

(3) $\ulcorner n+1 \urcorner \mathrel{=\!=\!=}$ is the austerity $\varsigma\, \ddot{v}\, \forall\, \ddot{v} \bullet {\downarrow}{\downarrow}{\downarrow}\, \ddot{v} \bullet \ddot{v}\ddot{v} \bullet \ddot{v}\ddot{v} \bullet \ddot{v} \downarrow \ulcorner n \urcorner \bullet \ddot{v}\ddot{v} \bullet \ddot{v}\ddot{v} \bullet \ulcorner n \urcorner$ of $\{x|\forall y(x \in y \to \ulcorner n \urcorner \in y)\}$.

4.6.2. Theorem on expression codes

(1) If formal expression $\mathcal{E} \varepsilon \Omega_+ = \llbracket 1, 2, \ldots \rrbracket$, $\ulcorner \mathcal{E} \urcorner \varepsilon \Omega_+$ and $\ulcorner \mathcal{E} \urcorner$ denotes a set of **£** as per Remarks on expression names 4.6.2.3.

(2) The domain of $\ulcorner\ \urcorner$ is the metaset $\Omega_+ = \llbracket 1, 2, 3, \ldots \rrbracket$ of counting metanumbers, and its codomain is a metasubset of $\Omega_+$.

(3) As it turns out, on account of §9, $\ulcorner n + 1 \urcorner$ denotes $\{x|x = \ulcorner n \urcorner\}$. Given Definition 18.0.5 the denotation of $\ulcorner 1 \urcorner$ is the universal set $U =\!\!= \{x|x = x\}$. So the denotata of the members of the codomain of $\ulcorner\ \urcorner$ are given by the series $U, \{U\}, \{\{U\}\}, \ldots$, whose members one mnemonically may take to be *universe counters*. Write $U^1$ for U and $U^{n+1}$ for $\{U^n\}$.

(4) So if $\mathcal{E}$ is a formal expression, $\ulcorner \mathcal{E} \urcorner$ is $U^{\mathcal{E}}$.

(5) Diagonalization in $\mathfrak{L}$ is not as in *Peano arithmetic*, nor as in the associated modal provability logic GL with precisely $\Box(\Box p \rightarrow p) \rightarrow \Box p)$ as characteristic axiom. So the rounded symbolization $\ulcorner A \urcorner$ of Definition 4.6.1, distinct from usual Gödel codes as in $\ulcorner A \urcorner$, is used for the name of formula A.

(6) The author thinks that it is useful to consider the truth set of Definition 4.5.3.7 to be a *predicate*, to conform with the literature on theories of truth, and the apposition form $\mathcal{T} \ulcorner A \urcorner$ is therefore preferred to $\ulcorner A \urcorner \in \mathcal{T}$ in the exposition. *Caveat*: despite such an interpretation of $\mathcal{T}$ as a predicate, $\mathcal{T}$ is a set.

(7) The diagonalization construction allows quantifying into named contexts without a detour into substitution, and so $\forall x \exists y \mathcal{T} \ulcorner x = y \urcorner$ is e.g. a well formed sentence.

# 5 Semantics

> “*Development of mathematics resembles a fast revolution of a wheel: sprinkles of water are flying in all directions. Fashion – it is the stream that leaves the main trajectory in the tangential direction. The streams of epigone works draw most attention, and they constitute the main mass, but inevitably disappear after a while because they parted with the wheel. To remain on the wheel, it is necessary to apply the effort in the direction perpendicular to the main stream.*”
>
> ---
>
> Vladimir Igorevich Arnold, as quoted in (Khesin and Tabachnikov 2012, p. 379).

The first subsection gives an outline of the metasettheory which will be presupposed, and indicates its strength at the level of ramified analysis. The succeeding sections are geared to obtain a Herzberger style semantic with a semantic process on ordinal numbers, starting at 0 with an initial real number, i.e. a set of natural numbers - which via the *metanumber standpoint* is interpreted also as a set of formulas; the semantic process ends with a different real number altogether at other real numbers and in particular at the closure ordinal. The semantic behaviour at the closure ordinal is instrumental in isolating satisfaction notions which are instrumental in justifying the librationist set theory which is developed in the following sections.

## 5.1 The metasettheory and its metagenus of metaclasses

To define the underlying theory of the metalanguage further, we use the term *metaset* instead of the word *set*, and the theory to be presupposed is *metasettheory* $ZF^{-}$, which is Zermelo-Fraenkel metasettheory without power set. In addition, constructible *metaclasses* are presupposed for some considerations. $L_{\beta_0}$ is the level at the constructible hierarchy which we need for the model of metasettheory $ZF^{-}$, where $\beta_0$ is the ordinal for ramified analysis, and the least ordinal such that L at that point is a model of *analysis*, i.e. second order arithmetic $\mathbf{Z}_2$, and of $ZF^{-}$. It is presupposed that $\overset{\text{meta}}{\mathrm{V}} = \overset{\text{meta}}{\mathrm{L}}{}_{\beta_0}$.

*Metagenus* $\mathcal{G} = \{\!\!\{ x \overset{\text{meta}}{\subset} \overset{\text{meta}}{\mathrm{L}} \vdots\, x \varepsilon L_{\beta_0 \cdot 2} \}\!\!\}$ provides any such *metaclasses* as we shall want in the meta theoretic rendering of librationism. Incidentally, that is an enormous overshoot if one only wants a metaclass in $\mathcal{G}$ which surjects from $\omega \in L_{\omega+1}$ to $V = L_{\beta_0}$, as there is such a surjection class already in the much smaller metagenus $\{\!\!\{ x \subset V \vdots\, x \varepsilon L_{\beta_0+2} \}\!\!\}$. Confer item 5.2.1.12 in the Barter postulate.

Metaset-theory $\Sigma_3 KP\Omega$, which is Kripke-Platek metaset-theory, with $\Sigma_3$-collection and with the variety $\Omega$ of natural metanumbers – as introduced in §4.2, is sufficient for the isolation of the closure ordinal. The strength of ŧƀ, and other extensions of **£**, is increased beyond the first $\Sigma_3-$admissible ordinal, however, and up to the ordinal $\beta_0$ for ramified analysis, by the content coded into the barter $\VDash$ of §4.6.2.

NOTE: Move this to below the section 5.2 on the barter. Or rewrite.
this is to allow the wide variety of the barter $\VDash$ on the metanumbers in $\Omega$, which e.g. justifies the reasonable result that neither $\{x|x \in x\} \in \{x|x \in x\}$ nor $\{x|x \in x\} \notin \{x|x \in x\}$ are librationistically valid.

Care should be taken to not confuse the metasets of the meta language used to introduce **£** with the sets **£**, or its extensions, postulate the existence of.

## 5.2 Barters

5.2.1. The Barter Postulate

(1) The *barter*, $\Vdash$, is a *semantic mapping* from $\Omega \times \mathrm{Ord}$ to its codomain ↀ.

(2) For $n \varepsilon \Omega$ and $\mathrm{Ord}(\alpha)$, $\Vdash^{\alpha}_{n}$ is written instead of $\Vdash(n, \alpha)$.

(3) For any metanumber $n$ and ordinal $\alpha$, $\Vdash^{\alpha}_{n}$ is a metaset of metasets of formulas; as formulas are identified with metanatural numbers, $\Vdash^{\alpha}_{n}$ can be considered to be a metaset of metareal numbers - for metasets of metanatural numbers are as much metareal numbers as sets of natural numbers are real numbers.

(4) The induced ordinal functions $\Vdash^{\alpha}_{0}, \Vdash^{\alpha}_{1}, \Vdash^{\alpha}_{2}, \dots$ to sets of formulas, which satisfy the conditions in this postulate, are the barters at $\alpha$ on 0, at $\alpha$ on 1, and so on.

(5) If the ordinal is fixed, at 1, say, one may consider $\Vdash^{1}_{x}$ a function with domain $\Omega$.

(6) For any barter at $\alpha$ on $n$, $\Vdash^{\alpha}_{n}$, and formula A, $\Vdash^{\alpha}_{n}$ A is written instead of $A \varepsilon \Vdash^{\alpha}_{n}$.

(7) The *constructible* real numbers in $\mathbb{R}_{L_{\beta_0}}$ are the metasets required in 5.2.1.3.

(8) Real number x is *stroked* $:= \downarrow AB \varepsilon x$ just if neither $A \varepsilon x$ nor $B \varepsilon x$.

(9) Real number x is *exhaustive* $:= \forall v A \varepsilon x$ just if $A^{b}_{v} \varepsilon x$ for all b substitutable for v in A.

(10) Real number x is *full* just if x is *stroked* & *exhaustive*.

(11) ↀ $= \{\!\![x \vdots\, x \varepsilon \mathrm{Def}(\mathbb{R}_{L_{\beta_0}})$ & *full*$(x)]\!\!\}$.

(12) $\Vdash^{0}_{x}$ is presupposed to be a *given* surjection from $\Omega$ to ↀ.[1]

(13) For $n$, and $\alpha = 0$,

$$\Pi A(\Vdash^{0}_{n} \mathcal{T}\ulcorner A \urcorner \Leftrightarrow \mathcal{T}\ulcorner A \urcorner \varepsilon \Vdash^{0}_{n}).$$

(14) For $n$, and $\alpha > 0$, $\Vdash^{\alpha}_{n}$ is the L−least member of ↀ such that

$$\Pi A(\Vdash^{\alpha}_{n} \mathcal{T}\ulcorner A \urcorner \Leftrightarrow \Sigma\gamma(\gamma \prec \alpha \,\&\, \Pi\delta(\gamma \preceq \delta \prec \alpha \Rightarrow \Vdash^{\delta}_{n} A))),$$

*and*

$$\Pi A(\Vdash^{\alpha}_{n} \mathcal{T}\ulcorner A \urcorner \Leftrightarrow \mathcal{T}\ulcorner A \urcorner \varepsilon \Vdash^{\alpha}_{n}).$$

[1] There is such a surjection as $\Vdash^{0}_{x}$ in $L_{\beta_0+2}$, as an answer by Gabe Goldberg to a question I posed on MathOverflow on July 18, 2024, effectively defines a surjection from $\Omega$ to $\mathbb{R}_{L_{\beta_0}}$ which is in $L_{\beta_0+2}$ - the next step to obtain a surjection from $\Omega$ to ↀ only appeals to the subset axiom, which holds at $L_{\beta_0+2}$; it can be shown that there is no such surjection from $\Omega$ to $\mathbb{R}_{L_{\beta_0}}$ in a $L_{\gamma}$ where $\gamma \prec \beta_0 + 2$.

(15) Notice again that $\Vdash_{x}^{0}$ is onto ↀ, by presupposition. But $\Vdash_{x}^{\alpha}$ for ordinal $\alpha \succ 0$ is not onto ↀ, for there is no $n \varepsilon \Omega$ such that $\Vdash_{n}^{\alpha}$ has an element which has no sentence of the form $\mathcal{T}\ulcorner A \urcorner$ as a member.

(16) $\Vdash^{\alpha} A := \Pi n(n \varepsilon \Omega \Rightarrow \Vdash_{n}^{\alpha} A)$ is *the generalbarter* at $\alpha$.

(17) Given 5.2.1.14 and 5.2.1.16 and $\alpha \succ 0$,

$$\Vdash^{\alpha} \mathcal{T}\ulcorner A \urcorner \Leftrightarrow \boldsymbol{\Sigma}\gamma(\gamma \prec \alpha \ \& \ \boldsymbol{\Pi}\delta(\gamma \preceq \delta \prec \alpha \Rightarrow \Vdash^{\delta} A)).$$

(18) The barter $\Vdash$ on $n$, where $n \varepsilon \Omega$, i.e. $\Vdash_{n}$, is *the* $n$*-barter*.

5.2.2. Theorem (**£** is existentially well behaved)

$$\Vdash^{\alpha} \exists v A \Leftrightarrow \Vdash^{\alpha} A_{v}^{b} \text{ for some b substitutable for v in A.}$$

*Proof:*

$$\Vdash^{\alpha} \exists v A$$

$$\Updownarrow$$

$$\Vdash^{\alpha} \neg \forall v \neg A$$

$$5.2.1.8 \Updownarrow$$

$$\nVdash^{\alpha} \forall v \neg A$$

$$5.2.1.9 \Updownarrow$$

For a b substitutable for v in A, $\nVdash^{\alpha} \neg A_{v}^{b}$.

$$5.2.1.8 \Updownarrow$$

For a b substitutable for v in A, $\Vdash^{\alpha} A_{v}^{b}$.

□

5.2.3. Presentation resolve

When the context allows, write $\mathcal{T}A$ as an abbreviation for $\mathcal{T}\ulcorner A \urcorner$.

5.2.4. Theorem (The barters of metanumbers for logic)

For any metanumber $n$ and any ordinal metanumber $\alpha$:

(1) $(\Vdash^{\alpha}_{n} \downarrow AB \Leftrightarrow$ neither $\Vdash^{\alpha}_{n} A$ nor $\Vdash^{\alpha}_{n} B)$.

(2) $(\Vdash^{\alpha}_{n} \forall xA \Leftrightarrow \boldsymbol{\Pi} b($b substitutable for v in A: $\Vdash^{\alpha}_{n} A^{b}_{x}))$.

*Proof.* Exercise. □

5.2.5. Theorem (The barters of metanumbers for truth)

For any metanumber $n$, and any ordinal metanumber $\alpha > 0$:

$$\Vdash^{\alpha}_{n} \mathcal{T}A \Leftrightarrow \boldsymbol{\Sigma}\gamma(\gamma \prec \alpha \ \&\ \boldsymbol{\Pi}\delta(\gamma \preceq \delta \prec \alpha \Rightarrow \Vdash^{\delta}_{n} A)).$$

*Proof.* Exercise. □

5.2.6. Theorem (The barters of metanumbers for sets)

For any metanumber $n$, and any ordinal metanumber $\alpha > 0$:

$$\Vdash^{\alpha}_{n} (a \in \{x|A\} \leftrightarrow \mathcal{T}\ulcorner A^{a}_{x} \urcorner), \text{ where a is substitutable for x in A.}$$

*Proof.* Exercise. □

## 5.3 Alethic comprehension

For any metanumber $n$, and any ordinal metanumber $\alpha > 0$,

5.3.1. Theorem on alethic comprehension with parameters from $\vec{v}$:

$$\Vdash^{\alpha}_{n} \forall\vec{v}\forall x(x \in \{y|A\} \leftrightarrow \mathcal{T}\ulcorner A \urcorner^{x}_{y}), \text{ where x is substitutable for y in A.}$$

*Proof.* Appeal to Theorem 5.2.6. □

## 5.4 Closure

5.4.1. Definition of cover, stabilization and closure

(1) $IN^{\alpha}(A)$ just if $\boldsymbol{\Pi}\beta(\alpha \preceq \beta \Rightarrow \Vdash^{\beta} \mathcal{T}A)$.

(2) $OUT^{\alpha}(A)$ just if $\boldsymbol{\Pi}\beta(\alpha \preceq \beta \Rightarrow \nVdash^{\beta} \mathcal{T}A)$.

(3) $IN(A)$ just if $\boldsymbol{\Sigma}\alpha IN^{\alpha}(A)$.

(4) $OUT(A)$ just if $\boldsymbol{\Sigma}\alpha OUT^{\alpha}(A)$.

(5) STAB(A) just if IN(A) Ⓡ OUT(A).

(6) UNSTAB(A) just if $\not{\text{STAB}}$(A).

(7) $\alpha$ *covers* A just if IN(A) $\Rightarrow$ $\text{IN}^{\alpha}$(A).

(8) $\alpha$ *stabilizes* A just if $\alpha$ covers A, and $\Vdash^{\alpha} \mathcal{T}\text{A} \Rightarrow \text{IN(A)}$.

(9) The *closure* ordinal is the least stabilizing ordinal.

5.4.2. Theorem (Herzberger 1980)

There is a closure ordinal.

*Proof:* Assume first that IN(A), to presuppose

5.4.3. Definition $$\text{H(A)} = \mu\alpha(\text{IN}^{\alpha}(\text{A})).$$

1. We show that there is a covering ordinal:

   We have

$$\mathbf{\Pi}\text{A}(\text{IN(A)} \Rightarrow \mathbf{\Sigma}\beta(\beta = \text{H(A)})). \tag{5.4.4}$$

   So

$$\mathbf{\Pi}\text{A}\mathbf{\Sigma}\beta(\text{IN(A)} \Rightarrow \beta = \text{H(A)}). \tag{5.4.5}$$

   $\Pi_2$–collection and quantifier rules give us

$$\mathbf{\Pi}\text{B}\mathbf{\Sigma}\text{Y}\mathbf{\Pi}\text{A}(\text{A}\varepsilon\text{B} \Rightarrow \mathbf{\Sigma}\beta(\beta\varepsilon\text{Y}\&(\beta = \text{H(A)})))). \tag{5.4.6}$$

   Instantiate with B $=$ [A : IN(A)] to obtain

$$\mathbf{\Sigma}\text{Y}\mathbf{\Pi}\text{A}(\text{IN(A)} \Rightarrow \mathbf{\Sigma}\beta(\beta\varepsilon\text{Y}\&(\beta = \text{H(A)})). \tag{5.4.7}$$

   I

   Let Z be a witness for (5.4.7), and define the least covering ordinal by means of $\Pi_2$–separation,

$$\varkappa = [\nu : \nu\varepsilon\text{Z}\ \&\ \text{Ordinal}(\nu)\ \&\ \mathbf{\Sigma}\text{A}(\text{IN(A)}\ \&\ \nu = \text{H(A)})]. \tag{5.4.8}$$

$\Pi_2$-collection was invoked in the step from (5.4.5) to (5.4.6), and as $\Pi_n$-collection implies $\mathbf{\Sigma}_{n+1}$ collection in Kripke–Platek theories, this justifies the choice of an underlying variety theory at least as strong as $\Sigma_3 KP\Omega$ for the meta language.[2]

2. We prove that there is a stabilizing ordinal:

   Let $[f(n): n\varepsilon\Omega]$, by an adaptation of Cantor's pairing function, be an enumeration of all and only elements of UNSTAB, where each element recurs infinitely often, so that if B=$f(m)$ and $m \prec n\varepsilon\Omega$, then there is a natural number o, $n \prec o\varepsilon\Omega$, such that $f(o) = B$. Let $g(0) = \varkappa$ of Equation 5.4.8, and $g(n+1) =$ the least $\nu > g(n)$ such that

   $$\Vdash^{\nu} f(n) \Leftrightarrow \nVdash^{g(n)} f(n).$$

   Let $\text{ß} = [\gamma\colon \mathbf{\Sigma} m \mathbf{\Sigma} \nu (m\varepsilon\Omega \,\&\, \nu = g(m) \,\&\, \gamma\varepsilon\nu)]$. It is obvious that ß is a limit ordinal which covers all elements of UNSTAB. It is also clear that if $m \prec n\varepsilon\Omega$ then $g(m) \prec g(n)$. Since ß covers UNSTAB, it suffices to show that $\Vdash^{\text{ß}} \mathcal{T}B$ entails that B is in STAB, to establish that ß is a stabilization ordinal.

   Suppose $\Vdash^{\text{ß}} \mathcal{T}B$. It follows that

   a) $\mathbf{\Sigma}\nu\mathbf{\Pi}\xi(\nu \preceq \xi \prec \text{ß} \Rightarrow \Vdash^{\xi} B)$

   Since g is increasing with ß as its range, we will for some natural number $m\varepsilon\Omega$ have that $\nu \preceq g(m) \prec \text{ß}$, so that

   b) $\mathbf{\Pi}\xi(g(m) \preceq \xi \prec \text{ß} \Rightarrow \Vdash^{\xi} B)$

   Suppose $B \notin STAB(\Xi)$. By our enumeration of unstable elements where each term recurs infinitely often, we have that $B = f(n)$ for some natural number n, $m \prec n \in \Omega$. It follows that $g(m) \prec g(n) \prec \text{ß}$. From a) and b) we can infer that $\Vdash^{g(n)} B$, since we supposed that $\Vdash^{\text{ß}} \mathcal{T}B$. From the construction of function g, $\nVdash^{g(n+1)} \neg B$, contradicting b). It follows that $\Vdash^{\text{ß}} \mathcal{T}B$ only if $B \in STAB(\Xi)$, so ß is a stabilization ordinal.

3. The proof that there is a closure ordinal finishes with an appeal to Definition 5.4.1.9. □

5.4.9. Definition

$\varpi$ is the closure ordinal.

5.4.10. Theorem

The closure ordinal $\varpi$ is the least $\Sigma_3$-admissible ordinal.

[2](Welch 2011) shows that KP + $\Sigma_3$-Determinacy is sufficient for the semantics of a commensurate system AQI (*Arithmetical Quasi Induction*) introduced in (Burgess 1986), and (Hachtman 2019) shows this equivalent to KP + $\Pi^1_2$-Monotone Induction. So a $\Sigma_3$-admissible ordinal is not necessary, but it may be needed for the proof we use, which connects the coding of the formal language with the natural metanumbers of the meta theory. Welch has pointed out in private communication that a $\Sigma_2$-admissible ordinal, without further assumptions, can be proven to be insufficient.

*Proof.* Combine Definition 5.4.9, and the result, of part 1 of the proof of Theorem 5.4.2, that a variety theory at least as strong as $\Sigma_3\mathrm{KP}\Omega$ is needed for the meta language. □

## 5.5 Fairs and satisfaction for theses, and their appropriateness

5.5.1. Definition (Satisfactions in terms of fairs)

(1) Barter $\VDash^{\xi}_{n}$ is a *fair* just if $\xi$ is the closure ordinal $\varpi$.

(2) Note e.g. for $s = \{x|x \in x\}$ that some fairs have $\VDash^{\varpi}_{m}$ $s \in s$ and others have $\VDash^{\varpi}_{n}$ $s \notin s$.

(3) A is *maximally* satisfied just if $\mathcal{T}\ulcorner A \urcorner$ is an element in all fairs.

(4) A is *optimally* satisfied just if A is an element in all fairs.

(5) A is satisfied just if $\neg\mathcal{T}\ulcorner \neg A \urcorner$ is an element in all fairs.

(6) A is *paradoxically* satisfied just if both A and $\neg$A are optimally satisfied.

5.5.2. Definition (Satisfaction notation in terms of fairs)

(1) $\Vdash^{M} A$ — A is a *maxim* just if $\mathbf{\Pi} n(n \varepsilon \Omega \Rightarrow \mathcal{T}\ulcorner A \urcorner \varepsilon \VDash^{\varpi}_{n})$.

(2) $\Vdash^{O} A$ — A is an *optimum* just if $\mathbf{\Pi} n(n \varepsilon \Omega \Rightarrow A \varepsilon \VDash^{\varpi}_{n})$.

(3) $\Vdash A$ — A is a *thesis* just if $\mathbf{\Pi} n(n \varepsilon \Omega \Rightarrow \neg\mathcal{T}\ulcorner \neg A \urcorner \varepsilon \VDash^{\varpi}_{n})$.

(4) $\Vdash^{P} A$ — A is a *paradox* just if $\Vdash A$ and $\Vdash \neg A$.

5.5.3. Observation: Notice, with reference to Definition 5.5.2, that, given the results of §5.4,

$$\Vdash^{M} A \Rightarrow \Vdash^{O} A \ \& \ \Vdash^{O} A \Rightarrow \Vdash A.$$

5.5.4. Definition of *appropriate* schematic maxim, optimum and thesis designations:

(1) Let A, B, C ... be sentence schemas, so that A e.g. may be given by $p \rightarrow (q \rightarrow p)$. In cases as the latter one may say that $\Vdash^{M}$ A because $\Vdash^{M} p \rightarrow (q \rightarrow p)$ for all *substitutions*.

(2) $\Vdash^{M}$ A is an appropriate maxim designation just if for all substitutions $\Vdash^{M}$ A.

(3) $\Vdash^{O}$ B is an appropriate optimum designation just if $\Vdash^{M}$ B is not an appropriate maxim designation and for all substitutions $\Vdash^{O}$ A.

(4) $\Vdash$ C is an appropriate thesis designation precisely if $\Vdash^{O}$ C is not an appropriate optimum designation and $\Vdash^{M}$ C is not an appropriate maxim designation, but for all subsitutions $\Vdash$ C.

5.5.5. Resolve of appropriateness:

The maxims, optima and theses in the text are as a rule supposed to be appropriately stated.

## 5.6 If A is a thesis, then A is a paradox or a maxim

5.6.1. The Grammatical Conjecture

Sentences cannot be paradoxical according to some but not all barters.

5.6.2. The Bifurcation Theorem

For $\alpha \preceq \varpi$:

$$\mathbf{\Pi}\mathrm{n}\Big(\mathrm{n}\varepsilon\Omega \Rightarrow \mathcal{T}\ulcorner\mathrm{A}\urcorner\,\varepsilon \VDash^{\alpha}_{\mathrm{n}}\Big)$$

$$\Updownarrow$$

$$\mathbf{\Pi}\mathrm{n}\Big(\mathrm{n}\varepsilon\Omega \Rightarrow \neg\mathcal{T}\ulcorner\neg\mathrm{A}\urcorner\,\varepsilon \VDash^{\alpha}_{\mathrm{n}}\Big) \ \& \ \mathbf{\Sigma}\mathrm{n}\Big(\mathrm{n}\varepsilon\omega \;\&\; \mathcal{T}\ulcorner\mathrm{A}\urcorner\,\varepsilon \VDash^{\alpha}_{\mathrm{n}}\Big)$$

*Proof.* The Updownarrow is composed of a Downarrow and an Uparrow, and they are considered separately in the ordinal cases.

$\alpha = 0$: The Downarrow holds on account of Postulate 6.3.2.2 and the fact that barters exist given the Barter Postulate, 5.2.1. The Grammatical Conjecture 5.6.1 justifies the Uparrow.

$\alpha = \beta + 1$ for some ordinal $\beta$: The Downarrow holds on account of the validity of Postulate 6.3.2.2, and the fact that there are barters, given Barter Postulate, 5.2.1. The Uparrow holds as $\alpha$ is a successor ordinal, and so on account of Postulates 6.3.2.3, 6.3.2.6 and 6.3.2.7, validates the step from $\mathbf{\Pi}\mathrm{n}\Big(\mathrm{n}\varepsilon\Omega \Rightarrow \neg\mathcal{T}\ulcorner\neg\mathrm{A}\urcorner\,\varepsilon \VDash^{\alpha}_{\mathrm{n}}\Big)$ to $\mathbf{\Pi}\mathrm{n}\Big(\mathrm{n}\varepsilon\Omega \Rightarrow \mathcal{T}\ulcorner\mathrm{A}\urcorner\,\varepsilon \VDash^{\alpha}_{\mathrm{n}}\Big)$, and as the conjunct $\mathbf{\Sigma}\mathrm{n}\Big(\mathrm{n}\varepsilon\omega \;\&\; \mathcal{T}\ulcorner\mathrm{A}\urcorner\,\varepsilon \VDash^{\alpha}_{\mathrm{n}}\Big)$ supports the joint entailment needed for the Uparrow.

Let finally $\alpha = \lambda$ for a limit ordinal $\lambda \preceq \varpi$: The Downarrow again holds on account of Postulate 6.3.2.2 and the the fact that there are barters, given Barter Postulate, 5.2.1. The Grammatical Conjecture 5.6.1 justifies the Uparrow. □

5.6.3. Theorem

$$\Vdash^{\mathrm{M}} \mathrm{A} \ \Leftrightarrow \Vdash \mathrm{A} \ \& \ \nVdash \neg\mathrm{A}.$$

*Proof.* If not, on account of the Bifurcation Theorem 5.6.2, either

$$\mathbf{\Pi}\mathrm{n}\Big(\mathrm{n}\varepsilon\Omega \Rightarrow \mathcal{T}\ulcorner\mathrm{A}\urcorner\,\varepsilon \VDash^{\alpha}_{\mathrm{n}}\Big) \ \& \ (\mathbf{\Sigma}\mathrm{n}\Big(\mathrm{n}\varepsilon\omega \;\&\; \mathcal{T}\ulcorner\neg\mathrm{A}\urcorner\,\varepsilon \VDash^{\alpha}_{\mathrm{n}}\Big) \ \circledR \ \mathbf{\Pi}\mathrm{n}\Big(\mathrm{n}\varepsilon\Omega \Rightarrow \neg\mathcal{T}\ulcorner\mathrm{A}\urcorner\,\varepsilon \VDash^{\alpha}_{\mathrm{n}}\Big))$$

or

$$\mathbf{\Pi}\mathrm{n}\Big(\mathrm{n}\varepsilon\Omega \Rightarrow \neg\mathcal{T}\ulcorner\neg\mathrm{A}\urcorner\,\varepsilon \VDash^{\alpha}_{\mathrm{n}}\Big) \ \& \ \mathbf{\Sigma}\mathrm{n}\Big(\mathrm{n}\varepsilon\omega \;\&\; \mathcal{T}\ulcorner\mathrm{A}\urcorner\,\varepsilon \VDash^{\alpha}_{\mathrm{n}}\Big) \ \& \ \mathbf{\Sigma}\mathrm{n}\Big(\mathrm{n}\varepsilon\omega \;\&\; \neg\mathcal{T}\ulcorner\mathrm{A}\urcorner\,\varepsilon \VDash^{\alpha}_{\mathrm{n}}\Big).$$

But these are impossible, given Grammatical Conjecture 5.6.1 and Truth maxim 6.3.2.2. □

5.6.4. Corollary on optimal arbitration:

$$\Vdash^{\mathrm{O}} (\mathcal{T}\ulcorner\mathrm{A}\urcorner \vee \mathcal{T}\ulcorner\neg\mathrm{A}\urcorner) \ \circledR \ (\Vdash \mathrm{A} \Leftrightarrow \Vdash \neg\mathrm{A}).$$

*Proof.* If the right disjunct is false so that ($\Vdash \mathrm{A} \nLeftrightarrow \Vdash \neg\mathrm{A}$), either $\Vdash \mathrm{A} \ \& \ \nVdash \neg\mathrm{A}$ or else $\Vdash \neg\mathrm{A} \ \& \ \nVdash \mathrm{A}$, and so on account of Theorem 5.6.3 $\VDash^{\mathrm{M}} \mathrm{A}$ or $\VDash^{\mathrm{M}} \neg\mathrm{A}$ – and either of the latter entail the left side of the corollary. So that proves the disjunction. □

5.6.5. Theorem on maximal arbitration:

$$\Vdash^{M} (\mathcal{T}\ulcorner A\urcorner \vee \mathcal{T}\ulcorner\neg A\urcorner)\ \circledR\ (\Vdash A \Leftrightarrow \Vdash \neg A).$$

*Proof.* Given Corollary 5.6.4 and Postulate 6.5.1.1 it suffices to prove

$$\Vdash^{O} (\mathcal{T}A \vee \mathcal{T}\neg A) \Rightarrow \Vdash^{M} (\mathcal{T}A \vee \mathcal{T}\neg A).$$

Use a disjunctive syllogism to obtain $\Vdash^{O} (\mathcal{T}\ulcorner\mathcal{T}\ulcorner A\urcorner\urcorner \vee \mathcal{T}\ulcorner\mathcal{T}\ulcorner\neg A\urcorner\urcorner)$ with Postulate 6.5.1.2 from $\Vdash^{O} (\mathcal{T}\ulcorner A\urcorner \vee \mathcal{T}\ulcorner\neg A\urcorner)$. Use theorem $\Vdash^{O} \mathcal{T}\ulcorner\mathcal{T}\ulcorner B\urcorner\urcorner \rightarrow \mathcal{T}\ulcorner(\mathcal{T}\ulcorner B\urcorner \vee C)\urcorner$ and disjunctive syllogism with $\Vdash^{O} (\mathcal{T}\mathcal{T}\ulcorner A\urcorner \vee \mathcal{T}\mathcal{T}\ulcorner\neg A\urcorner)$ to obtain $\Vdash^{O} \mathcal{T}\ulcorner\mathcal{T}\ulcorner A\urcorner \vee \mathcal{T}\ulcorner\neg A\urcorner\urcorner$. Use Definitions 5.5.1.3, 5.5.1.4, 5.5.1.3 and 5.5.1.4 and results noted to conclude $\Vdash^{M} \mathcal{T}\ulcorner A\urcorner \vee \mathcal{T}\ulcorner\neg A\urcorner$. □

5.6.6. Theorem (Maxims are optimal)

$\VDash^{M} A \Rightarrow \VDash^{O} A$

*Proof.* By Definition 5.5.2, $\VDash^{M}$ A only if for any metanumber $\mathfrak{n}$, $\mathcal{T}A\varepsilon \VDash^{\varpi}_{\mathfrak{n}}$ – so that $\VDash^{O} \mathcal{T}A$. Given Postulate 6.5.1.2, $\VDash^{O}$ A. □

## 5.7 Orthodoxy, definiteness and paradoxicality

5.7.1. Definition

(1) A is *orthodox* just if $\Vdash^{O} \forall\vec{v}(\mathcal{T}A \vee \mathcal{T}\neg A)$.

(2) Set a is orthodox just if $x \in a$ is orthodox.

(3) A is *definite* just if $\Vdash A$ or $\Vdash \neg A$.

(4) A is *apocryphal* just if orthodox and *indefinite*.

(5) Set a is apocryphal just if $b \in a$ is apocryphal for some set b.

5.7.2. Remark

Some definite sentences are determinate, and some are indeterminate.

5.7.3. Remark

Set $s \mathrel{=\!=} \{v|v \in v\}$ is apocryphal. For sentence $s \in s$ is apocryphal, by cause of its orthodoxy and the fact that it is indefinite as $\nVdash s \in s$ and $\nVdash s \notin s$.

5.7.4. Definition of paradoxical formula

Formula A is paradoxical just if it is not orthodox.

5.7.5. Remark

Given Definitions 5.7.1.1 and 5.7.4, A is paradoxical just if $\nVdash^{\mathrm{O}}\ \forall\vec{\mathrm{v}}(\mathcal{T}\mathrm{A}(\vec{\mathrm{v}}) \vee \mathcal{T}\neg\mathrm{A}(\vec{\mathrm{v}}))$; so there is, given Definition 5.5.1.4, a generalbarter $\VDash$ such that $\VDash^{\varpi}\ \exists\vec{\mathrm{v}}(\neg\mathcal{T}\mathrm{A}(\vec{\mathrm{v}}) \wedge \neg\mathcal{T}\neg\mathrm{A}(\vec{\mathrm{v}}))$. On account of Theorem 12.1.1, there is a, possibly empty, vector $\vec{\mathrm{a}}$ such that

$$\VDash^{\varpi}\ \neg\mathcal{T}\mathrm{A}(\vec{\mathrm{a}}) \wedge \neg\mathcal{T}\neg\mathrm{A}(\vec{\mathrm{a}}).$$

5.7.6. Definition

Set a is paradoxical just if it not orthodox just if $\nVdash^{\mathrm{O}}\ \forall\mathrm{x}(\mathcal{T}\mathrm{x} \in \mathrm{a} \vee \mathcal{T}\mathrm{x} \notin \mathrm{a})$; so, given Definition 5.5.1.4, there is a generalbarter $\VDash$ such that $\VDash^{\varpi}\ \exists\mathrm{x}(\neg\mathcal{T}\mathrm{x} \in \mathrm{a} \wedge \neg\mathcal{T}\mathrm{x} \notin \mathrm{a})$. Consequently, given Theorem 12.1.1, for some term b, $\VDash^{\varpi}\ (\neg\mathcal{T}\mathrm{b} \in \mathrm{a} \wedge \neg\mathcal{T}\mathrm{b} \notin \mathrm{a})$.

5.7.7. Fact

The proof of Theorem 5.6.5 shows that $\Vdash^{\mathrm{O}}\ (\mathcal{T}\mathrm{A} \vee \mathcal{T}\neg\mathrm{A}) \Rightarrow\ \Vdash^{\mathrm{M}}\ (\mathcal{T}\mathrm{A} \vee \mathcal{T}\neg\mathrm{A})$, so, given Definition 5.7.1.1 and Postulate 6.5.1.1, A is orthodox just if $\Vdash^{\mathrm{M}}\ (\mathcal{T}\mathrm{A} \vee \mathcal{T}\neg\mathrm{A})$. But the latter should not be used for defining orthodoxy, as the induced revision of Definition 5.7.4 would not give the intended extension for the term 'paradoxical'.

## 5.8 The non-triviality assumptions

5.8.1. Definition

A logical theory is *trivial* just if all of its sentences are derivable.

The assumption that there *are* barter initials for variants of **£** amounts to assuming that the system under consideration is not trivial, and, consequently, consistent. It follows from (Bjørdal 2012), under the assumption that $\Sigma_3$KP is consistent, that the empty set is a barter initial for **£** simpliciter. So if $\Sigma_3$KP is consistent, then **£** is not trivial.

§ 25 shows that ꝥ$\mathcal{H}\mathcal{R}(\mathbf{G})$ has an account of **NBGC** + **TA** if **NBGC** + **TA** is consistent.

# 6 Maxims, paradoxes and optima

*"Pure mathematics is, in its way, the poetry of logical ideas."*

`Albert Einstein` in (Einstein 1935).

## 6.1 Axioms & warrants, theorems & proofs.

6.1.1. Definition

(1) A *warrant* of an axiom is a semantic demonstration of it from Definition 5.2.1.

(2) A *proof* of a theorem is a demonstration of it from axioms and other theorems.

## 6.2 Logic maxims

6.2.1. Postulate (Classical logic maxims)

(1) $\Vdash^{M} A \rightarrow (B \rightarrow A)$

(2) $\Vdash^{M} (A \rightarrow (B \rightarrow C)) \rightarrow ((A \rightarrow B) \rightarrow (A \rightarrow C))$

(3) $\Vdash^{M} (\neg B \rightarrow \neg A) \rightarrow (A \rightarrow B)$

(4) $\Vdash^{M} \forall x(A \rightarrow B) \rightarrow (\forall xA \rightarrow \forall xB)$

(5) $\Vdash^{M} A \rightarrow \forall vA$, provided v is not free in A

(6) $\Vdash^{M} \forall vA \rightarrow A^{b}_{v}$ , provided b is substitutable for v in A

(7) $\Vdash^{M} \Gamma$ is one of (6.2.1.1–6.2.1.7) only if $\Vdash^{M} \forall v\Gamma$ is one of (6.2.1.1–6.2.1.7).

6.2.2. Remark

The role of the maximal mode, which allows the deduction from $\Vdash^{M} (A \rightarrow B)$ and $\Vdash^{M} A$ to $\Vdash^{M} B$, is played by mode 7.2.7.

6.2.3. Remark

An induction, upon 6.2.1.7 and 7.2.7, proves *generalization* is a derived inference mode. Compare the proof of Theorem 45.4 of (Hunter 1971, pp. 174–175).

## 6.3 Maxims on truth

6.3.1. Definition Russell's paradoxical set

$$r \Longleftarrow \{x|x \notin x\}$$

Definition 6.3.1 is presupposed in the articulation of the truth maxim in Postulate 6.3.2.3.

6.3.2. Postulate (Truth maxims)

(1) $\Vdash^{M} \mathcal{T}(A \to B) \to (\mathcal{T}A \to \mathcal{T}B)$

(2) $\Vdash^{M} \mathcal{T}A \to \neg\mathcal{T}\neg A$

(3) $\Vdash^{M} (\mathcal{T}r \in r \vee \mathcal{T}r \notin r) \to (\mathcal{T}A \vee \mathcal{T}\neg A)$

(4) $\Vdash^{M} \mathcal{T}A \vee \mathcal{T}\neg A \vee (\mathcal{T}\neg\mathcal{T}\neg B \to \mathcal{T}B)$

(5) $\Vdash^{M} \mathcal{T}A \vee \mathcal{T}\neg A \vee (\mathcal{T}B \to \mathcal{T}\mathcal{T}B)$

(6) $\Vdash^{M} \mathcal{T}(\mathcal{T}A \to A) \to (\mathcal{T}A \vee \mathcal{T}\neg A)$

(7) $\Vdash^{M} \mathcal{T}(\mathcal{T}A \to \mathcal{T}\mathcal{T}A) \to (\mathcal{T}A \vee \mathcal{T}\neg A)$

(8) $\Vdash^{M} \exists v\mathcal{T}A \to \mathcal{T}\exists vA$.

(9) $\Vdash^{M} \mathcal{T}\forall vA \to \forall v\mathcal{T}A$

(10) $\Vdash^{M} \forall u(a \in u \to b \in u) \to (A^{a}_{v} \to A^{b}_{v})$, for a and b both substitutable for v in A.

(11) $\Vdash^{M} \mathfrak{D}(\{x|A\}) \to (\forall x\mathcal{T}A \to \mathcal{T}\forall xA)$

## 6.4 *Excursus 1*: Warrants of maxims on truth

*Warrant 6.3.2.1.* Suppose $\VDash^{\gamma} \mathcal{T}(A \to B)$ and $\VDash^{\gamma} \mathcal{T}A$. It follows that for some ordinal $\delta$ and any ordinal $\epsilon$ such that $\delta \preceq \epsilon \prec \gamma$, $\VDash^{\epsilon} (A \to B)$ and $\VDash^{\epsilon} A$. So, on account of Definition 5.2.1.8, $\VDash^{\epsilon} B$, and, consequently, $\VDash^{\gamma} \mathcal{T}B$. So for any $\gamma$, $\VDash^{\gamma} \mathcal{T}(A \to B) \to (\mathcal{T}A \to \mathcal{T}B)$.

$$\VDash^{\varpi} \mathcal{T}(\mathcal{T}(A \to B) \to (\mathcal{T}A \to \mathcal{T}B))$$

is a consequence of this, so $\Vdash^{M} \mathcal{T}(A \to B) \to (\mathcal{T}A \to \mathcal{T}B)$. ∎

*Warrant 6.3.2.2.* Assume $\nVDash^{\gamma}(\mathcal{T}A \to \neg\mathcal{T}\neg A)$. It follows that $\VDash^{\gamma} (\mathcal{T}A \wedge \mathcal{T}\neg A)$. As a consequence, $\VDash^{\gamma} \mathcal{T}A$ and $\VDash^{\gamma} \mathcal{T}\neg A$. It follows that for some ordinal $\delta$ and any ordinal $\epsilon$ such that $\delta \preceq \epsilon \prec \gamma, \VDash^{\epsilon} A$ and $\VDash^{\epsilon} \neg A$. But that is impossible, so $\VDash^{\gamma} (\mathcal{T}A \to \neg\mathcal{T}\neg A)$, and we are done. ∎

*Warrant 6.3.2.3.* The postulate's maxim somewhat extends (Bjørdal 2012). Let an ordinal $\delta$ be *monogamous* just if a successor ordinal, so for monogamous $\delta$, for any sentence B, $\VDash^{\delta} \mathcal{T}B$ just if $\VDash^{\delta} \neg\mathcal{T}\neg B$. $\Vdash^{M} (\mathcal{T}r \in r \vee \mathcal{T}r \notin r) \to (\mathcal{T}A \vee \mathcal{T}\neg A)$ holds simply because monogamous ordinals are monogamous ordinals. ∎

*Warrant 6.3.2.4.* Let an ordinal $\gamma$ be reflected, just if $\VDash^{\gamma} \mathcal{T}B$, provided $\VDash^{\gamma} \mathcal{T}\neg\mathcal{T}\neg B$. Any limit ordinal $\lambda$ is reflected, for if B holds at all ordinals as from some ordinal $\mu$ below $\lambda$ according to $\VDash$, then also $\neg\mathcal{T}\neg B$ holds at all ordinals as from $\mu$ below $\lambda$ according to $\VDash$. So limit ordinals are reflected, and successor ordinals are monogamous, in the sense of Postulate 6.3.2.3. The content of 6.3.2.4 is thus that all ordinals are reflected or monogamous, as for a monogamous successor ordinal $\delta$, ($\VDash^{\delta} (\mathcal{T}A \vee \mathcal{T}\neg A)$, and if $\delta$ is a reflected limit ordinal, $\VDash^{\delta} (\mathcal{T}\neg\mathcal{T}\neg B \to \mathcal{T}B)$. In either case, 6.3.2.4 is warranted. ∎

*Warrant 6.3.2*.5. Let an ordinal $\gamma$ be *transitive* just if for any A,

$$\exists\theta(\theta \prec \gamma \,\&\, \mathbf{\Pi}\xi(\theta \preceq \xi \Rightarrow \VDash^{\xi} \mathrm{A})) \Rightarrow \exists\theta(\theta \prec \gamma \,\&\, \mathbf{\Pi}\xi(\theta \preceq \xi \Rightarrow (\VDash^{\xi} \mathcal{T}\mathrm{A})).$$

Precisely limit ordinals are transitive ordinals.

The content of Postulate 6.3.2.5 is that ordinals are transitive, or monogamous, in the sense of Warrant 6.3.2.3. But that is true, as all ordinals larger than 0 are successor ordinals or limit ordinals. So 6.3.2.5 has been warranted. ∎

*Warrant 6.3.2*.6. At successor ordinals this holds, because there the consequent is true. Let $\lambda$ be a limit ordinal, and $\rho$ such that

$$\mathbf{\Pi}\xi(\rho \preceq \xi \prec \lambda) \Rightarrow \VDash^{\xi} \mathcal{T}\mathrm{A} \to \mathrm{A},$$

so that $\VDash^{\lambda} \mathcal{T}(\mathcal{T}\ulcorner\mathrm{A}\urcorner \to \mathrm{A})$. Suppose there is some ordinal $\sigma \prec \lambda$ and $\rho \preceq \sigma$ such that $\VDash^{\sigma}$ A. If so, $\VDash^{\lambda} \mathcal{T}\mathrm{A}$. If there is no ordinal $\sigma \prec \lambda$ and $\rho \prec \sigma$ such that $\VDash^{\sigma}$ A, then $\VDash^{\lambda} \mathcal{T}\neg\mathrm{A}$. So 6.3.2.6 has been warranted. ∎

*Warrant 6.3.2*.7. Postulate 6.3.2.7 holds at all successor ordinals, as the consequent always holds there. $\VDash^{\lambda} \mathcal{T}(\mathcal{T}\mathrm{A} \to \mathcal{T}\mathcal{T}\mathrm{A}) \Rightarrow \mathbf{\Sigma}\delta\mathbf{\Pi}\epsilon(\delta \preceq \epsilon \prec \lambda \Rightarrow \VDash^{\epsilon} \mathcal{T}\mathrm{A} \to \mathcal{T}\mathcal{T}\mathrm{A})$ if $\lambda$ is a limit ordinal. But all ordinals $\epsilon$ in the interval from and including $\delta$ and less than $\lambda$ will have a successor $\epsilon + 1$ which is also in the interval, so also $\VDash^{\epsilon+1} \mathcal{T}\mathrm{A} \to \mathcal{T}\mathcal{T}\mathrm{A}$. But the latter statement has the consequence that $\VDash^{\epsilon} \mathrm{A} \to \mathcal{T}\mathrm{A}$. So we have established that for any limit $\lambda$, $\VDash^{\lambda} \mathcal{T}(\mathcal{T}\mathrm{A} \to \mathcal{T}\mathcal{T}\mathrm{A}) \to \mathcal{T}(\mathrm{A} \to \mathcal{T}\mathrm{A})$. Given postulate 6.3.2.2 and contraposition, we obtain that $\VDash^{\lambda} \mathcal{T}(\mathcal{T}\mathrm{A} \to \mathcal{T}\mathcal{T}\mathrm{A}) \to \mathcal{T}(\mathcal{T}\neg\mathrm{A} \to \neg\mathrm{A})$. At this point is only takes postulate 6.3.2.6 to finish the warrant. ∎

*Warrant 6.3.2*.8. Suppose $\VDash^{\gamma} \exists\mathrm{v}\mathcal{T}\mathrm{A}$. On account of Definition 5.2.1.9, $\VDash^{\gamma} \mathcal{T}\mathrm{A}^{\mathrm{b}}_{\mathrm{v}}$ for a b substitutable for v in A. So, on account of Definition 5.2.1.17 it follows that for an ordinal $\delta$ and any ordinal $\epsilon$ such that $\delta \preceq \epsilon \prec \gamma$, $\VDash^{\epsilon} \mathrm{A}^{\mathrm{b}}_{\mathrm{v}}$ for a b substitutable for v in A. So on account of Definition 5.2.1.9, again, for an ordinal $\delta$ and any ordinal $\epsilon$ such that $\delta \preceq \epsilon \prec \gamma$, $\VDash^{\epsilon} \exists\mathrm{v}\mathrm{A}$. So on account of Definition 5.2.1.17, $\VDash^{\gamma} \mathcal{T}\exists\mathrm{v}\mathrm{A}$. ∎

*Warrant 6.3.2*.9. Let ordinal $\gamma$ be such that $\VDash^{\gamma} \mathcal{T}\forall\mathrm{v}\mathrm{A}$. There is, consequently, an ordinal $\delta$ such that for any ordinal $\epsilon$ fulfilling $\delta \preceq \epsilon \prec \gamma$, $\VDash^{\epsilon} \forall\mathrm{v}\mathrm{A}$. So either $\gamma = \delta + 1 = \epsilon + 1$ or $\gamma$ is a limit ordinal such that $\VDash^{\epsilon} \forall\mathrm{v}\mathrm{A}$ for all ordinals $\epsilon$ such that $\delta \preceq \epsilon \prec \gamma$. In either case, $\VDash^{\epsilon} \forall\mathrm{v}\mathrm{A}$ holds at any $\epsilon$ smaller than $\gamma$ and at least as large as $\delta$. It follows from Definition 5.2.1.9, that $\VDash^{\epsilon} \mathrm{A}^{\mathrm{b}}_{\mathrm{v}}$, at any $\epsilon$ smaller than $\gamma$ and at least as large as $\delta$, for all b substitutable for v in A. So $\VDash^{\gamma} \mathcal{T}\mathrm{A}^{\mathrm{b}}_{\mathrm{v}}$, for all b substitutable for v in A. So from Definition 5.2.1.9, again, $\VDash^{\gamma} \forall\mathrm{v}\mathcal{T}\mathrm{A}$. So $\VDash^{\beta} \mathcal{T}\forall\mathrm{v}\mathrm{A} \to \forall\mathrm{v}\mathcal{T}\mathrm{A}$ holds at any ordinal $\beta$. So $\VDash^{\varpi} \mathcal{T}(\mathcal{T}\forall\mathrm{v}\mathrm{A} \to \forall\mathrm{v}\mathcal{T}\mathrm{A})$, and consequently $\Vdash^{\mathrm{M}} \mathcal{T}\forall\mathrm{v}\mathrm{A} \to \forall\mathrm{v}\mathcal{T}\mathrm{A}$. ∎

*Warrant 6.3.2*.10. The warrant is in the proof of Theorem 9.3.1.5. ∎

*Warrant 6.3.2.11.* – Notice that 6.3.2.11 is the *Barcan postulate* for orthodox formulas.

Assume

$$\not\Vdash^{M} \mathfrak{D}(\{x|A\}) \rightarrow (\forall x \mathcal{T} A \rightarrow \mathcal{T} \forall x A). \tag{6.4.1}$$

It follows, by Definitions 5.2.1 and 5.5.1, that we for the generalbarter on $n$, $\Vdash_{n}$, have:

$$\Vdash_{n}^{\varpi} \neg \mathcal{T} (\mathfrak{D}(\{x|A\}) \rightarrow (\forall x \mathcal{T} A \rightarrow \mathcal{T} \forall x A)). \tag{6.4.2}$$

Definition 5.2.1.17 has the consequence:

$$\begin{gathered} \mathbf{\Pi}\gamma(\gamma \prec \varpi \;\Rightarrow\; \mathbf{\Sigma}\delta(\gamma \preceq \delta \prec \varpi \;\; \& \\ \Vdash_{n}^{\delta} \mathfrak{D}(\{x|A\}) \;\wedge \forall x \mathcal{T} A \wedge \neg \mathcal{T} \forall x A)). \end{gathered} \tag{6.4.3}$$

As the nethermost ordinal has exceptional features, proceed by instantiating with 1 in Equation (6.4.3) to obtain the consequence

$$\mathbf{\Sigma}\delta\Big(1 \preceq \delta \prec \varpi \;\&\; \Vdash_{n}^{\delta} \mathfrak{D}(\{x|A\}) \;\wedge \forall x \mathcal{T} A \wedge \neg \mathcal{T} \forall x A\Big). \tag{6.4.4}$$

The following reasoning with case distinction establishes that (6.4.4) cannot hold, and so, consequently, that the assumption expressed by Equation (6.4.1) cannot hold.

Case 1/2 - $\delta$ is a limit:

Suppose

$$\Vdash_{n}^{\delta} \mathfrak{D}(\{x|A\}) \wedge \; \forall x \mathcal{T} A \wedge \neg \mathcal{T} \forall x A. \tag{6.4.5}$$

Then, for all constants $c$, and all ordinals $\psi$ larger than a $\xi$ smaller than $\delta$,

$$\Vdash_{n}^{\psi} \mathfrak{D}(\{x|A\}) \wedge A_{x}^{c},$$

so as well

$$\Vdash_{n}^{\psi} \mathfrak{D}(\{x|A\}) \wedge \forall x A.$$

Also, however,

$$\Vdash_{n}^{\delta} \neg \mathcal{T} \forall x A,$$

so that for some $\psi \preceq \phi \preceq \delta$,

$$\Vdash_{n}^{\phi} \neg A_{x}^{c}.$$

So

$$(\Vdash^{\delta}_{n} \mathfrak{D}(\{x|A\}) \wedge \forall x\mathcal{T}A \wedge \neg\mathcal{T}\forall xA$$

cannot hold at a limit ordinal $\delta$.

Case 2/2 - $\delta = \beta + 1$ is a successor ordinal:

Suppose

$$\Vdash^{\delta}_{n} \mathfrak{D}(\{x|A\}) \wedge \ \forall x\mathcal{T}A \wedge \neg\mathcal{T}\forall xA.$$

Then

$$\Vdash^{\beta}_{n} \exists x\neg A,$$

so that there, by Theorem 5.2.2, is a constant c for which

$$\Vdash^{\beta}_{n} \neg A^{c}_{x}.$$

However, as

$$\Vdash^{\delta}_{n} \forall x\mathcal{T}A,$$

and $\delta = \beta + 1$,

also

$$\Vdash^{\beta}_{n} A^{c}_{x}.$$

So

$$\Vdash^{\delta}_{n} \mathfrak{D}(\{x|A\}) \wedge \ \forall x\mathcal{T}A \wedge \neg\mathcal{T}\forall xA$$

cannot hold at a successor ordinal $\delta$.

Cases 1/2 and 2/2 entail that for any ordinal $\alpha \succ 0$ and any $n \in \Omega$,

$$\Vdash^{\alpha}_{n} \mathfrak{D}(\{x|A\}) \rightarrow (\forall x\mathcal{T}A \rightarrow \mathcal{T}\forall xA).$$

By arithmetical generalization and the fact that $\varpi \succ 0$ it follows that

$$\Vdash^{\varpi} \mathcal{T}(\mathfrak{D}(\{x|A\}) \rightarrow (\forall x\mathcal{T}A \rightarrow \mathcal{T}\forall xA)),$$

so that

$$\Vdash^{M} \mathfrak{D}(\{x|A\}) \rightarrow (\forall x\mathcal{T}A \rightarrow \mathcal{T}\forall xA),$$

which warrants Postulate 6.3.2.11's posit. ∎

6.4.6. Remark

The semantic justification for some of the maxims of Postulates 6.3.2.1 – 6.3.2.11 can be lifted from (Bjørdal 2012)(340–341).

6.4.7. Remark

Postulates 6.3.2.6 and 6.3.2.7 originate with (Turner 1990).

6.4.8. Remark

The maxims of Postulates 6.3.2.7 and 6.3.2.8 were not included in (Bjørdal 2012), as the author thought they were both derivable. The warrant of Postulate 6.3.2.8 shows that this was correct for its maxim schema, but the warrant of Postulate 6.3.2.7 suggests that Postulate 6.3.2.3 is needed for its semantical justification.

6.4.9. Remark

Although the converses of Postulates 6.3.2.5 and 6.3.2.6 hold at limit ordinals, they are not maxims, for we may at a sucessor $\sigma$ have that

$$\VDash^{\sigma} (\mathcal{T}\neg A \vee \mathcal{T}A) \wedge \neg\mathcal{T}(\mathcal{T}A \rightarrow A),$$

and it happen for $\{x|x \notin x\} \in \{x|x \notin x\}$ at $\sigma$ or $\sigma + 1$. This contrasts with Remark 69.3.1. (ii) in (Cantini 1996)(396).

6.4.10. Exercise

Let A be *deferent* just if for the barter at $\varpi$ on $n$, $\VDash^{\varpi}_{n}$, $\VDash^{\varpi}_{n} \mathcal{T}A$ $\textcircled{r}$ $\VDash^{\varpi}_{n} \mathcal{T}\neg A$. Show that just deferent formulas are orthodox.

6.4.11. Exercise

Prove that $\Vdash^{O} \forall\vec{v}(\mathcal{T}A \vee \mathcal{T}\neg A) \Rightarrow \Vdash^{M} \forall\vec{v}(\mathcal{T}A \vee \mathcal{T}\neg A)$.

Remark on Exercise 6.4.11: Defining a formula A as orthodox just if $\Vdash^{M} \forall\vec{v}(\mathcal{T}A \vee \mathcal{T}\neg A)$, instead of using Definition 5.7.1.1, is not advisable. For defining a formula as paradoxical just if not orthodox, as in Definition 5.7.4, would then induce an unacceptable extension for the term "paradoxical".

## 6.5 Optima on truth

6.5.1. Postulate

(1) $\Vdash^{O} \mathcal{T}A \rightarrow A$

(2) $\Vdash^{O} \mathcal{T}A \leftrightarrow \mathcal{T}\mathcal{T}A$

(3) $\Vdash^{O} \mathcal{T}\neg\mathcal{T}\neg A \rightarrow \mathcal{T}A$

(4) $\Vdash^{O} \mathcal{T}(\mathcal{T}A \rightarrow \mathcal{T}B) \rightarrow \mathcal{T}(A \rightarrow B)$

(5) $\Vdash^{O} \mathcal{T}(A \rightarrow \mathcal{T}A) \leftrightarrow \mathcal{T}(\mathcal{T}A \rightarrow A)$

(6) $\Vdash^{O} \forall x\mathcal{T}A \rightarrow \mathcal{T}\forall xA$.

(7) $\Vdash^{O} r \notin r$.

(8) $\Vdash^{O} r \notin \{x|x \notin r\}$.

Compare Theorem 2 of (Bjørdal 2012, p. 342).

## 6.6 *Excursus 2:* Warrants of optima on truth

*Warrant 6.5.1*.1. Given Theorem 5.4.2, Definition 5.4.9 and Definition 5.4.1, it follows that for the closure ordinal $\varpi$, $\VDash^{\varpi} \mathcal{T}A$ only if $IN(A)$. As the closure ordinal is a covering ordinal, $IN(A)$ only if $IN^{\varpi}(A)$, so that $\VDash^{\varpi} \mathcal{T}A$ only if $\mathbf{\Pi}\beta(\varpi \preceq \beta \Rightarrow \VDash^{\beta} \mathcal{T}A)$. Consequently $\VDash^{\varpi} \mathcal{T}A$ only if $\VDash^{\varpi+1} \mathcal{T}A$. As $\VDash^{\varpi+1} \mathcal{T}A$ only if $\VDash^{\varpi} A$, $\VDash^{\varpi} \mathcal{T}A$ only if $\VDash^{\varpi} A$. Given Theorem 5.2.4.1, $\VDash^{\varpi} \mathcal{T}A \rightarrow A$. $\blacksquare$

*Warrant 6.5.1*.2. Given Postulate 6.5.1.1, we only need to show $\Vdash^{O} \mathcal{T}A \rightarrow \mathcal{T}\mathcal{T}A$. Given Postulates 6.5.1.7 and 6.5.1.8 and Theorem 5.3.1, $\VDash^{O} \neg\mathcal{T}r \in r \wedge \neg\mathcal{T}r \notin r)$. Given Postulate 6.3.2.5 and Theorem 5.6.6, $\Vdash^{O} \mathcal{T}r \in r \vee \mathcal{T}r \notin r \vee (\mathcal{T}A \rightarrow \mathcal{T}\mathcal{T}A)$. So $\Vdash^{O} \mathcal{T}A \rightarrow \mathcal{T}\mathcal{T}A$. $\blacksquare$

*Warrant 6.5.1*.3. Argue as the warrant of Postulate 6.5.1.2. $\blacksquare$

*Warrant 6.5.1*.4. $\Vdash^{O} \neg\mathcal{T}(A \wedge \neg B) \rightarrow \neg\mathcal{T}\neg\mathcal{T}(A \rightarrow B)$ is an instance of Theorem 6.5.1.4, and $\Vdash^{M} \neg\mathcal{T}\neg\mathcal{T}(A \rightarrow B) \rightarrow \neg\mathcal{T}\neg(\mathcal{T}A \rightarrow \mathcal{T}B)$ is a consequence of Postulate 6.3.2.1. $\blacksquare$

*Warrant 6.5.1*.5. LS to RS: Suppose $\VDash^{\varpi} \mathcal{T}(A \rightarrow \mathcal{T}A)$ and $\nVDash^{\varpi} \mathcal{T}(\mathcal{T}A \rightarrow A)$. Then for $\alpha$ smaller than $\varpi, \mathbf{\Pi}\beta(\alpha \preceq \beta \Rightarrow \VDash^{\beta} (A \rightarrow \mathcal{T}A))$, and for all $\alpha \prec \varpi, \mathbf{\Sigma}\beta(\alpha \preceq \beta \;\&\; \VDash^{\beta} \mathcal{T}A \wedge \neg A)$. Absurd! RS to LS: Suppose $\VDash^{\varpi} \mathcal{T}(\mathcal{T}A \rightarrow A)$ and $\nVDash^{\varpi} \mathcal{T}(A \rightarrow \mathcal{T}A)$. Then for $\alpha$ smaller than $\varpi, \mathbf{\Pi}\beta(\alpha \preceq \beta \Rightarrow \VDash^{\beta} (\mathcal{T}A \rightarrow A))$, and for all $\alpha \prec \varpi, \mathbf{\Sigma}\beta(\alpha \preceq \beta \;\&\; \VDash^{\beta} A \wedge \neg\mathcal{T}A)$. Absurd! Confer the argument for the same result in the Theorem 2 (iv) of (Bjørdal 2012, p. 342). $\blacksquare$

*Warrant 6.5.1*.6. Theorem 5.2.4.1 and Definition 5.5.2.2 entail $\nVdash^{O} \forall x\mathcal{T}A \rightarrow \mathcal{T}\forall xA$ only if for some metanumber n, $\nVDash^{\varpi}_{n} \forall x\mathcal{T}A \rightarrow \mathcal{T}\forall xA$. Given Theorem 5.2.4.1, $\nVDash^{\varpi}_{n} \forall x\mathcal{T}A \rightarrow \mathcal{T}\forall xA$ entails $\VDash^{\varpi}_{n} \forall x\mathcal{T}A$ and $\VDash^{\varpi}_{n} \neg\mathcal{T}\forall xA$. But $\VDash^{\varpi}_{n} \forall x\mathcal{T}A$ entails $\mathbf{\Sigma}\beta(\beta \prec \varpi \;\&\; \mathbf{\Pi}\gamma(\beta \preceq \gamma \prec \varpi \Rightarrow \VDash^{\gamma}_{n} A^{b}_{x}))$ for all b substituable for x in A, and $\VDash^{\varpi}_{n} \neg\mathcal{T}\forall xA$ entails $\mathbf{\Pi}\beta(\beta \prec \varpi \Rightarrow \mathbf{\Sigma}\gamma(\beta \preceq \gamma \;\&\; \neg\forall xA))$, which is impossible. $\blacksquare$

*Warrant 6.5.1.7*. Given Definition 6.3.1 and Theorem 5.3.1, for any n, $\VDash^{\varpi}_{n} r \in r \leftrightarrow \mathcal{T}r \notin r$. Given Postulate 6.5.1.1, for any n, $\Vdash^{\varpi}_{n} \mathcal{T}r \notin r \rightarrow r \notin r$, so that for any n, $\VDash^{\varpi}_{n} r \notin r$.

$\blacksquare$

*Warrant 6.5.1*.8. Given Definition 6.3.1 and Theorem 5.3.1, $\VDash^{\varpi}_{n} r \in \{x|x \notin r\} \leftrightarrow \mathcal{T}r \notin r$. But from Warrant 6.5.1.7, for any n, $\VDash^{\varpi}_{n} r \notin r$ and $\VDash^{\varpi}_{n} r \in r \leftrightarrow \mathcal{T}r \notin r$, so $\VDash^{\varpi}_{n} \neg\mathcal{T}r \notin r$ for any metanumber n. As $\VDash^{\varpi}_{n} r \in \{x|x \notin r\} \leftrightarrow \mathcal{T}r \notin r$ for any metanumber n, also $\VDash^{\varpi}_{n} r \notin \{x|x \notin r\}$ for any metanumber n. $\blacksquare$

# 7 Inference modes

*Some things can only be seen from one side.*

---

## 7.1 The simple inference modes

Only $\neg$, $\mathcal{T}$ and one occurence of a formula variable are allowed in the formulas in the antecedent and the consequent of the simple inference modes. Moreover, $\mathcal{T}$ may only occur once in the antecedent, and in the consequent.

7.1.1. Postulates on simple thetical inference modes:

1 $\Vdash A \Rightarrow \Vdash \mathcal{T}A$

2 $\Vdash A \Leftarrow \Vdash \mathcal{T}A$

3 $\Vdash \mathcal{T}A \Rightarrow \Vdash \neg\mathcal{T}\neg A$

4 $\Vdash \mathcal{T}A \Leftarrow \Vdash \neg\mathcal{T}\neg A$

The corresponding valid, simple *maximal* inference modes of Postulate 7.1.2 can be justified by the valid simple *thetical* inference modes in Postulate 7.1.1.

7.1.2. Postulates on simple maximal inference modes:

1 $\Vdash^{M} A \Rightarrow \Vdash^{M} \mathcal{T}A$

2 $\Vdash^{M} A \Leftarrow \Vdash^{M} \mathcal{T}A$

3 $\Vdash^{M} \mathcal{T}A \Rightarrow \Vdash^{M} \neg\mathcal{T}\neg A$

4 $\Vdash^{M} \mathcal{T}A \Leftarrow \Vdash^{M} \neg\mathcal{T}\neg A$

## 7.2 Involved inference modes

7.2.1. Postulates on quantificational thetical modes:

1 $\Vdash \forall v\mathcal{T}A \Rightarrow \Vdash \mathcal{T}\forall vA$

2 $\Vdash \neg\mathcal{T}\forall vA \Rightarrow \Vdash \neg\forall v\mathcal{T}A$

3 $\Vdash \mathcal{T}\exists vA \Rightarrow \Vdash \exists v\mathcal{T}A$

*Warrant.* Consider first Postulate 7.2.1.1, and assume that $\Vdash \forall v\mathcal{T}A$. So $\Pi\alpha(\alpha \preceq \varpi \Rightarrow \Vdash^{\alpha+1} \forall v\mathcal{T}A)$, thus for a substitutable for v in A, $\Pi\alpha(\alpha \preceq \varpi \Rightarrow \Vdash^{\alpha+1} \mathcal{T}A^{a}_{v})$. It follows that $\Pi\alpha(\alpha \preceq \varpi \Rightarrow \Vdash^{\alpha} A^{a}_{v})$, so that $\Pi\alpha(\alpha \preceq \varpi \Rightarrow \Vdash^{\alpha} \forall vA)$ and $\Vdash \mathcal{T}\forall vA$. Postulates 7.2.1.2 and 7.2.1.2 are warranted similarly. ∎

7.2.2. POSTULATE of the Barcan mode:

$$\Vdash^{M} \forall v \mathcal{T} A \Rightarrow \Vdash^{M} \mathcal{T} \forall v A.$$

*Warrant.* Appeal to Postulates 7.2.1.1 and 7.2.1.2. ∎

7.2.3. POSTULATE FAILURE 1:

$$\Vdash \neg \exists v \mathcal{T} A \nRightarrow \Vdash \neg \mathcal{T} \exists v A.$$

*Warrant.* Confer the discussion about these limitative results in Theorem 12.1.3, plus its Corollary 12.1.5 and Theorem 12.1.5 in § 12.1. ∎

7.2.4. POSTULATE FAILURE 2:

$$\Vdash^{M} \mathcal{T} \exists v A \nRightarrow \Vdash^{M} \exists v \mathcal{T} A$$

*Warrant.* Given Postulate 7.2.3, inference mode 7.2.1.3 does not enter a combination, as that between 7.2.1.1 and 7.2.1.2, to avoid the invalidity. Again, confer with the discussion about such limitative results in Theorem 12.1.3, plus its Corollary 12.1.5 and Theorem 12.1.5 in § 12.1. ∎

7.2.5. POSTULATE ON THETICAL DISTRIBUTIVE MODES

1 $\Vdash^{M} (A \rightarrow B) \Rightarrow (\Vdash A \Rightarrow \Vdash B).$

2 $\Vdash^{M} (A \rightarrow B) \Rightarrow (\Vdash \neg B \Rightarrow \Vdash \neg A).$

3 $\Vdash (A \rightarrow B) \Rightarrow (\Vdash^{M} A \Rightarrow \Vdash B).$

7.2.6. POSTULATE (A USEFUL MODE)

$$( \VDash^{M} (A \rightarrow (B \rightarrow C)) \;\&\; \VDash (A \rightarrow B)) \Rightarrow \VDash B.$$

7.2.7. POSTULATE (THE MAXIM MODE)

$$\Vdash^{M} (A \rightarrow B) \Rightarrow (\Vdash^{M} A \Rightarrow \Vdash^{M} B).$$

7.2.8. REMARK

Postulate 7.2.7 is entailed by Postulates 7.2.5.1 and 7.2.5.2.

7.2.9. POSTULATE (COMPLEX MODES)

1 $\Vdash^{M} \mathcal{T} A \Rightarrow \Vdash^{M} \mathcal{T}(\mathcal{T} A \leftrightarrow A) \wedge \mathcal{T}(\mathcal{T} \neg A \leftrightarrow \neg A)$ (The Tarski mode)

2 $\Vdash^{M} \mathcal{T}A \to A \Rightarrow \Vdash^{M} \mathcal{T}A \lor \mathcal{T}\neg A$

3 $\Vdash^{M} \mathcal{T}\neg\mathcal{T}\neg A \Rightarrow \Vdash^{M} \mathcal{T}A$

4 $\Vdash^{M} \mathcal{T}(\mathcal{T}A \to \mathcal{T}B) \Rightarrow \Vdash^{M} \mathcal{T}(A \to B)$

5 $\Vdash A$ & $\Vdash B \Rightarrow \Vdash \neg\mathcal{T}\neg A \land \neg\mathcal{T}\neg B$

6 $\Vdash^{M} \mathfrak{O}(A(x)) \Rightarrow (\Vdash^{M} \exists xA \Rightarrow \Vdash^{M} A^{a}_{x}$ for some a substitutable for x in A).

7 $\Vdash^{M} A^{a}_{v}$ for any constant a $\Rightarrow \Vdash^{M} \forall vA$

*Warrant 7.2.9.*1. Clearly $\Vdash^{M} \mathcal{T}A \Rightarrow \Vdash^{M} \mathcal{T}(A \land \mathcal{T}A) \land \mathcal{T}(A \land \neg\mathcal{T}\neg A)$. It is librationistically derivable that $\Vdash^{M} \mathcal{T}((A \land \mathcal{T}A) \to \mathcal{T}(A \leftrightarrow \mathcal{T}A))$ and $\Vdash^{M} \mathcal{T}((A \land \neg T\neg A) \to \mathcal{T}(\neg A \leftrightarrow \mathcal{T}\neg A))$, so Postulate 7.2.7 suffices to finish. ∎

*Proof.* (7.2.9.4) Suppose $\VDash^{\varpi} \mathcal{T}(\mathcal{T}A \to \mathcal{T}B)$. (i) Let $\rho$ be be a ordinal as from which $\mathcal{T}A \to \mathcal{T}B$ holds, so that

$$\mathbf{\Pi}\xi(\rho \preceq \xi \prec \varpi \Rightarrow \VDash^{\xi} (\mathcal{T}A \to \mathcal{T}B).$$

Thus $\VDash^{\rho+1} (\mathcal{T}A \to \mathcal{T}B)$, and therefore $\VDash^{\rho} (A \to B)$. Consequently, succeeding successors will have $\mathcal{T}A \to \mathcal{T}B$ and $A \to B$. (ii) Let limit ordinal $\lambda \prec \varpi$, above $\rho$, have $\mathcal{T}A \to \mathcal{T}B$, and $A \to B$ below, as from $\rho$. As $\lambda \prec \varpi$, from the assumption on $\rho$, $\VDash^{\lambda} (\mathcal{T}A \to \mathcal{T}B)$. As $\Vdash^{\lambda+1} (\mathcal{T}A \to \mathcal{T}B)$, also $\Vdash^{\lambda} (A \to B)$. (iii) By a repetition of (i) and (ii) it follows that $A \to B$ holds as from $\rho$ below $\varpi$, so that $\VDash^{\varpi} \mathcal{T}(A \to B)$. □

*Proof.* (7.2.9.6) This is established on page 52, in the proof of Theorem 12.2.1. □

# 8 The failure of Löb's rule and formula for $\mathcal{T}$ in £

> *"Wie das Wort „schön" der Ästhetik und „gut" der Ethik, so weist „wahr" der Logik die Richtung. Zwar haben alle Wissenschaften Wahrheit als Ziel; aber die Logik beschäftigt sich noch in ganz anderer Weise mit ihr. Sie verhält sich zur Wahrheit etwa so wie die Physik zur Schwere oder zur Wärme. Wahrheiten zu entdecken, ist Aufgabe aller Wissenschaften: der Logik kommt es zu, die Gesetze des Wahrseins zu erkennen."*
>
> Gottlob Frege, in (Frege 1918, p. 58).

8.0.1. Definition of Löb's rule for predicate $\mathfrak{T}$ in theory $\mathbb{T}$:

If $\vdash_{\mathbb{T}}$ is the theorem sign and [] the coding function for $\mathbb{T}$, Löb's rule for predicate $\mathfrak{T}$ in $\mathbb{T}$ is

$$\vdash_{\mathbb{T}} \mathfrak{T}[\mathrm{B}] \to \mathrm{B} \;\Rightarrow\; \vdash_{\mathbb{T}} \mathrm{B}.$$

8.0.2. Reminder: `Bew` is the provability predicate of (Gödel 1931) for Peano arithmetic.

8.0.3. Observation:

Löb's *theorem*, as per (Löb 1955), is that Löb's *rule* holds for predicate `Bew` in Peano arithmetic.

8.0.4. Theorem to the effect that Löb's rule fails for $\mathcal{T}$ in £:

$$\Vdash^{\mathrm{M}} \mathcal{T}\ulcorner \mathrm{B} \urcorner \to \mathrm{B} \not\Rightarrow \Vdash^{\mathrm{M}} \mathrm{B}.$$

*Proof.* A consequence of §6.2 and §7 is that

$$\Vdash^{\mathrm{M}} \mathrm{A} \to \mathrm{A}, \tag{8.0.5}$$

so that on account of Postulate 7.1.2.1,

$$\Vdash^{\mathrm{M}} \mathcal{T}\ulcorner \mathrm{A} \to \mathrm{A} \urcorner. \tag{8.0.6}$$

Consequently, on account or Maxim 6.3.2.2 and the Maxim mode 7.2.7:

$$\Vdash^{\mathrm{M}} \neg\mathcal{T}\ulcorner \mathrm{A} \wedge \neg\mathrm{A} \urcorner. \tag{8.0.7}$$

From Equation 8.0.7, maxim $\Vdash^{\mathrm{M}} \neg\mathcal{T}\ulcorner \mathrm{A} \wedge \neg\mathrm{A} \urcorner \to (\mathcal{T}\ulcorner \mathrm{A} \wedge \neg\mathrm{A} \urcorner \to (\mathrm{A} \wedge \neg\mathrm{A})$ and the Maxim mode 7.2.7 it follows that

$$\Vdash^{\mathrm{M}} \mathcal{T}\ulcorner \mathrm{A} \wedge \neg\mathrm{A} \urcorner \to \mathrm{A} \wedge \neg\mathrm{A}. \tag{8.0.8}$$

Given the non-triviality assumptions of §5.8,

$$\nVdash^{\mathrm{M}} \mathrm{A} \wedge \neg\mathrm{A}. \tag{8.0.9}$$

The proof finishes by substituting $\mathrm{A} \wedge \neg\mathrm{A}$ with B in Equations 8.0.8 and 8.0.8. □

8.0.10. Definition of Löb's formula for predicate $\mathfrak{T}$ in theory $\mathbb{T}$; cfr. Definition 8.0.1:

$$\mathfrak{T}[\mathfrak{T}[\mathrm{B}] \to \mathrm{B}] \to \mathfrak{T}[\mathrm{B}].$$

8.0.11. Theorem: Löb's formula for predicate $\mathcal{T}$ in **£** is not a maximal theorem of **£**.

*Proof.* Assume $\Vdash^{\mathrm{M}} \mathcal{T}\ulcorner\mathcal{T}\ulcorner\mathrm{B}\urcorner \to \mathrm{B}\urcorner \to \mathcal{T}\ulcorner\mathrm{B}\urcorner$. Suppose $\Vdash^{\mathrm{M}} \mathcal{T}\ulcorner\mathrm{B}\urcorner \to \mathrm{B}$, so that on account of Postulate 7.1.2.1, $\Vdash^{\mathrm{M}} \mathcal{T}\ulcorner\mathcal{T}\ulcorner\mathrm{B}\urcorner \to \mathrm{B}\urcorner$. Given the assumption plus the Maxim Mode 7.2.7, $\Vdash^{\mathrm{M}} \mathcal{T}\ulcorner\mathrm{B}\urcorner$. Thus $\Vdash^{\mathrm{M}}$ B on account of Postulate 7.1.2.2. So the presence of Löb's formula as a maximal theorem in **£** would have entailed the validity of Löb's rule for **£**. But Löb's rule fails for $\mathcal{T}$ in **£** by Theorem 8.0.4, it is not a theorem schema of **£** that Löb's formula is a maxim. □

# 9 The theory of identity

*"To be is to have your essence if everything has it."*

*Source:* Theorem 9.1.3, page 40.

A streamlining of sections 4 and 5 of (Bjørdal 2012, pp. 342–345) is obtained from the inference modes 7.2.9.1 – 7.2.9.4, and as a result **£** does not, as e.g. the comparable systems studied by (Cantini 1996), need additional axiomatic principles for having a well behaved notion of identity in this section, or natural number in §11.

## 9.1 Co-essentiality

In a lasting contribution (Whitehead and Russell 1910) improves upon *Leibniz' law*, as they show that a bi-conditional corresponding with our Definition 9.1.2 is adequate for identity. In Principa Mathematica the theorem corresponding with Definition 9.1.2 is *13.101, which is proven via its predicative Definition *13.1 and its *Axiom of Reducibility* *12.1.

It bears mentioning that the second Principia author published the thorough monograph (Russell 1900) on Leibniz, though this does not establish that he contributed theorem *13.101. For a cursory reading of that historical treatise suggests that Russell did not make such discoveries while writing that text.

We define the identity relation by means of a notion of *co-essentiality*, which is similar to the relation named *membership congruency* by Abraham A. Fraenkel and Yehoshua Bar-Hillel, and discussed in (A. A. Fraenkel et al. 1973, p. 27), though not used in the previous edition (A. A. Fraenkel and Bar-Hillel 1958).

9.1.1. Definition

Sets a and b are *co-essential* just if $\forall u(a \in u \rightarrow b \in u)$.

(Forster 2019) relates that (Hailperin 1944) "gave the first of a number of finite axiomatisations of NF now known. Many of them exploit the function $x \mapsto \{y|y \in x\}$ which is injective and total and is an $\in$ -isomorphism. This function was known to Whitehead, who suggested to Quine that $\{y|x \in y\}$ should be called the "essence" of x (a terminology clearly suggested by a view of sets as properties-in-extension)." Incidentally, Quine was Whitehead's student while studying at Harvard, but Quine got his doctorate twelve years before the publication of (Hailperin 1944).

9.1.2. Definition (Identity as co-essentiality)

$$a = b \;==\; \forall u(a \in u \rightarrow b \in u)$$

9.1.3. Theorem

*"To be is to have your essence if everything has it."*

*Proof.* Assume that $\exists x(x = b)$. Given Definition 9.1.2 it follows that $\exists x \forall u(x \in u \rightarrow b \in u)$. From the latter it follows instantiation and logic that

$$\exists x(x \in \{y|y \in b\} \rightarrow b \in \{y|y \in b\}).$$

Quantifier rules entail that

$$\forall x(x \in \{y|y \in b\}) \rightarrow b \in \{y|y \in b\},$$

so that by discharging the assumption $\exists x(x = b)$ one gets

$$\exists x(x = b) \rightarrow (\forall x(x \in \{y|y \in b\}) \rightarrow b \in \{y|y \in b\}).$$

The other direction is trivial, so that

$$\exists x(x = b) \leftrightarrow (\forall x(x \in \{y|y \in b\}) \rightarrow b \in \{y|y \in b\}).$$

The theorem follows from the last sentence by assuming Whitehead's theory on set theoretic essence, stated just below Definition 9.1.1. □

The justification for the analogous definition ∗13·01 in Principia Mathematica will not justify Definition 9.1.2 in **£**. For the symmetry of Definition 9.1.2, is in **£** shown by the proof of Theorem 9.3.1.4 below, and that proof does not appeal to *predicativity* or the Axiom of reducibility, as does the proof of ∗13·01 by (Whitehead and Russell 1910).

## 9.2 The identity lemmas

9.2.1. Lemma

$$\Vdash^M \mathcal{T}(\forall u(a \in u \rightarrow b \in u) \rightarrow \mathcal{T}\forall u(a \in u \rightarrow b \in u)).$$

*Proof.* By instantiation we have:

$$\begin{aligned}&\Vdash^M \forall u(a \in u \rightarrow b \in u) \rightarrow \\ &\qquad (a \in \{v|\forall u(a \in u \rightarrow v \in u)\} \rightarrow b \in \{v|\forall u(a \in u \rightarrow v \in u)\}).\end{aligned}$$

But $\Vdash^M a \in \{v|\forall u(a \in u \rightarrow v \in u)\}$, so that

$$\Vdash^M \forall u(a \in u \rightarrow b \in u) \rightarrow b \in \{v|\forall u(a \in u \rightarrow v \in u)\}.$$

Finish with Alethic Comprehension and Postulate 7.1.2.1. □

9.2.2. Lemma

$$\Vdash^{M} \mathcal{T}(a = b \to \mathcal{T} a = b)$$

*Proof.* Use Definition 9.1.2 and Lemma 9.2.1. □

9.2.3. Lemma

$$\mathcal{T}(\mathcal{T} \ulcorner a \neq b \urcorner \to a \neq b)$$

*Proof.* Use Lemma 9.2.2, Postulates 6.3.2.2, 7.1.2.1 and logic. □

9.2.4. Lemma

$$\Vdash^{M} \mathcal{T}(\mathcal{T} \ulcorner \neg\forall u(a \in u \to b \in u) \urcorner \to \neg\forall u(a \in u \to b \in u))$$

*Proof.* Use Lemma 9.2.3 and Definition 9.1.2. □

9.2.5. Lemma

$$\Vdash^{M} \mathcal{T}((\mathcal{T} \forall u(a \in u \to b \in u)) \to \forall u(a \in u \to b \in u))$$

*Proof.* Combine Lemma 9.2.4 with Postulate 6.5.1.5 to obtain

$$\Vdash^{M} \mathcal{T}(\neg\forall u(a \in u \to b \in u) \to \mathcal{T}\neg\forall u(a \in u \to b \in u)).$$

An instance of Theorem 6.3.2.2 is

$$\Vdash^{M} \mathcal{T}(\mathcal{T}\neg\forall u(a \in u \to b \in u) \to \neg\mathcal{T}\forall u(a \in u \to b \in u)).$$

A hypothetical syllogism and contraposition now suffices to finish the proof. □

## 9.3 The adequacy of identity as co-essentiality

9.3.1. Theorem (Orthodoxy, equivalence and fungibility)

(1) $\Vdash^{M} \mathcal{T} a = b \vee \mathcal{T} a \neq b$ — Orthodoxy

(2) $\Vdash^{M} a = a$ — Reflexivity

(3) $\Vdash^{M} a = b \wedge b = c \to a = c$ — Transitivity

(4) $\Vdash^{M} a = b \to b = a$ — Symmetry

(5) $\Vdash^{M} a = b \to (A^{a}_{v} \to A^{b}_{v})$, with a and b substitutable for v in A. — Fungibility

*Proof.*

1. Use Lemma 9.2.3 and Postulate 6.3.2.6.
2. Trivial
3. Trivial, given Definition 9.1.2

4. Clearly,

$$\begin{aligned}&\Vdash^{M} \forall v(a \in v \rightarrow b \in v) \rightarrow \\ &\qquad (a \in \{w|\forall v(w \in v \rightarrow a \in v)\} \rightarrow b \in \{w|\forall v(w \in v \rightarrow a \in v)\}).\end{aligned} \tag{9.3.2}$$

But

$$\Vdash^{M} a \in \{w|\forall v(w \in v \rightarrow a \in v)\},$$

so that by alethic comprehension,

$$\Vdash^{M} \forall v(a \in v \rightarrow b \in v) \rightarrow \mathcal{T}\forall v(b \in v \rightarrow a \in v). \tag{9.3.3}$$

An instance of Lemma 9.2.2 states:

$$\Vdash^{M} \forall v(b \in v \rightarrow a \in v) \rightarrow \mathcal{T}\forall v(b \in v \rightarrow a \in v). \tag{9.3.4}$$

By invoking 6.5.1.5 on equation 9.3.4 we obtain

$$\Vdash^{M} \mathcal{T}\ulcorner\forall v(b \in v \rightarrow a \in v)\urcorner \rightarrow \forall v(b \in v \rightarrow a \in v). \tag{9.3.5}$$

Finish with a hypothetical syllogism with equations 9.3.3 and 9.3.5, and lastly an appeal to co-essentiality Definition 9.1.2.

5. The promissory note issued in sentence Warrant 6.3.W10 of Postulate 6.3.2.10 on page 29 of § 6 is satisfied, and the mentioned Postulate is warranted.

Suppose, for a and b substitutable for v in A, and the generalbarter at $\varpi \Vdash^{\varpi}$,

$$\nVdash^{\varpi} \mathcal{T}(\forall u(a \in u \rightarrow b \in u) \rightarrow (A^{a}_{v} \rightarrow A^{b}_{v})).$$

On account of the validity of the mode of 7.2.9.3 we get

$$\nVdash^{\varpi} \mathcal{T}\neg\mathcal{T}(\forall u(a \in u \rightarrow b \in u) \wedge A^{a}_{v} \wedge \neg A^{b}_{v}).$$

It follows from Definition 5.2.1.1 that

$$\Vdash^{\varpi} \neg\mathcal{T}\neg\mathcal{T}(\forall u(a \in u \rightarrow b \in u) \wedge A^{a}_{v} \wedge \neg A^{b}_{v}).$$

On account of Postulate 6.3.2.1,

$$\Vdash^{\varpi} \neg\mathcal{T}\neg(\mathcal{T}\forall u(a \in u \rightarrow b \in u) \wedge \mathcal{T}A^{a}_{v} \wedge \mathcal{T}\neg A^{b}_{v}).$$

On account of Lemma 9.2.5, we get

$$\Vdash^{\varpi} \neg\mathcal{T}\neg(\forall u(a \in u \rightarrow b \in u) \wedge \mathcal{T}A^{a}_{v} \wedge \neg\mathcal{T}A^{b}_{v}).$$

From alethic comprehension and existential generalization we obtain

$$\Vdash^{\varpi} \neg\mathcal{T}\neg(\forall u(a \in u \rightarrow b \in u) \land \exists u(a \in u \land b \notin u)),$$

which is absurd. So Postulate 6.3.2.10 is tautological, and we are done. $\square$

# 10 Alphabetologicality

*"Mathematics, rightly viewed, possesses not only truth, but supreme beauty — a beauty cold and austere, like that of sculpture, without appeal to any part of our weaker nature, without the gorgeous trappings of painting or music, yet sublimely pure, and capable of a stern perfection such as only the greatest art can show."*

Bertrand Russell in "The Study of Mathematics", quoted from (Russell 1907, p. 33; Russell 1910, p. 73; Russell 1917, p. 60 and Russell 1918, p. 60).

Postulates 10.0.1 and 10.0.2 express, given Definition 10.0.3, that identity is an equivalence relation which is neutral with respect to alphabetological variants.

10.0.1. Postulate (The Lindenbaum-Tarski closure for identity)

If classical logic proves that $\forall x(A(x) \leftrightarrow B(x))$, then

$$\Vdash^{M} \{x|A(x)\} = \{x|B(x)\}.$$

10.0.2. Postulate on alphabetical indifference

$$\{x|A(x)\} = \{x|B(x)\} \rightarrow \{x|A(x)\} = \{y|B(x)^{y}_{x}\},$$

where y is substitutable for x in B.

10.0.3. Definition of alphabetologicality

Two sets are *alphabetological variants* of each other just if they are identical on account of Postulates 10.0.1 and 10.0.2.

Postulates 10.0.1 and 10.0.2 compensate somewhat for the loss of extensionality in **£**, as per § 19, and secure such theorems as:

$$\Vdash^{M} \{x|A(x)\} = \{y|A(y) \wedge \exists z(B(z) \vee \neg B(z))\}.$$

# 11 Arithmetic

*"Die ganzen Zahlen hat der liebe Gott gemacht, alles andere ist Menschenwerk."*

Leopold Kronecker to the Berliner Naturforscher-Versammlung in 1886, according to (H. Weber 1893, p. 19).

11.0.1. Definition

(1) $\emptyset = \{x | x \neq x\}$.

(2) $a' = \{x | x = a \vee x \in a.\}$

(3) $\omega = \{x | \forall y(\emptyset \in y \wedge \forall z(z \in y \rightarrow z' \in y) \rightarrow x \in y)\}$.

11.0.2. Theorem

(1) $\Vdash^{M} \emptyset \in \omega$

(2) $\Vdash^{M} \forall x(x \in \omega \rightarrow x' \in \omega)$

(3) $\omega$ is orthodox

(4) $\Vdash^{M} \forall y(\emptyset \in y \wedge \forall z(z \in y \rightarrow z' \in y) \rightarrow \forall x(x \in \omega \rightarrow x \in y))$

(5) $\Vdash^{M} A(\emptyset) \wedge \forall x(A(x) \rightarrow A(x')) \rightarrow \forall y(y \in \omega \rightarrow A(y))$

*Proof:*

1. Combine alethic comprehension and the fact that

$$\Vdash^{M} \mathcal{T} \forall y(\emptyset \in y \wedge \forall z(z \in y \rightarrow z' \in y) \rightarrow \emptyset \in y)$$

2. This follows from alethic comprehension and

$$\Vdash^{M} \forall x(\mathcal{T}(\forall y(\emptyset \in y \wedge \forall z(z \in y \rightarrow z' \in y) \rightarrow x \in y)) \rightarrow \mathcal{T}(\forall y(\emptyset \in y \wedge \forall z(z \in y \rightarrow z' \in y) \rightarrow x' \in y))).$$

3. From logic:

$$\Vdash^{M} \emptyset \in \omega \wedge \forall x(x \in \omega \rightarrow x' \in \omega) \rightarrow (\forall y(\emptyset \in y \wedge \forall x(x \in y \rightarrow x' \in y) \rightarrow a \in y) \rightarrow a \in \omega).$$

By combining 1 and 2 we have

$$\Vdash^{M} \forall y(\emptyset \in y \wedge \forall x(x \in y \rightarrow x' \in y) \rightarrow a \in y) \rightarrow a \in \omega).$$

Postulates 6.3.2.1 and 7.1.2.1, and alethic comprehension, give us

$$\Vdash^{M} a \in \omega \rightarrow \mathcal{T} a \in \omega.$$

7.2.7 along with Postulates 6.3.2.1, 6.3.2.2 and 6.3.2.6 give us

$$\Vdash^{M} \mathcal{T}a \in \omega \vee \mathcal{T}a \notin \omega$$

As a was arbitrary, $\Vdash^{M} \forall x(\mathcal{T}x \in \omega \vee \mathcal{T}x \notin \omega)$, and the proof is finished.

4. Immediate, given 3, as it is equivalent with

$$\Vdash^{M} \forall x(x \in \omega \rightarrow \forall y(\emptyset \in y \wedge \forall z(z \in y \rightarrow z' \in y) \rightarrow x \in y)).$$

5. For the following, compare (Cantini 1996, p. 356):

Presuppose for arbritrary sentence $A(x)$ that

11.0.3. Definition

$$A'(x) \Longleftrightarrow A(\emptyset) \wedge \forall x(A(x) \rightarrow A(x')) \rightarrow A(x))$$

is fulfilled.

By logic,

$$\Vdash^{M} A'(\emptyset) \ \& \ \Vdash^{M} \forall x(A'(x) \rightarrow A'(x')).$$

The inference mode of Postulate 7.1.2.1 and maxim of Postulate 6.3.2.9 entail

$$\Vdash^{M} \mathcal{T}A'(\emptyset) \ \& \ \Vdash^{M} \forall x\mathcal{T}(A'(x) \rightarrow A'(x')).$$

By quantifier distribution and Postulate 6.3.2.1 we get

$$\Vdash^{M} \mathcal{T}A'(\emptyset) \ \& \ \Vdash^{M} \forall x(\mathcal{T}A'(x) \rightarrow \mathcal{T}A'(x')).$$

Alethic comprehension gives us

$$\Vdash^{M} \emptyset \in \{y|A'(y)\} \ \& \ \Vdash^{M} \forall x(x \in \{y|A'(y)\} \rightarrow x' \in \{y|A'(y)\}).$$

Adjunction gives us

$$\Vdash^{M} \emptyset \in \{y|A'(y)\} \wedge \forall x(x \in \{y|A'(y)\} \rightarrow x' \in \{y|A'(y)\}).$$

4 and the inference of mode 7.2.7 give us

$$\Vdash^{M} \forall x(x \in \omega \rightarrow x \in \{y|A'(y)\}).$$

From 3 and 7.2.9.1 we have

$$\Vdash^{M} \forall x(\mathcal{T}x \in \omega \rightarrow x \in \omega),$$

so that

$$\Vdash^{M} \forall x(\mathcal{T}x \in \omega \rightarrow x \in \{y|A'(y)\}).$$

Alethic comprehension gives us

$$\Vdash^{M} \forall x(\mathcal{T}x \in \omega \rightarrow \mathcal{T}A'(x)),$$

which, combined with 7.2.9.4 establishes

$$\Vdash^{M} \forall x(x \in \omega \rightarrow A'(x))$$

Finish with an appeal to Definition 11.0.3, and rearrangement. □

# 12 Some shortcomings and redresses

*"If people do not believe that mathematics is simple, it is only because they do not realize how complicated life is."*

John von Neumann, 1947, to the Association for Computing Machinery, according to the recollection of (Alt 1972, p. 694).

Some drawbacks and imperfections of librationism are presented, and discussed.

## 12.1 Deficiencies related to existential instantiation

Despite the important Theorem 5.2.2, which justifies

12.1.1. Theorem

$$\Vdash^{\alpha} \exists x(\neg\mathcal{T}A \wedge \neg\mathcal{T}\neg A) \Rightarrow \text{for some term b}, \Vdash^{\alpha} (\neg\mathcal{T}A^{b}_{x} \wedge \neg\mathcal{T}\neg A^{b}_{x}),$$

and consequently

12.1.2. Theorem (Optimal existential instantiation)

$$\text{If } \Vdash^{O} \exists xA, \text{then } \Vdash^{O} A^{b}_{x} \text{ for some term b}$$

There is, nevertheless, as pointed to in Postulate 7.2.4, the following limitative result:

12.1.3. Theorem (Maximal lack of existential instantiation)

$$\text{It may happen that } \Vdash^{M} \exists xA, \text{ and for no term a}, \Vdash^{M} A^{a}_{x}.$$

*Proof.* As the proof of Theorem 12.1.5. □

12.1.4. Corollary

Maximal existential instantiation, in the form

$$\Vdash^{M} \exists xA \Rightarrow \mathbf{\Sigma}a \Vdash^{M} A^{a}_{x},$$

is not valid.

12.1.5. Theorem

The schema

$$\Vdash^{M} \mathcal{T}\exists xA \Rightarrow \Vdash^{M} \exists x\mathcal{T}A$$

is not a valid inference mode.

*Proof.* Let A be $(x = \emptyset \leftrightarrow r \in r)$.

Suppose that $\Vdash^{M} \exists x \mathcal{T} A$. If so $\Vdash^{\varpi} \mathcal{T} \exists x \mathcal{T} A$, and there is an ordinal $\gamma$ such that $\Vdash^{\beta} \exists x \mathcal{T} A$ holds whenever $\gamma \prec \beta \prec \varpi$. Let limit ordinal $\lambda$ satisfy $\gamma \prec \lambda \prec \varpi$, so that $\Vdash^{\lambda} \exists x \mathcal{T} A$. On account of 5.2.1.8 and 5.2.1.9, for a term a, and ordinal $\delta$, $a = \emptyset \leftrightarrow r \in r$ holds at all ordinals $\theta$ which satisfy $\delta \prec \theta \prec \lambda$. But this is impossible, as $r \in r$ holds at some of those ordinals, and $r \notin r$ holds at others, whereas identity is orthodox. □

Theorem 12.1.5 entails that the schema of Postulate 7.2.1.3 does not hold as a maxim, for, as the proof of Postulate 7.2.1.3 just showed, some instances of the schema $\mathcal{T} \exists v A \rightarrow \exists v \mathcal{T} A$ are paradoxical truths.

## 12.2 An orthodox redress

12.2.1. Theorem

$$\Vdash^{M} \mathfrak{O}(A(x)) \Rightarrow (\Vdash^{M} \exists x A \Rightarrow \Vdash^{M} A^{a}_{x}, \text{ for some a substitutable for x in A}).$$

*Proof.* The following proof of Theorem 12.2.1 was announced in Postulate 7.2.9.6 on page 36, and the latter Postulate amounts to the same as the former Theorem:

$$\text{Assume that } A(x) \text{ is orthodox, i.e. } \Vdash^{M} \mathcal{T} A(x) \vee \mathcal{T} \neg A(x).$$

By soundness,

$$\Vdash^{M} \exists x A \Rightarrow \text{for the barter at } \varpi, \VDash^{\varpi} \mathcal{T} \exists x A.$$

$$\text{As } \varpi \text{ is a stabilising ordinal, } \VDash^{\varpi} \exists x A.$$

$$\text{Given Definition 5.2.1 and Theorem 5.2.2, for a a, } \VDash^{\varpi} A^{a}_{x}.$$

$$\text{As } A(x) \text{ is orthodox, } \VDash^{\varpi} \mathcal{T} A^{a}_{x}.$$

$$\text{So } \Vdash^{M} A^{a}_{x}.$$

$$\text{So 7.2.9.6, and thus is valid.}$$

□

## 12.3 The maximal Barcan failure

It will be shown that the truth theoretic Barcan schema does not hold as a maxim.

The precursor to this negative result, in a truth theoretic context, is *McGee's paradox*, in (McGee 1985), which we adapt to our set theoretic context. Compare (Cantini 1996, pp. 380–382) and (Bjørdal 2012, p. 537).

First we decide upon some notions:

12.3.1. Reminder: For r in 12.3.2.5, recall Definition 6.3.1.

12.3.2. Definition

(1) $a' == \{x | x \in a \vee x = a\}$.

(2) $\{a, b\} == \{x | x = a \vee x = b\}$.

(3) $\{a\} == \{a, a\}$.

(4) $a_\omega == \{u | \forall x (\langle \emptyset, a \rangle \in x \wedge \forall y, z (\langle y, z \rangle \in x \rightarrow \langle y', \{v | v \in z\} \rangle) \rightarrow u \in x)\}$.

(5) $t == \{x | x = r \wedge x \notin x \wedge \neg \mathcal{T} x \in x\}$.

(6) Use $\overline{0}$, $\overline{1}$, $\overline{2}$, …for the members of $\omega$.

(7) Let $t_{\overline{0}} == t$ and $t_{\overline{n+1}} == \{v | v \in t_{\overline{n}}\}$.

(8) $B(t_{\overline{i}}) == \exists w (\langle w, t_{\overline{i}} \rangle \in t_\omega) \rightarrow r \notin t_{\overline{i}}$

(9) $B(x) == \exists w (\langle w, x \rangle \in t_\omega) \rightarrow r \notin x$

12.3.3. Lemma

For any $a, a_\omega$ is orthodox.

*Proof.* Adapt the proof of Theorem 11.0.2.3. □

12.3.4. Lemma

$\Vdash^\lambda r = r \wedge r \notin r \wedge \neg \mathcal{T} r \in r$ just if $\lambda$ is a limit.

*Proof.* For any successor ordinal $\chi + 1$, $\Vdash^{\chi+1} \neg \mathcal{T} r \in r \leftrightarrow r \in r$. So the Lemma holds as at precisely limit ordinals $\lambda$, $\Vdash^\lambda r \notin r \wedge \neg \mathcal{T} r \in r$. □

12.3.5. Lemma

For any limit ordinal $\alpha \prec \varpi$, and $\beta = \alpha + \omega$ :

1. $\Vdash^\beta \forall x \mathcal{T} B(x)$
2. $\Vdash^\beta \neg \mathcal{T} \forall x B(x)$.

*Proof.*

1. If $\Vdash^\beta \neg \exists w (\langle w, x \rangle \in t_\omega)$, it follows that $\Vdash^\beta \mathcal{T} B(x)$ on account of Lemma 12.3.3. If, on the other hand, $\Vdash^\beta \exists w (\langle w, x \rangle \in t_\omega)$ we have that $\Vdash^\beta \mathcal{T} B(x)$, as there is a $\gamma \geq \alpha + i$ such that
$$\forall \delta (\alpha \prec \gamma \preceq \delta \prec \beta \Rightarrow \Vdash^\delta B(x)).$$
So for any x, $\Vdash^\beta \mathcal{T} B(x)$, and consequently $\Vdash^\beta \forall x \mathcal{T} B(x)$.

2. Otherwise, $\Vdash^\beta \mathcal{T} \forall x B(x)$, and we would have $\Vdash^\delta \forall x B(x)$ as from some ordinal $\delta$ below $\beta$ and above $\alpha$. Let $\delta == \alpha + (n + 1)$, for finite ordinal $n \geq 0$, be such an ordinal. By instantiation, $\Vdash^\delta B(t_{\overline{n}})$, this entails that $\Vdash^{\alpha+(n+1)} B(t_{\overline{n}})$. As $\Vdash^{\alpha+(n+1)} \exists w (\langle w, t_{\overline{n}} \rangle \in t_\omega)$, it follows that $\Vdash^{\alpha+(n+1)} r \notin t_{\overline{n}}$. As a consequence, $\Vdash^{\alpha+1} r \notin t_{\overline{0}}$. But the latter entails $\Vdash^{\alpha}$ $(r \neq r \vee r \in r \vee \mathcal{T} r \in r)$ which contradicts Lemma 12.3.4, as $\alpha$ is presupposed to be a limit ordinal. □

Despite Postulate 6.5.1.6, with its warrant in §6.6,

12.3.6. Theorem

$$\nVdash^{M} \forall x \mathcal{T} B(x) \rightarrow \mathcal{T} \forall x B(x).$$

*Proof.* Given Lemma 12.3.5,

$$\Vdash^{\varpi} \neg \mathcal{T} \neg (\forall x \mathcal{T} B(x) \wedge \neg \mathcal{T} \forall x B(x)).$$

An appeal to Definition 5.5.2.1 suffices to finish the proof. □

12.3.7. Theorem

For some formula $A(x)$,

$$\Vdash \forall x \mathcal{T} A(x) \rightarrow \mathcal{T} \forall x A(x) \ \ \& \ \ \Vdash \forall x \mathcal{T} A(x) \wedge \neg \mathcal{T} \forall x A(x)$$

*Proof.* Let $B(x)$ in Theorem 12.3.6 be $A(x)$, and combine with Postulate 6.5.1.6. □

## 12.4 Omega-consistency as confinement

(McGee 1985) isolated a theory of truth which is consistent but $\omega$-inconsistent. (Friedman and Sheard 1987) proposed a more substantial a theory of truth, which inherits that $\omega$-inconsistency property, and (Halbach 1994) found that its proof-theoretic strength is the same as the theory of ramified analysis for all finite levels.

An essential ingredient in the proof of McGee's negative result fails in **£**, viz. the equation

$$\forall x(x \in \omega \rightarrow \mathcal{T} \ulcorner A(x) \urcorner) \rightarrow \mathcal{T} \ulcorner \forall x(x \in \omega \rightarrow A(x) \urcorner), \tag{12.4.1}$$

which is used at (McGee 1985, p. 399). Notice that

$$\Vdash^{M} \forall x(x \in \omega \rightarrow \mathcal{T} \ulcorner A(x) \urcorner) \leftrightarrow \forall x \mathcal{T} \ulcorner x \in \omega \rightarrow A(x) \urcorner,$$

as $\omega$ is orthodox. So the maximal version of 12.4.1 follows from the appropriate maximal version of the Barcan-formula, but the latter's maximhood fails for **£** on account of Theorem 12.3.6. For such reasons, the Mcgee argument for omega inconsistency cannot be carried through for **£**.

## 12.5 More orthodox redresses

12.5.1. THEOREM (ORTHODOX ATTESTOR)

If $A(x)$ is orthodox, then

$$\Vdash^{M} \mathcal{T}\exists x A(x) \Rightarrow \Vdash^{M} \exists x \mathcal{T} A(x).$$

*Proof.* Appeal to Theorem 12.2.1, and existential generalization. □

12.5.2. THEOREM (THE BARCAN FORMULA HOLDS FOR ORTHODOX FORMULAS)

$$\Vdash^{M} \mathfrak{O}(B(x)) \Rightarrow \Vdash^{M} (\forall x \mathcal{T} B(x) \rightarrow \mathcal{T}\forall x B(x)).$$

*Proof.* As on page 36. □

# 13 Classicalities and deviations

> *"On a signalé beaucoup d'antinomies, et le désaccord a subsisté, personne n'a été convaincu; d'une contradiction, on peut toujours se tirer par un coup de pouce! Je veux dire par un distinguo."*
>
> Henri Poincaré, in (Poincaré 1917, p. 162).

The *Grundlagenkrise* which struck the mathematical and philosophical communities as a consequence of the paradoxes, showed one could not presuppose all pretheorethically plausible comprehension principles in set theory or semantics.

In the following some facts which relate to desiderata fulfilled by **£** will be expressed. The reader may compare with the desiderata of (Sheard 2003),, in (Halbach and Horsten 2002), and (Leitgeb 2007) or others, concerning theories on paradoxes. Some of the facts on desiderata follow from § 12.4, § 13.2 and § 13.3.

## 13.1 Facts on desiderata met fully, or partially, by £

13.1.1. Fact There are no type restricions imposed, and there is no language hierarchy.

13.1.2. Fact *Any* condition $A(v)$ on v corresponds with a set $\{v|A(v)\}$.

13.1.3. Fact (Degrees of compositionality)

Truth is compositional over $\Vdash^{M}$.

If $\Vdash^{M}$ A, $\Vdash$ B and $\Vdash$ $\neg$B, then $\Vdash$ A $\wedge$ B and $\Vdash$ A $\wedge$ $\neg$B.

But truth is *not* compositional over $\Vdash$, as for some A and B, $\Vdash$ A and $\Vdash$ B but not $\Vdash$ A $\wedge$ B.

13.1.4. Fact The contrasistency of **£** is fundamental for its intuitive *responsible* theory of comprehension, which serves well also in paradoxical contexts.

13.1.5. Fact Truth is a set, and so truth is as well a predicate. So it is a consequence from the alethic comprehension principle of § 5.3 that truth-paradoxes and set-paradoxes are treated in the same way in **£**.

13.1.6. Fact On account of results in § 13.3, **£** is *serene* in the sense that $\Vdash^{M}$ A only if classical logic does not prove $\neg$A, and if classical logic proves A then $\Vdash^{M}$ A. Moreover, $\Vdash$ B if B is a thesis of classical logic, and if $\Vdash$ B then classical logic does not prove $\neg$B.

13.1.7. Fact **£** is *unswerving* in the sense that if A is a paradoxical sentence, then **£** should have $\Vdash$ A or $\Vdash$ $\neg$A, and indeed it as a rule has both.

13.1.8. Fact **£** is ontologically parsimonious, for according to it there are precisely, and just $\aleph_0$ cardinal numbers, viz. $\aleph_0$ itself, that is the $\omega$ of Definition §11.0.1.3, and its members, which are the finite cardinal numbers.

13.1.9. Fact It follows immediately from Fact 13.1.8 that **£** satisfies the related, and very important desideratum of *being able to give an answer to the metaphysical question: How many objects are* *there?* The librationist answer is that there is $\omega$, i.e. a countable infinity of objects. Notice that one cannot provide an answer to that metaphysical question if one thinks in terms of classical set theory, for if one in such a framework thinks that there are $\kappa$ objects in the world, one must for Cantorian reasons also think that there is a larger cardinality $\lambda$ of objects in the world … .

13.1.10. Fact The outer veridical logic of **£** is the set of theses which encapsulate truth statements. $\Vdash^{M} \mathcal{T}\ulcorner A \urcorner \vee \neg\mathcal{T}\ulcorner A \urcorner$ is for example an instance of *The Law of Excluded Middle* of classical logic in the *outer veridical logic* of **£**. The inner veridical logic of **£** is the set of theses which are encapsulated by a truth predicate. $\Vdash^{M} \mathcal{T}\ulcorner A \vee \neg A \urcorner$ is for example an instance of *The Law of Excluded Middle* of classical logic in the *inner verdidical logic* of **£**. It is a consequence of Fact 13.1.6 that the inner and outer veridical logics of **£** are classical, as they should be.

13.1.11. Fact The variety of truth conditionals summed up in Exercise 13.2.2 has the consequence that the *outer veridical* and *inner veridical* logics of $\mathcal{T}$, see Fact 13.1.10, coincide in **£**, in the sense of Definitions 13.2.1.7 and 13.2.1.10.

13.1.12. Fact As related in § 7, **£** has novel inferential modes. The conjunction of these may seem to be an amputation of the classical inferential principle *modus ponens*. But they are in reality an extension of the classical inference rule modus ponens, as *the maxim mode* 7.2.7 serves all the purposes as modus ponens serves in classical logics, and all classical logical theses are maxims of **£**. The inference modes **£** has beyond the maxim mode helps engender novel paradoxical theses which are out of reach for classical logic.

13.1.13. Fact A naive desideratum is that **£** should obtain all truth-biconditionals, as in Definitions 13.2.1.1 and 13.2.1.7, with their weak counterparts, by means of the inference modes which **£** endorses, as per § 7. **£** compensates for the fact that the statements of Definitions 13.2.1.1, 13.2.1.2, 13.2.1.3 and 13.2.1.4 are not true with the truth of the statements of Definitions 13.2.1.5 and 13.2.1.6, and with the fact that the inferential modes exhibited in Definitions 13.2.1.7 and 13.2.1.10 can be used. A consequence of this is that *revenge paradoxicalities* are not a threat.

13.1.14. Fact By § 12.4, **£** is omega-complete, so it only allows standard interpretations.

## 13.2 The truth–conditionals

13.2.1. Definition

1 `Hale` material truth adequacy: $\Vdash^{M} \mathcal{T}\ulcorner A \urcorner \leftrightarrow A$

2 `Hale` material truthwards adequacy: $\Vdash^{M} \mathcal{T}\ulcorner A \urcorner \leftarrow A$

3 `Hale` material truthly adequacy: $\Vdash^{M} \mathcal{T}\ulcorner A \urcorner \rightarrow A$

4 `Weak` material truth adequacy: $\Vdash \mathcal{T}\ulcorner A \urcorner \leftrightarrow A$

5 `Weak` material truthwards adequacy: $\Vdash \mathcal{T}\ulcorner A \urcorner \leftarrow A$

6 `Weak` material truthly adequacy: $\Vdash \mathcal{T}\ulcorner A \urcorner \rightarrow A$

7 `Hale` formal truth adequacy: $\Vdash^{M} \mathcal{T}\ulcorner A \urcorner \Leftrightarrow \Vdash^{M} A$

8 `Hale` formal truthwards adequacy: $\Vdash^{M} \mathcal{T}\ulcorner A \urcorner \Leftarrow \Vdash^{M} A$

9 `Hale` formal truthly adequacy: $\Vdash^{M} \mathcal{T}\ulcorner A \urcorner \Rightarrow \Vdash^{M} A$

10 `Weak` formal truth adequacy: $\Vdash \mathcal{T}\ulcorner A \urcorner \Leftrightarrow \Vdash A$

11 `Weak` formal truthwards adequacy: $\Vdash \mathcal{T}\ulcorner A \urcorner \Leftarrow \Vdash A$

12 `Weak` formal truthly adequacy: $\Vdash \mathcal{T}\ulcorner A \urcorner \Rightarrow \Vdash A$

13.2.2. Exercise

Show that **£** obeys the formal and as well the weak material truthwards and truthly adequacies of Definition 13.2.1. The first four adequacies in the list fail on account of paradoxicalities

## 13.3 £ is *serene* and *extraclassical*, and it is not paraconsistent

Let $\mathbb{T}$ be a theory and $\Vdash$ its thesis indicator.

13.3.1. Definition

$$\mathbb{T} \text{ is } \textit{adjunctive} \text{ just if } \Vdash A \;\&\; \Vdash B \Rightarrow \Vdash A \wedge B.$$

13.3.2. Definition

$$\mathbb{T} \text{ is } \textit{dejunctive} \text{ just if } \Vdash A \wedge B \Rightarrow \Vdash A \;\&\; \Vdash B.$$

13.3.3. Definition

$$\mathbb{T} \text{ is } \textit{consistent} \text{ just if for no } p,\ \Vdash p \wedge \neg p.$$

13.3.4. Definition

$$\mathbb{T} \text{ is } \textit{contradictory} \text{ just if it is inconsistent.}$$

13.3.5. Definition

$$\mathbb{T} \text{ is } \textit{contrasistent} \text{ just if for some } p, \text{ both } \Vdash p \text{ and } \Vdash \neg p.$$

13.3.6. Definition

$$A \text{ is an } \textit{antithesis} \text{ of } \mathbb{T} \text{ just if } \neg A \text{ is a thesis of } \mathbb{T}.$$

13.3.7. Definition

$$\mathbb{S} \text{ is a moderation of } \mathbb{T} \text{ just if no thesis of } \mathbb{S} \text{ is an antithesis of } \mathbb{T}.$$

13.3.8. Definition

$$\mathbb{X} \text{ is an extension of } \mathbb{T} \text{ just if all theses of } \mathbb{T} \text{ are theses of } \mathbb{X}.$$

13.3.9. Definition

$$\mathbb{X} \text{ properly extends } \mathbb{T} \text{ just if } \mathbb{X} \text{ extends } \mathbb{T} \text{ and } \mathbb{T} \text{ does not extend } \mathbb{X}.$$

13.3.10. Definition

$\tau$ denotes *classical* logic.

13.3.11. Definition

$\mathbb{T}$ is *preclassical* just if some of $\tau$'s theses are not theses of $\mathbb{T}$

13.3.12. Definition

$\mathbb{T}$ is *progressive* just if it is a proper extension of $\tau$.

13.3.13. Definition

$\mathbb{T}$ is *moderate* just if it is a moderation of $\tau$.

13.3.14. Definition

$\mathbb{T}$ is *serene* just if it is progressive and moderate.

13.3.15. Definition

$\mathbb{T}$ is *extraclassical* just if it is serene and contrasistent.

13.3.16. Lemma

**£** is an extension of $\tau$.

*Proof.* Appeal to § 6.2. □

13.3.17. Lemma

$\tau$ is not an extension of **£**

*Proof.* **£** has paradoxical thesis $r \in r$, which is not theses of classical logic. □

13.3.18. Lemma

**£** is progressive.

*Proof.* **£** is a proper extension of $\tau$ given Fact 13.3.9 and Lemma 13.3.17. An appeal to Definition 13.3.12 suffices to finish the proof. □

13.3.19. Lemma

£ is moderate.

*Proof.* If $\tau$ proves $\neg A$, $\Vdash^{M} \neg A$ as **£** is progressive. By Theorem 5.6.3, $\tau$ proves $\neg A$ only if $\Vdash \neg A$ & $\nVdash A$. So a fortiori, if $\tau$ proves $\neg A$, $\nVdash A$. By contraposition, if $\Vdash A$ then $\tau$ does not prove $\neg A$. A is arbitrary, so no thesis of **£** is an antithesis of classical logic $\tau$. Consequently, **£** is a moderation of classical logic $\tau$. An appeal to Definition 13.3.13 finishes the proof. ☐

13.3.20. Theorem

**£** is serene.

*Proof.* By Definition 13.3.14, as **£** is progressive and moderate by Lemmas 13.3.18 & 13.3.17. ☐

13.3.21. Theorem: Paraconsistent theories are preclassical

*Proof.* Canonically, «A logic is paraconsistent iff its logical consequence relation ($\models$, either semantic or proof theoretic) is not explosive. Paraconsistency is a property of a consequence relation. The argument ex contradictione quodlibet (ECQ) is paraconsistently invalid: in general, it is not the case that $A, \neg A \models B$.»(Priest, Tanaka, et al. 2022, § 1)

As a consequence of the quotation: if the paraconsistent logic under consideration has modus ponens as an inference rule, $\nVdash A \to (\neg A \to B)$; so that logic Test1 is preclassical.

Let famous *Logic of Paradox*, of (Asenjo 1966), which was amended, and published broadly upon by Graham Priest, and others, does not validate modus ponens, but it does not validate

$$\nVdash A \to (\neg A \to B)$$

either, and so it is for the same reason, as for the unidentified paraconsisten modal logic of the previous paragraph without modus ponens, a preclassical logic. ☐

13.3.22. Theorem

Preclassical theories are not serene.

*Proof.* As preclassical theories cannot be progressive, i.e. proper extensions of classical logic, for there are theorems of classical logic which fail in preclassical theories. ☐

13.3.23. Theorem

**£** is not paraconsistent.

*Proof.* Invoke Theorems 13.3.20, 13.3.21 and 13.3.22. ☐

13.3.24. Exercise: $\mathbb{T}$ is contrasistent just if for some $p$, $\mathbb{T} \Vdash p$ and $\mathbb{T} \Vdash \neg p$

13.3.25. Exercise: If $\mathbb{T}$ is contrasistent and adjunctive then $\mathbb{T}$ is contradictory

13.3.26. Exercise: If $\mathbb{T}$ is contradictory and dejunctive then $\mathbb{T}$ is contrasistent

13.3.27. Exercise: If $\mathbb{T}$ is dejunctive and adjunctive, $\mathbb{T}$ is contradictory iff contrasistent

13.3.28. Exercise Adjunction is not a valid inference mode in **£**

13.3.29. REMARK THE PAIRS CONSISTENCY & COSISTENCY AND CONTRASISTENCY & CONTRADICTION CONFLATE IN CLASSICAL CONTEXTS, FOR CLASSICAL SYSTEMS ARE ADJUNCTIVE AND DEJUNCTIVE

13.3.30. REMARK WITH PROPER COMPREHENSION, MOST PARACONSISTENT THEORIES ARE NOT EVEN MODERATE, AS THEN SOME CONTRADICTION IS A THESIS

13.3.31. REMARKThe well-known non-adjunctive paraconsistent logic of (Jaskowski 1999) and (Jaskowski 1948), (Jaskowski 1948; Jaskowski 1999) is moderate even with liberal comprehension principles. But it is not *classic*, given Theorems 13.3.21 and 13.3.22.

## 13.4 Incompatability and complementarity

13.4.1. DEFINITION (INCOMPATABILITY)

Theses A and B of a consistent theory $\mathbb{T}$ are *incompatible* just if $\mathbb{T}$ proves A, B, and $\neg(A \wedge B)$.

13.4.2. THEOREM

£ has incompatible theses.

*Proof.* £ proves $R \in R$ and $R \notin R$. But £ is serene, given § 13.3. So the theses $R \in R$ and $R \notin R$ of £ are incompatible. Consequently, as £ is serene, £ proves $\neg(R \in R \wedge R \notin R)$. □

13.4.3. DEFINITION OF COMPLEMENTARY THESE:

A and $\neg$A in a theory $\mathbb{T}$ are complementary just if they are incompatible theses of $\mathbb{T}$.

13.4.4. DEFINITION OF COMPLEMENTARY THEORIES:

$\mathbb{T}$ is complementary just if it has complementary theses.

13.4.5. COROLLARY

£ is a complementary theory.

# 14 Substitutions and Diagonalization

*“… нельзя быть математиком, не будучи в то же время и поэтом в душе.”*

Sofia Kovalevskaya in a letter dated 1890 and reproduced in (Kovalevskaya 1951, p. 314). [a]

[a]“…it is impossible to be a mathematician without also being a poet in spirit.” (Kovalevskaya 1951, p. 314) translated to (Kovalevskaya 1978, p. 35).

When possible we adapt from (Smoryński 1977, pp. 836–837). Given Definitions 4.5.3.2 and 4.6.1, $\ulcorner \ddot{v}_0 \urcorner \mathrel{=\!=} \ulcorner \ddot{v} \urcorner = \{U\} \mathrel{=\!=} U^2$. Notice that on account of Definition 4.5.4.4, if $\ulcorner \ddot{v}_n \urcorner = U^m$ then $\ulcorner \ddot{v}_{n+1} \urcorner \mathrel{=\!=} \ulcorner \ddot{v}_n \frown \bullet \urcorner = U^{2m+1}$.

## 14.1 sub

Define recursively, with recourse to Definition 4.5.10 and other definitions above, presupposing that y is a universe counter, that a and b are terms and that A and B are formulas. $\mathcal{E}$ is an arbitrary formal expression, and $\mathcal{E}(\ddot{v}_i)$ is an arbitrary formal expression wherein the variable $\ddot{v}_i$ occurs freely:

(1) $\textbf{sub}(\ulcorner \ddot{v}_i \urcorner, i, y) \mathrel{=\!=} y$.

(2) $\textbf{sub}(\ulcorner \ddot{v}_i \urcorner, j, y) \mathrel{=\!=} \ulcorner \ddot{v}_i \urcorner$ if $i \neq j$.

(3) $\textbf{sub}(\ulcorner b \urcorner, i, y) \mathrel{=\!=} U^m \;\&\; \textbf{sub}(\ulcorner a \urcorner, i, y) \mathrel{=\!=} U^n$
$\Rightarrow \textbf{sub}(\ulcorner ba \urcorner, i, y) \mathrel{=\!=} U^{m \frown n}$.

(4) $\textbf{sub}(\ulcorner A \urcorner, i, y) \mathrel{=\!=} U^m \;\&\; \textbf{sub}(\ulcorner B \urcorner, i, y) \mathrel{=\!=} U^n \big)$
$\Rightarrow \textbf{sub}(\ulcorner \downarrow AB \urcorner, i, y) \mathrel{=\!=} U^{3 \frown m \frown n}$.

(5) $\textbf{sub}(\ulcorner \forall v_i A \urcorner, i, y) \mathrel{=\!=} \ulcorner \forall v_i A \urcorner$.

(6) $\textbf{sub}(\ulcorner \varsigma v_i A \urcorner, i, y) \mathrel{=\!=} \ulcorner \varsigma v_i A \urcorner$.

(7) $\textbf{sub}(\ulcorner \forall v_j A \urcorner, i, y) \mathrel{=\!=} \ulcorner \forall v_j \textbf{sub}(\ulcorner A \urcorner, i, y) \urcorner$ if $i \neq j$.

(8) $\textbf{sub}(\ulcorner \varsigma v_j A \urcorner, i, y) \mathrel{=\!=} \ulcorner \varsigma v_j \textbf{sub}(\ulcorner A \urcorner, i, y) \urcorner$ if $i \neq j$.

(9) $\textbf{sub}(\ulcorner \mathcal{E} \urcorner, i, y) \mathrel{=\!=} U^m \Rightarrow \textbf{sub}(\ulcorner \# \mathcal{E} \urcorner, i, y) \mathrel{=\!=} U^{6 \frown m}$.

(10) $\textbf{sub}(\ulcorner \mathcal{T} \urcorner, i, y)$ is $U^7$.

(11) $\textbf{sub}(\ulcorner \mathcal{F} \urcorner, i, y)$ is $U^8$.

14.1.1. The **sub**stitution Theorem

For any term t and any expression $\mathcal{E}(\ddot{v}_i)$, $\textbf{sub}(\ulcorner \mathcal{E}(\ddot{v}_i) \urcorner, i, \ulcorner t \urcorner) \mathrel{=\!=} \ulcorner \mathcal{E}(t) \urcorner$.

*Proof.* (I) The simplest expression $\mathcal{E}$ for which this applies is SUB 14.1.1, where we indeed have that $\textsc{sub}(\ulcorner \ddot{v}_i \urcorner, i, \ulcorner t \urcorner) = \ulcorner t \urcorner$. (II) SUB 14.1.2 does not have the correct form. (III) The atomic cases of SUB 14.1.3 give rise to the possible variants $\textsc{sub}(\ulcorner \ddot{v}_i a \urcorner, i, \ulcorner t \urcorner)$, $\textsc{sub}(\ulcorner b\ddot{v}_i \urcorner, i, \ulcorner t \urcorner)$ and $\textsc{sub}(\ulcorner \ddot{v}_i \ddot{v}_i \urcorner, i, \ulcorner t \urcorner)$. $\ulcorner t \urcorner$ is $U^t$, given Theorem 4.6.2.4. Suppose a has no occurrence of $\ddot{v}_i$, so that $\textsc{sub}(\ulcorner a \urcorner, i, \ulcorner t \urcorner)$ is $U^a$. It follows from 14.1.3 that $\textsc{sub}(\ulcorner \ddot{v}_i a \urcorner, i, \ulcorner t \urcorner)$ is $U^{t\frown a} = U^{ta} = \ulcorner ta \urcorner$; the other variants $\textsc{sub}(\ulcorner b\ddot{v}_i \urcorner, i, \ulcorner t \urcorner)$ and $\textsc{sub}(\ulcorner \ddot{v}_i \ddot{v}_i \urcorner, i, \ulcorner t \urcorner)$ have corresponding verifying consequences. So these considerations establish that if the substitution theorem holds for $\ulcorner a \urcorner$ and for $\ulcorner b \urcorner$, then it as well holds for $\ulcorner ba \urcorner$. (IV) The formula $\downarrow AB$ gives no information on possible occurrences of $\ddot{v}_i$, so this case is not to be considered; however, if we suppose that we have $\downarrow A(\ddot{v}_i)B$, we know that $\textsc{sub}(\ulcorner A(\ddot{v}_i) \urcorner, i, \ulcorner t \urcorner) = \ulcorner A(t) \urcorner = U^{A(t)}$ so that $\textsc{sub}(\ulcorner \downarrow A(\ddot{v}_i)B \urcorner, i, \ulcorner t \urcorner) = \ulcorner \downarrow A(t)B \urcorner$. Such considerations show that the substitution theorem holds for $\ulcorner \downarrow AB \urcorner$ if it holds for $\ulcorner A \urcorner$ and for $\ulcorner B \urcorner$. (V) (5) does not have $\ddot{v}_i$ free. (VI) (6) does not have $\ddot{v}_i$ free. (VII) The condition for (7) $\textsc{sub}(\ulcorner \forall v_j A \urcorner, i, y)$ with $i \neq j$ shows that if $\ddot{v}_i$ is free in A, by the induction hypothesis, $\textsc{sub}(\ulcorner \forall v_j A(\ddot{v}_i) \urcorner, i, \ulcorner t \urcorner)$ is $\ulcorner \forall v_j A(t) \urcorner$. (VIII) The condition for (8) $\textsc{sub}(\ulcorner \varsigma v_j A \urcorner, i, y)$ with $i \neq j$ shows that if $\ddot{v}_i$ is free in A, by the induction hypothesis, $\textsc{sub}(\ulcorner \varsigma v_j A(\ddot{v}_i) \urcorner, i, \ulcorner t \urcorner)$ is $\ulcorner \varsigma v_j A(t) \urcorner$. (IX–XI) In the cases (9–11) $\ddot{v}_i$ does not have free occurrences. □

## 14.2 Sb

14.2.1. DEFINITION

$\mathbf{Sb}(x, y) = \textbf{sub}(x, i, y)$ if $i = \mu j(j < x$ such that $v_i$ occurs free in x), and
$\mathbf{Sb}(x, y) = x$ if there is no such i.

## 14.3 Diagonalization

14.3.1. THEOREM: THE DIAGONAL LEMMA.

If $A(v_i)$ has precisely the indicated variable free, then there is a sentence D such that

$$\Vdash^M D \leftrightarrow A\ulcorner D \urcorner.$$

*Proof.* Let m be $\ulcorner A(\mathbf{Sb}(x, x)) \urcorner$ and D be $A(\mathbf{Sb}(m, m))$. We have:

$$\Vdash^M D \leftrightarrow A(\mathbf{Sb}(m, m)) \leftrightarrow A(\mathbf{Sb}(\ulcorner A(\mathbf{Sb}(x, x)) \urcorner, m)) \leftrightarrow A(\ulcorner A(\mathbf{Sb}(m, m)) \urcorner)) \leftrightarrow A(\ulcorner D \urcorner).$$

□

# 15 Liar theses and their origins

*“ξεῖνε, Φιλίτας εἰμί· λόγων ὁ ψευδόμενός με*
*ὤλεσε καὶ νυκτῶν φροντίδες ἑσπέριοι.”*

---

The epitaph on Philetas of Coz, according to Athenaeus of Naucratis in *Deipnosophistae* 9.401 e.[a]

[a]“Philetas of Cos am I,
’Twas the Liar who made me die,
And the bad nights caused thereby.”

---

English rendering by St. George William Joseph Stock in his *Stoicism*, London: Archibald Constable, 1908, p. 36.

As the epigraph above helps testifying, the Liar paradoxes have a distinguished history back to antiquity in philosophy, and have continued to occupy a central place in philosophical logic and the philosophy of language. Needless to say, but the Liar paradoxes are related to the set theoretical paradoxes - which are a main focus of this treatise.

The first subsection below contains a theorem and an observation on Liar theses, and the next ones give two distinct origins for Liar theses.

## 15.1 A theorem and an observation on Liar theses

15.1.1. Theorem (The Liar theses)

If a sentence L is a *Liar* paradox in the sense that $\Vdash^{M} \mathrm{L} \leftrightarrow \neg\mathcal{T}\ulcorner \mathrm{L} \urcorner$, then paradoxical Liar theses as the following are derivable: $\Vdash \mathrm{L}$, $\Vdash \neg \mathrm{L}$, $\Vdash \mathcal{T}\ulcorner \mathrm{L} \urcorner$, $\Vdash \mathcal{T}\ulcorner \neg\mathrm{L} \urcorner$, $\Vdash \neg\mathcal{T}\ulcorner \neg\mathrm{L} \urcorner$ and $\Vdash \neg\mathcal{T}\ulcorner \mathrm{L} \urcorner$.

*Proof.* It suffices to do $\neg\mathcal{T}\ulcorner \mathrm{L} \urcorner$, and leave the others as exercises: As $\Vdash \mathcal{T}\ulcorner \mathrm{L} \urcorner \rightarrow \mathrm{L}$, given Postulate 6.5.1.1, substitution yields $\Vdash \mathcal{T}\ulcorner \mathrm{L} \urcorner \rightarrow \neg\mathcal{T}\ulcorner \mathrm{L} \urcorner$, and so logic entails $\Vdash \neg\mathcal{T}\ulcorner \mathrm{L} \urcorner$. □

15.1.2. Observation

Each element in variety [L, $\mathcal{T}\ulcorner \mathrm{L} \urcorner$, $\neg\mathcal{T}\ulcorner \neg\mathrm{L} \urcorner$] is incompatible with any member of [¬L, $\mathcal{T}\ulcorner \neg\mathrm{L} \urcorner$, $\neg\mathcal{T}\ulcorner \mathrm{L} \urcorner$] in **£**, and vice versa. Moreover, each element in variety [L, $\mathcal{T}\ulcorner \mathrm{L} \urcorner$, $\neg\mathcal{T}\ulcorner \neg\mathrm{L} \urcorner$] is complementary to precisely one member of [¬L, $\mathcal{T}\ulcorner \neg\mathrm{L} \urcorner$, $\neg\mathcal{T}\ulcorner \mathrm{L} \urcorner$] in **£**, and vice versa.

## 15.2 Russell’s set and Alethic Comprehension as origin

According to (Ramsey 1925, p. 20), there is an essential difference between *syntactic paradoxes* which ″involve only logical or mathematical terms such as class and number″, and *semantic paradoxes*, which ″… are not purely logical, and cannot be stated in logical terms alone; for they all contain some reference to thought, language, or symbolism″. Ramsey offered no precise distinction between ″logical or mathematical terms″ and ″terms which contain some reference to thought, language, or symbolism″, but he considered Russell’s a canonical representative of syntactic paradoxes, and the Liar an archetypical semantic paradox.

In (A. A. Fraenkel and Bar-Hillel 1958, p. 5), the authors adjudged:

> "Since (Ramsey 1925) it has become customary to distinguish between logical and semantic (sometimes also called syntactic or epistemological) antinomies."

It is here argued, to the contrary, that one should take paradoxes, as the Liar-paradox, to be so inextricably intertwined with set theoretical paradoxes as to not consider them to be different in kinds.

Others reached the same conclusion, but on the basis of considerations different from the ones adduced further below:

(Scott 1974) argued that Zermelo's set theory axioms were justified by type theoretic reasoning:

> "The truth is that there is only one way of avoiding the paradoxes: namely, the use of some form of the theory of types. That was at the basis of both Russell's and Zermelo's intuitions. Indeed the best way to regard Zermelo's theory is as a simplification and extension of Russell's. (We mean Russell's *simple* theory of types, of course.) The simplification was to make the types *cumulative*." (Scott 1974, p. 208)

Alonzo Church, who, incidentally, was my teacher in a graduate seminary in logic, with an oral exam, in the spring of 1989, at UCLA, virtually equated Russell's theory of types and Alfred Tarski's resolution of the Liar paradox, in (Church 1976), as he stated:

> "In the light of this it seems justified to say that Russell's resolution of the semantical antinomies is not a different one than Tarski's but is a special case of it."(Church 1976, p. 756)

The interest of Scott's and Church's points of view, for our purposes here, is that they take Tarski's resolution of what Ramsey would consider a *semantic* paradoxes to be the same as Russell's or Zermelo's resolution of a set theoretic paradox which Ramsey would consider to be a *syntactical* paradox.

In **£** there are bridge principles, as for example per Theorem 15.2.1 and Definition 15.2.4, between given, supposedly syntactical paradoxes, and supposedly semantical paradoxes.

15.2.1. Theorem

There is a *Liar sentence* L obeying $\Vdash^{\mathrm{M}} \mathrm{L} \leftrightarrow \neg\mathcal{T}\ulcorner \mathrm{L} \urcorner$.

*Proof.* Recall Definition 15.2.1 of r as $\{x|x \notin x\}$. By alethic comprehension,

$$\Vdash^{\mathrm{M}} r \in r \leftrightarrow \mathcal{T}\ulcorner r \notin r) \urcorner. \tag{15.2.2}$$

By negating both sides of the biconditional in 15.2.2, we get

$$\Vdash^{\mathrm{M}} r \notin r \leftrightarrow \neg\mathcal{T}\ulcorner r \notin r) \urcorner. \tag{15.2.3}$$

15.2.4. Definition

$$\mathrm{L} =\!= r \notin r,$$

Substituting with L of Definition 15.2.4 in equation 15.2.3 gives the more canonical, and announced, equation for the Liar sentence:

$$\Vdash^{\mathrm{M}} \mathrm{L} \leftrightarrow \neg\mathcal{T} \ulcorner \mathrm{L} \urcorner . \qquad \square \quad (15.2.5)$$

The Liar Paradox 15.2.5 is resolved semantically and philosophically as Russell's paradox.

Liar sentences, and variants, with provenances from classical Greek philosophy, should be taken as given by maxims of Theorems as 15.2.1.

15.2.6. Observation

Notice that there are other definitions than Definition 15.2.4, e.g. some introduced by means of the Diagonal Lemma, which support the biconditional 15.2.5 of the Liar Paradox, and such instances will as well as well support Theorem 15.1.1.

## 15.3 The Diagonal Lemma as origin

Given §15, it suffices to establish

15.3.1. Theorem

There is a sentence $\Lambda$ obtained by means of the function **Sb** which is such that

$$\Vdash^{\mathrm{M}} \Lambda \leftrightarrow \neg\mathcal{T} \ulcorner \Lambda \urcorner .$$

*Proof.* Let m in the proof of Theorem 14.3.1 be $\ulcorner \neg\mathcal{T}(\mathbf{Sb}(\mathrm{x},\mathrm{x})) \urcorner$ and $\Lambda$ be $\neg\mathcal{T}(\mathbf{Sb}(\mathrm{m},\mathrm{m}))$. $\square$

15.3.2. Remark on Theorem 15.3:

If one were to suppose that $\Vdash^{\mathrm{M}} \Lambda$ it would follow that $\Vdash^{\mathrm{M}} \mathcal{T} \ulcorner \Lambda \urcorner$ from Theorem 16.1.2.1, and $\Vdash^{\mathrm{M}} \neg\mathcal{T} \ulcorner \Lambda \urcorner$ from Corollary 15.3. If one were to suppose that $\Vdash^{\mathrm{M}} \neg\Lambda$ it would follow that $\Vdash^{\mathrm{M}} \mathcal{T} \ulcorner \neg\Lambda \urcorner$ from Theorem 16.1.2.1, and $\Vdash^{\mathrm{M}} \mathcal{T} \ulcorner \Lambda \urcorner$ from Corollary 15.3, so that one on account of Postulate 6.3.2.2 would have that $\Vdash^{\mathrm{M}} \mathcal{T} \ulcorner \neg\Lambda \urcorner$ and $\Vdash^{\mathrm{M}} \neg\mathcal{T} \ulcorner \neg\Lambda \urcorner$. But **£** is maximally adjunctive in the sense that

$$\Vdash^{\mathrm{M}} \mathrm{A} \;\&\; \Vdash^{\mathrm{M}} \mathrm{B} \Rightarrow \Vdash^{\mathrm{M}} \mathrm{A} \wedge \mathrm{B},$$

so that $\Vdash^{\mathrm{M}} \mathcal{T} \ulcorner \neg\Lambda \urcorner \wedge \neg\mathcal{T} \ulcorner \neg\Lambda \urcorner$ would follow from supposing $\Vdash^{\mathrm{M}} \neg\Lambda$ or from supposing $\Vdash^{\mathrm{M}} \Lambda$. But $\nVdash^{\mathrm{M}} \mathcal{T} \ulcorner \neg\Lambda \urcorner \wedge \neg\mathcal{T} \ulcorner \neg\Lambda \urcorner$, so consequently $\nVdash^{\mathrm{M}} \neg\Lambda$ and $\nVdash^{\mathrm{M}} \Lambda$. Given Theorem 15.1.1, however, which holds as well for the Liar sentence $\Lambda$ of Corollary 15.3, both $\Vdash \neg\Lambda$ and $\Vdash \Lambda$. So the sentence $\neg\Lambda$, which, given the librationist derivability conditions, is reminiscent of Gödel's incompleteness sentence, is *maximally* incomplete, in the sense that neither $\neg\Lambda$ nor $\Lambda$ is a maxim. But $\Lambda$ is nevertheless a paradoxical thesis, and an antithesis. This suggests that **£** is not incomplete for Gödelian reasons.

# 16 How provabilitylike is truth? How truthlike is provability?

> *"Kurt Gödel, when colleague John Bachall presented himself as a physicist at an Institute for Advanced Studies faculty dinner: "I don't believe in natural science.""*
>
> Ed Regis, in (Regis 1988, p. 58).

The similarity and differences between the canonic provability conditions (or Hilbert-Bernays-Löb derivability conditions) and the librationist deducibility provisions are discussed in §16.1. The adequacy of the librationist deducibility provisions are discussed briefly in Remarks 16.1.4, and shed light upon by Theorem 16.1.5.

## 16.1 Proof and canonic provability conditions compared with Librationist truth and the deducibility provisions for £

In connection with proofs of incompleteness for Peano arithmetic, and similar theories, one often appeals to the notion of *derivability* or *provability conditions*; to express the canonical ones, let «Bew» be Gödel's provability predicate, «PA» abbreviate Peano Arithmetic and $\ulcorner\urcorner$ stand for the classical Gödel codes:

16.1.1. Reminder (The *canonical* – i.e. Hilbert-Bernays-Löb – derivability conditions)

1 $\vdash^{PA} A \Rightarrow \vdash^{PA} \mathtt{Bew}(\ulcorner A \urcorner)$.

2 $\vdash^{PA} \mathtt{Bew}(\ulcorner A \urcorner) \rightarrow \mathtt{Bew}(\ulcorner \mathtt{Bew}(\ulcorner A \urcorner) \urcorner)$.

3 $\vdash^{PA} \mathtt{Bew}(\ulcorner A \rightarrow B \urcorner) \rightarrow (\mathtt{Bew}(\ulcorner A \urcorner \rightarrow \mathtt{Bew}(\ulcorner B \urcorner))$.

Under the translation from $\vdash^{PA}$ to $\Vdash^{M}$ and $\mathtt{Bew}(\ulcorner A \urcorner)$ to $\mathcal{T} \ulcorner A \urcorner$, the librationist thruth predicate $\mathcal{T}$ of **£** obeys the following slightly weaker conditions:

16.1.2. Theorem (The librationist deducibility provisions)

1 $\Vdash^{M} A \Rightarrow \Vdash^{M} \mathcal{T} \ulcorner A \urcorner$.

2 $\Vdash^{O} \mathcal{T} \ulcorner A \urcorner \rightarrow \mathcal{T} \ulcorner \mathcal{T} \ulcorner A \urcorner \urcorner$.

3 $\Vdash^{M} \mathcal{T} \ulcorner A \rightarrow B \urcorner \rightarrow (\mathcal{T} \ulcorner A \urcorner \rightarrow \mathcal{T} \ulcorner B \urcorner)$.

*Proof.* Theorem 16.1.2.1 holds on account of Postulate 7.1.2.1. Theorem 16.1.2.2 is one half of Postulate 6.5.1.2. Theorem 16.1.2.3 is Truth maxim 6.3.2.1. □

16.1.3. Remark on Theorem 16.1.2.2:

It is crucial that $\Vdash$ has superscript $^{O}$, as opposed to $^{M}$; cfr. §5.5 and §6.5.

16.1.4. Remarks on the adequacy of the librationist deducibility provisions

Despite the fact that Theorem 16.1.2.2, with §5.5 and its Definition 5.5.4 and its Resolve of appropriateness 5.5.5 entail that $\Vdash \mathcal{T}\ulcorner \mathrm{D} \urcorner \wedge \neg\mathcal{T}\ulcorner\mathcal{T}\ulcorner \mathrm{D}\urcorner\urcorner$ for *some* sentence D; still and still, for *any* sentence A, $\Vdash \mathcal{T}\ulcorner \mathrm{A}\urcorner \rightarrow \mathcal{T}\ulcorner\mathcal{T}\ulcorner \mathrm{A}\urcorner\urcorner$. Furthermore, notice that on account of Theorem 16.1.2.1, **£** retains $\Vdash^{\mathrm{M}} \mathcal{T}\ulcorner \mathrm{A}\urcorner \Rightarrow \Vdash^{\mathrm{M}} \mathcal{T}\ulcorner\mathcal{T}\ulcorner \mathrm{A}\urcorner\urcorner$, which in a certain sense weighs up for the paradoxical fact that $\Vdash \mathcal{T}\ulcorner \mathrm{D}\urcorner \wedge \neg\mathcal{T}\ulcorner\mathcal{T}\ulcorner \mathrm{D}\urcorner\urcorner$ for some sentence D. These remarks naturally amount to the following

16.1.5. Theorem:

> Under the translation presupposed, the librationist deducibility provisions
> coincide with the canonical provability conditions for orthodox sentences.

*Proof.* This follows from the fact that $\Vdash^{\mathrm{M}} \mathcal{T}\ulcorner \mathrm{A}\urcorner \rightarrow \mathcal{T}\ulcorner\mathcal{T}\ulcorner \mathrm{A}\urcorner\urcorner$ if A is orthodox. □

## 16.2 Disparity between Gödel's coding and librationist coding, and the effect upon prerequsites for the converses of Theorems 16.1.1.1 and 16.1.2.1

The coding in (Gödel 1931) is by recourse to arithmetic, and an effect is that a requirement for the converse of Theorems 16.1.1.1 to hold is that the presupposed arithmetic is $\omega$-consistent; for this, see In **£** the coding is introduced already in the formation rules, and not via an arithmetical construction; the converse of Theorem 16.1.2.1 holds on account of Postulate 7.1.2.2, and the latter holds on account of the semantics for **£** without any appeals to assumptions like $\omega$-consistency.

## 16.3 Unsuprisingly, but even the diagonal lemma with the librationist deducibility provisions do not entail Löb's rule or Löb's formula

In §8, the proof of Theorem 8.0.4 established that Löb's rule for truth in **£** fails, and the proof of Theorem 8.0.11 showed that Löb's formula fails in **£**. Nevertheless, it is of interest to see how these fact also follow from the fact that the librationist deducibility provisions are slightly weaker than the canonical Hilbert-Bernays-Löb derivability conditions.

The argument from the Diagonal Lemma to Löb's theorem fails in **£** on account of the presence of paradoxical sentences. Given two sentences A and B such that $\Vdash^{\mathrm{M}}$ A and $\Vdash^{\mathrm{M}} \neg$B, the fixpoints $\Vdash^{\mathrm{M}} \mathrm{C^A} \leftrightarrow (\mathcal{T}\ulcorner \mathrm{C^A}\urcorner \rightarrow \mathrm{A})$ and $\Vdash^{\mathrm{M}} \mathrm{C^B} \leftrightarrow (\mathcal{T}\ulcorner \mathrm{C^B}\urcorner \rightarrow \mathrm{B})$ are evaluable as follows: Firstly, $\Vdash^{\mathrm{M}} \mathrm{C^A}$ as $\Vdash^{\mathrm{M}}$ A entails $\Vdash^{\mathrm{M}} (\mathcal{T}\ulcorner \mathrm{C^A}\urcorner \rightarrow \mathrm{A})$. In the second case we have: (I) $\Vdash^{\mathrm{M}} \mathrm{C^B} \leftarrow \neg\mathcal{T}\ulcorner \mathrm{C^B}\urcorner$ as $\Vdash^{\mathrm{M}} \mathrm{C^B} \leftrightarrow (\mathcal{T}\ulcorner \mathrm{C^B}\urcorner \rightarrow \mathrm{B})$, and (II) $\Vdash^{\mathrm{M}} \mathrm{C^B} \rightarrow \neg\mathcal{T}\ulcorner \mathrm{C^B}\urcorner$ as $\Vdash^{\mathrm{M}} \neg$B and $\Vdash^{\mathrm{M}} \mathrm{C^B} \leftrightarrow (\mathcal{T}\ulcorner \mathrm{C^B}\urcorner \rightarrow \mathrm{B})$. So (III) $\Vdash^{\mathrm{M}} \mathrm{C^B} \leftrightarrow \neg\mathcal{T}\ulcorner \mathrm{C^B}\urcorner$.

Let us thereupon consider G. Boolos' rendition in (George Boolos 1993, p. 56) of the argument, from the fixpoint there written in a simplified style as

(1) $\mathrm{PA} \vdash \mathrm{I} \leftrightarrow (\mathrm{PI} \rightarrow \mathrm{S})$,

to Löb's theorem.

Compare first with fixpoint $\Vdash^{\mathrm{M}} \mathrm{C^A} \leftrightarrow (\mathcal{T}\ulcorner \mathrm{C^A}\urcorner \rightarrow \mathrm{A})$ plus $\Vdash^{\mathrm{M}}$ A: All the steps in the argument chain from the just stated (1) until the final step

(14) $\mathrm{PA} \vdash \mathrm{S}$

on (George Boolos 1993, p. 57) can be retraced, which is no surprise as we presupposed that $\Vdash^{M}$ A.

Compare next with fixpoint $\Vdash^{M} C^{B} \leftrightarrow (\mathcal{T}\ulcorner C^{B}\urcorner \rightarrow B)$ plus $\Vdash^{M} \neg B$: In this case we can retrace all the steps from (1–6) on (George Boolos 1993, p. 56) until (7) on (George Boolos 1993, p. 56) with $\Vdash^{M}$ instead of PA $\vdash$, $C^{B}$ instead of I, $\mathcal{T}\ulcorner C^{B}\urcorner$ instead of PI and B instead of S:

(7£) $\Vdash^{M} \mathcal{T}\ulcorner C^{B}\urcorner \rightarrow (\mathcal{T}\ulcorner\mathcal{T}\ulcorner C^{B}\urcorner\urcorner \rightarrow B)$.

But we do not have the correspondent of (8) in (George Boolos 1993, p. 56), as on the contrary

($\not 8$**£**$^{M}$) $\nVdash^{M} \mathcal{T}\ulcorner C^{B}\urcorner \rightarrow \mathcal{T}\ulcorner\mathcal{T}\ulcorner C^{B}\urcorner\urcorner$.

Although **£** doesn't have the maxim $\Vdash^{M} \mathcal{T}\ulcorner C^{B}\urcorner \rightarrow \mathcal{T}\ulcorner\mathcal{T}\ulcorner C^{B}\urcorner\urcorner$, it has the optimal theorem $\Vdash^{O} \mathcal{T}\ulcorner C^{B}\urcorner \rightarrow \mathcal{T}\ulcorner\mathcal{T}\ulcorner C^{B}\urcorner\urcorner$. But this permits that in cases also $\Vdash \mathcal{T}\ulcorner C^{B}\urcorner \,\&\, \neg\mathcal{T}\ulcorner\mathcal{T}\ulcorner C^{B}\urcorner\urcorner$ is a theorem.

So:

(8£$^{O}$) $\Vdash \mathcal{T}\ulcorner C^{B}\urcorner \rightarrow \mathcal{T}\ulcorner\mathcal{T}\ulcorner C^{B}\urcorner\urcorner$.

Given (7£) and (8£$^{O}$) and Postulate 7.2.6, it follows that $\Vdash$ B. But it was assumed that $\VDash^{M} \neg B$. As (7£) cannot be negotiated with, we must conclude that $\Vdash \mathcal{T}\ulcorner C^{B}\urcorner \wedge \neg\mathcal{T}\ulcorner\mathcal{T}\ulcorner C^{B}\urcorner\urcorner$.

Consequently, also the correspondent of (9) on (George Boolos 1993, p. 56) fails, as instead

($\not 9$£) $\nVdash^{M} \mathcal{T}\ulcorner C^{B}\urcorner \rightarrow \mathcal{T}\ulcorner B\urcorner$.

## 16.4 The librationist deducibility provisions block Gödelian incompleteness

Let us stress again that **£** is incomplete; for example, neither $s \in s$ nor $s \notin a$ is a theorem of **£** if $s = \{x|x \in x\}$. Importantly, the failure of the Löb-rule, and the associated stronger Löb formula, in £ has the consequence that the incompleteness pointed to in the previous sentence of this section may obtain in £. For $\Vdash^{M} s \in s \leftrightarrow \mathcal{T}\ulcorner s \in s\urcorner$, and had the Löb rule reigned in **£** then we would have had $\Vdash^{M} s \in s$.

As the Löb-formula fails for the truth predicate of **£**, the fact that it is derivable under the canonical derivability conditions strengthens the plausibility of a suggestion in Remarks 16.1.4 to the effect that the slightly librationist derivability provisions are *more adequate* for our set and truth theoretic framework than are the canonical derivability conditions.

To drive the last point home, consider that (Löb 1955) established the result now known as Löb's *theorem*, that A is provable in arithmetic if $Bew(\ulcorner A\urcorner) \rightarrow A$ is, where Bew is Gödel's provability predicate.Theorem 2.1 of (Macintyre and Simmons 1973, pp. 170–171) may be the first location which has a complete proof of the stronger result that the Löb *formula*

$$Bew(\ulcorner Bew(\ulcorner A\urcorner) \rightarrow A\urcorner) \rightarrow Bew(\ulcorner A\urcorner)$$

is a theorem of arithmetic.[3] But $\Vdash \mathcal{T}\ulcorner A\urcorner$ is not entailed by $\Vdash \mathcal{T}\ulcorner\mathcal{T}\ulcorner A\urcorner \rightarrow A\urcorner$ in **£**, and the Löb-formula $\mathcal{T}\ulcorner\mathcal{T}\ulcorner A\urcorner \rightarrow A\urcorner \rightarrow \mathcal{T}\ulcorner A\urcorner$ is not a theorem of **£**.

[3] Löb noticed that the Löb-formula $Bew(\ulcorner Bew(\ulcorner A\urcorner) \rightarrow A\urcorner) \rightarrow Bew(\ulcorner A\urcorner)$ is derivable, however, or it was an improbable oversight if he didn't: For he established Löb's theorem which states that A is derivable if $Bew(\ulcorner A\urcorner) \rightarrow A$

## 16.5 More notes on provability and incompleteness in £

Proofs of Gödel's incompleteness theorem typically appeal to an assumption that the theory is not *contrasistent*, in the sense of Definition 13.3.5, and in the classical frameworks the presupposed failure of contrasistency is equivalent with the failure of *consistency*, as in Definition 13.3.3. Arguably, the consistent and contrasistent set theory **£** finesses Gödelian incompleteness, but that does not mean that **£** is semantically complete.

Let for comparison «Bew» be Gödel's provability predicate and $\ulcorner\urcorner$ stand for Gödel codes.

Notice Theorem 16.1.2.2 has $\Vdash^{O}$ and not $\Vdash^{M}$. Despite the fact that Theorem 16.1.2.2 and that the nature of proper *optima* given §5.5 entail that $\Vdash \mathcal{T}\ulcorner D\urcorner \wedge \neg\mathcal{T}\ulcorner\mathcal{T}\ulcorner D\urcorner\urcorner$ for some sentence D, still and still, as well for *any* sentence A, $\Vdash \mathcal{T}\ulcorner A\urcorner \rightarrow \mathcal{T}\ulcorner\mathcal{T}\ulcorner A\urcorner\urcorner$. Moreover, notice that on account of Theorem 16.1.2.1, **£** retains $\Vdash^{M} \mathcal{T}\ulcorner A\urcorner \Rightarrow \Vdash^{M} \mathcal{T}\ulcorner\mathcal{T}\ulcorner A\urcorner\urcorner$, which weighs up for the fact that $\Vdash \mathcal{T}\ulcorner D\urcorner \wedge \neg\mathcal{T}\ulcorner\mathcal{T}\ulcorner D\urcorner\urcorner$ for some sentence D. Notice also that if A is not paradoxical, then $\Vdash^{M}$ $\mathcal{T}\ulcorner A\urcorner \rightarrow \mathcal{T}\ulcorner\mathcal{T}\ulcorner A\urcorner\urcorner$ and so for sentences which are not paradoxical the librationist derivability conditions correspond precisely with the canonical provability conditions.

Importantly, as shown in §14, the Löb-rule and Löb formula do not hold for truth in £. As a consequence the incompleteness pointed to in the last paragraph of this section may obtain in £, although $\Vdash^{M} s \in s \leftrightarrow \mathcal{T}\ulcorner s \in s\urcorner$ as well as $\Vdash^{M} s \notin s \leftrightarrow \mathcal{T}\ulcorner s \notin s\urcorner$. Arguably, the librationist deducibility conditions for $\mathcal{T}$ are on account of the Löb failure more adequate than the canonical provability conditions associated with Gödel's `Bew` predicate used in the proof of the Incompleteness Theorem.

Given all of this one may take $\Vdash^{M} \mathcal{T}\ulcorner A\urcorner$ to not only express that A is a true maxim, but as well that A is *provable* as a maxim. An actual maxim of **£** as $\Vdash^{M} \neg x\mathcal{T}\ulcorner A \wedge \neg A\urcorner$ may therefore be taken to expresses that **£** believes that it is not trivial, and that **£** believes that it does not prove a contradiction.

The author does not know that there is a sentence A such that $\nVdash^{M}$ A in all the librationist theories under consideration in this manuscript, and such that we *should* want that $\Vdash^{M}$ A, nor that there is a sentence B such that $\nVdash$ B and such that we should want that $\Vdash$ B. Certainly, if C is the statement that there is a certain inaccessible cardinal larger or equal to the first hypothetized *1-inaccesible* cardinal, it will be the case, with the assumptions made, that even $\nVdash^{M} C^{\mathbf{V}}$, where $\mathbf{V}$ is as in § 25, and the notation $C^{\mathbf{V}}$ as in Definition 4.5.21. It is not obvious to the author, however, at this point, that we should want $\Vdash^{M} C^{\mathbf{V}}$.

But for the record, given § 25: if D is the statement that there is an *inaccesible* cardinal, and that there for any *inaccessible* cardinal is a larger *inaccessible* cardinal, then $\Vdash^{M} D^{\mathbf{V}}$. However, those inaccessible postulated there count as *0-inaccessible* cardinals, and the theory **NBG+TA** which is modelled has its least standard model in $\mathbf{V}_{\zeta}$ if $\zeta$ is the first *1-inaccessible* cardinal, i.e. the first regular limit of *inaccessible* cardinals.

---

is derivable. But (Löb, 1955) p. 115 mandated a derivability condition II which states that if formula T is derivable from the formula S, formula $\mathrm{Bew}(\ulcorner S\urcorner) \rightarrow \mathrm{Bew}(\ulcorner T\urcorner)$ is a theorem. As A is derivable from $\mathrm{Bew}(\ulcorner A\urcorner) \rightarrow A$, by Löb' Theorem, the Löb-formula $\mathrm{Bew}(\ulcorner \mathrm{Bew}(\ulcorner A\urcorner) \rightarrow A\urcorner) \rightarrow \mathrm{Bew}(\ulcorner A\urcorner)$ is a theorem.

At a more pedestrian level, when compared with inaccessible cardinals, if $s = \{x|x \in x\}$ neither $\Vdash s \in s$ nor $\Vdash s \notin s$. A morale is that it is not always an incompleteness that neither a sentence nor its negation is a thesis of **£**, for we have no reasons to presuppose that the sentence $s \in s$ is true, or that the sentence $s \notin s$ is true.

## 16.6 Responsible naiveté without revenge

The revenge problem is avoided as **£** is unswerving, in the sense of Fact 13.1.7, and as it has complementary theses in the sense of Definition 13.4.3.

Consider the Liar sentence L of equation 15.2.5. If $\VDash L$, it follows that $\VDash \neg\mathcal{T}\ulcorner L\urcorner$ via the equation 15.2.5. However, it as well follows, and is, as we shall presuppose, *responded*, that $\VDash \mathcal{T}\ulcorner L\urcorner$, from $\VDash L$ and Postulate 7.1.1. Moreover, given $\VDash \neg\mathcal{T}\ulcorner L\urcorner$, $\VDash \neg L$ is responded, and the response happens again on account of Postulate 7.1.1.

It is not a desirable option to prefer $\nVDash L$ and $\nVDash \neg L$, for it stands to reason that either the Liar sentence is true, or the Liar sentence is not true.

From the responsible extraclassical point of view $\VDash \neg L$ and $\VDash \neg\mathcal{T}\ulcorner L\urcorner$ state things as they are, as do $\VDash \mathcal{T}\ulcorner L\urcorner$ and $\VDash L$, and the sentences are taken to respond to reality appropriately.

## 16.7 $\ulcorner A\urcorner$ is true just if A states the truth

The following perspective upon the semantics might be useful for some purposes.

16.7.1. Definition The closure ordinal $\varpi$ is *the truth*.

16.7.2. Definition The *way* of sentence A is $[\delta : \delta \preceq \varpi \;\&\; (\Xi, \delta) \Vdash A]$.

16.7.3. Definition A *states* the supremum of its way.

16.7.4. Definition A *expresses* its way.

16.7.5. Definition $\ulcorner A\urcorner$ is true just if A states the truth.

16.7.6. Definition $\ulcorner A\urcorner$ is false just if $\ulcorner\neg A\urcorner$ is true.

16.7.7. Definition The way of $\neg A$ is the truth minus the way of A.

16.7.8. Definition The way of $A \wedge B$ is the way of intersected with the way of B.

Here the sentence $\ulcorner A\urcorner$ *is true* is to be interpreted as $\Vdash \mathcal{T}\ulcorner A\urcorner$, and the sentence A *states the truth* as equivalent with $\Vdash A$.

Moreover, "just if" is here to be interpreted via the bidirectional entailment in

$$\Vdash \mathcal{T}\ulcorner A\urcorner \Leftrightarrow \Vdash A.$$

It is a fact that

$$\Vdash \neg\mathcal{T}\ulcorner A\urcorner \Leftrightarrow \Vdash \neg A,$$

so, consequently,

$$\Vdash^{M} \mathcal{T}\ulcorner A\urcorner \Leftrightarrow \Vdash^{M} A.$$

The connectives are not truth-functional in librationism, but they are *way-functional*, and can be accounted for by following classical interdefinability connections as in any Boolean algebra: The way of the negation $\neg A$ of A, is truth minus the way of A, and the way of the conjunction $A \wedge B$ is the intersection of the way of A and the way of B. The ways of sentences built up from other connectives follow from their definitions in terms $\neg$ and $\wedge$.

According to librationism, a true paradoxical sentence L and its true companion sentence $\neg L$ complement each other. For the way of L, as defined in Definition 16.7.2, is in such a case a set of ordinals with $\varpi$ as least upper bound, whereas as well the way of $\neg L$ is a set of ordinals with $\varpi$ as least upper bound; moreover, the ways of L and $\neg L$ do not overlap. Thus, by the Definition 16.7.4, L does not express the same as what $\neg L$ expresses, for L and $\neg L$ have different ways.

## 16.8 Argumenta *ad paradoxo*

That an assumption in **£** has the consequence that $\Vdash A$ and $\Vdash \neg A$ does not suffice as a proof by contradiction against the assumption. Instead, if the considerations leading to $\Vdash A$ and $\Vdash \neg A$ cannot be extended to arrive at $\Vdash A \wedge \neg A$, they constitute an *argumentum ad paradoxo* to show that A and $\neg A$ are complementary theses of **£**. The considerations in § 16.6 are *argumenta ad paradoxo*, which justify that metamathematical statements as L and $\neg L$ are true complementary theses of **£**.

# 17 The manifestation sets

*"Der gives yderst faae Sætninger i den høiere Analyse som ere bevisede med overbevisende Strænghed. Overalt finder man den ulykkelige Maade at slutte fra det Specielle til det Almindelige, og yderst mærkværdigt er det at der efter en saadan Fremgangsmaade dog kuns findes faae af de saakaldte Paradoxer. Det er virkelig meget interessant at eftersøge Grunden hertil."*

---

Niels Henrik Abel in (Abel 1902a, p. 22).[a][b]

[a] "Il n'y a que très peu de propositions, dans l'analyse supérieure, qui soient démontrées avec une rigueur décisive. Partout on trouve la malheureuse manière de conclure du particulier au général, et il est très singulier qu'avec une pareille méthode, il ne se trouve malgré tout que peu de ce qu'on appelle paradoxes. Il est vraiment très intéressant d'en rechercher la raison." –

(Abel 1902a, p. 22) translated to (Abel 1902b, p. 23).

[b] "There are very few theorems in advanced analysis which have been demonstrated in a logically tenable manner. Everywhere one finds this miserable way of concluding from the special to the general and it is extremely peculiar that such a procedure led to so few of the so-called paradoxes." –

(Abel 1902a, p. 22) translated to (Ore 1957, p. 114).

We explain the manifestation set construction in § 17.1, and will as from § 25 see that it facilitates **£**'s ability to be extended with strong set theoretic principles. In § 17.2 we show how we may obtain Quine atoms via orthodox manifestation sets. The foci in the succeeding sections will be upon *negative* results: In § 18.0.1 we account for the *autocombative* paradox. Next, in § 18, we elucidate the virtually universal paradoxicality of power sets. Finally the failure of extensionality in **£** is discussed in § 19, where it is shown that all orthodox sets are distinct from, as well as co-extensional with infinitely many co-extensional and pairwise distinct orthodox sets.

## 17.1 The manifestation set construction

The following construction improves (Bjørdal 2012, pp. 345–346), (Cantini 1996, p. 76), (Visser 1989, pp. 695–696) and earlier literature referred to there. One may, plausibly, find that Roger's theorem and Kleene's second recursion theorem are related, but the proof that there are manifestation sets does not rely upon any presuppositions on computability.

17.1.1. Definition Kuratowskian ordered pairs

$$\langle a,b\rangle \Longleftrightarrow \{\{a\},\{a,b\}\}$$

17.1.2. Definition of the Manifestation set of $A(w,x)$

$$\llbracket^{v}_{\substack{w\\x}} A(w,\{y|(y,x)\in x\})\rrbracket \stackrel{\text{def}}{=}$$
$$\{v|(v,\{(w,x)|A(w,\{y|(y,x)\in x\})\})\in\{(w,x)|A(w,\{y|(y,x)\in x\})\}\}.$$

17.1.3. Theorem on orthodox manifestation sets

If formula $A(w,x)$ is orthodox, then

$$\Vdash^{M} \forall u(u \in \llparenthesis\substack{v\\w\\x} A(w,\{y|(y,x) \in x\})\rrparenthesis \leftrightarrow A(u,\llparenthesis\substack{v\\w\\x} A(w,\{y|(y,x) \in x\})\rrparenthesis)).$$

*Proof.* As $A(w,x)$ is orthodox also $A(w,\{y|(y,x) \in x\})$ is orthodox. Given Definition 17.1.2,

$$\Vdash^{M} \forall u\Big(u \in \llparenthesis\substack{v\\w\\x} A(w,\{y|(y,x) \in x\})\rrparenthesis$$
$$\updownarrow$$
$$u \in \{v|(v,\{(w,x)|A(w,\{y|(y,x) \in x\})\}) \in \{(w,x)|A(w,\{y|(y,x) \in x\})\}\}\Big).$$

As it is assumed that $A(w,\{y|(y,x) \in x\})$ is orthodox, also $\{(w,x)|A(w,\{y|(y,x) \in x\}$ and $\{v|(v,\{(w,x)|A(w,\{y|(y,x) \in x\})\}) \in \{(w,x)|A(w,\{y|(y,x) \in x\})\}\}$ are orthodox. So:

$$\Vdash^{M} \forall u\Big(u \in \{v|(v,\{(w,x)|A(w,\{y|(y,x) \in x\})\}) \in \{(w,x)|A(w,\{y|(y,x) \in x\})\}\}$$
$$\updownarrow$$
$$(u,\{(w,x)|A(w,\{y|(y,x) \in x\})\}) \in \{(w,x)|A(w,\{y|(y,x) \in x\})\}\}\Big).$$

By orthodoxy:

$$\Vdash^{M} \forall u\Big((u,\{(w,x)|A(w,\{y|(y,x) \in x\})\}) \in \{(w,x)|A(w,\{y|(y,x) \in x\})\}\}$$
$$\updownarrow$$
$$A(u,\{v|(v,\{(w,x)|A(w,\{y|(y,x) \in x\}) \in \{(w,x)|A(w,\{y|(y,x) \in x\}\}\Big).$$

Given Definition 17.1.2:

$$\Vdash^{M} A(u,\{v|(v,\{(w,x)|A(w,\{y|(y,x) \in x\})\}) \in \{(w,x)|A(w,\{y|(y,x) \in x\}\}\}).$$
$$\updownarrow$$
$$A(u,\llparenthesis\substack{v\\w\\x} A(w,\{y|(y,x) \in x\})\rrparenthesis).$$

The final step of the proof is by appealing to the transitivity of proven conditionals. □

17.1.4. Theorem Comprehension for orthodox manifestation set with parameters

Some manifestation sets have parameters, so for orthodox $A(\vec{v},w,x)$:

$$\Vdash^{M} \forall u\forall\vec{v}(u \in \llparenthesis\substack{v\\w\\x} A(\vec{v},w,\{y|(y,x) \in x\})\rrparenthesis \leftrightarrow A(\vec{v},u,\llparenthesis\substack{v\\w\\x} A(\vec{v},w,\{y|(y,x) \in x\})\rrparenthesis)).$$

*Proof.* Adjust Definition 17.1.2. For the notation, consult Definition 4.5.19. □

17.1.5. Exercise on General comprehension for manifestation sets

Show that for arbitrary $A(w, x)$,

$$\Vdash^{M} \forall u(u \in \llparenthesis \substack{v\\w\\x} A(w, \{y|(y, x) \in x\}) \rrparenthesis \leftrightarrow \mathcal{T}\mathcal{T} A(u, \llparenthesis \substack{v\\w\\x} A(w, \{y|(y, x) \in x\}) \rrparenthesis).$$

## 17.2 Quine atoms

The most elementary Quine atom is the manifestation set $\llparenthesis = \rrparenthesis$ of formula $=$. By means of manifest comprehension,

$$\Vdash^{M} \forall u(u \in \llparenthesis \substack{v\\w\\x} w = \{y|(y, x) \in x\} \rrparenthesis \leftrightarrow \mathcal{T}\mathcal{T} u = \llparenthesis \substack{v\\w\\x} w = \{y|(y, x) \in x\} \rrparenthesis ). \tag{17.2.1}$$

As identity is an orthodox relation,

$$\Vdash^{M} \forall u(u \in \llparenthesis \substack{v\\w\\x} w = \{y|(y, x) \in x\} \rrparenthesis \leftrightarrow u = \llparenthesis \substack{v\\w\\x} w = \{y|(y, x) \in x\} \rrparenthesis ). \tag{17.2.2}$$

As identity is an equivalence relation,

$$\Vdash^{M} \llparenthesis \substack{v\\w\\x} w = \{y|(y, x) \in x\} \rrparenthesis = \llparenthesis \substack{v\\w\\x} w = \{y|(y, x) \in x\} \rrparenthesis. \tag{17.2.3}$$

So from equations 17.2.2 and 17.2.3,

$$\Vdash^{M} \llparenthesis \substack{v\\w\\x} w = \{y|(y, x) \in x\} \rrparenthesis \in \llparenthesis \substack{v\\w\\x} w = \{y|(y, x) \in x\} \rrparenthesis. \tag{17.2.4}$$

17.2.5. Exercise Prove that there are infinitely many distinct Quine atoms.

# 18 Powersets are paradoxical lest as $\mathcal{P}(\{\mathbf{v}|\mathbf{v} = \mathbf{v}\})$

> *"Dagegen scheint mir aber jede überflüssige Einengung des mathematischen Forschungstriebes eine viel grössere Gefahr mit sich zu bringen und eine um so grössere, als dafür aus dem Wesen der Wissenschaft wirklich keinerlei Rechtfertigung gezogen werden kann; denn das Wesen der Mathematik liegt gerade in ihrer* Freiheit.*"*
>
> `Georg Cantor` in (Cantor 1883a, p. 20) and (Cantor 1883b, p. 564).

We first establish

18.0.1. Theorem Power sets are paradoxical, lest of an unpardoxical universal set.

In contrast to orthodox manifestation sets, many are paradoxical. This is for example the case with the following quite heretical manifestations set å, the *autocombatant*, which is amongst the most surprising paradoxical sets of librationism,

18.0.2. Theorem (The paradoxical autocombative truths)

$\mathring{a}$ is the manifestation set $\llbracket^{\mathrm{v}}_{\substack{\mathrm{w}\\ \mathrm{x}}}\, \mathrm{w} \notin \{\mathrm{y}|(\mathrm{y},\mathrm{x}) \in \mathrm{x}\}\rrbracket$, for which: $\Vdash \forall \mathrm{v}(\mathrm{v} \in \mathring{a}$ & $\Vdash \forall \mathrm{v}(\mathrm{v} \notin \mathring{a})$.

We prove this on account of the lemma that

$$\Vdash^{\mathrm{M}} \forall \mathrm{v}(\mathrm{v} \in \mathring{a} \leftrightarrow \mathcal{T}\mathcal{T}\mathrm{v} \notin \mathring{a}). \tag{18.0.3}$$

*Proof.* If $\lambda$ is any limit below the closure ordinal $\varpi$, we will, for any term a, and the barter at $\lambda$, $\VDash^{\lambda}$, have that $\VDash^{\lambda}$ a $\notin$ å; otherwise a contradiction would follow as a $\notin$ **å** would hold at succeeding successor ordinals $\sigma$, $\sigma + 1$ and $\sigma + 2$ below $\lambda$. Consequently, on account of Equation (18.0.3), $\VDash^{\lambda+2}$ a $\in$ **å**. As a was arbitrary, $\VDash^{\lambda} \forall \mathrm{v}(\mathrm{v} \notin \mathbf{\mathring{a}})$ and $\VDash^{\lambda+2} \forall \mathrm{v}(\mathrm{v} \in \mathbf{\mathring{a}}))$. As a result, $\VDash^{\varpi} \neg\mathcal{T}\neg\forall \mathrm{v}(\mathrm{v} \in \mathbf{\mathring{a}}))$ and $\VDash^{\varpi} \neg\mathcal{T}\neg\forall \mathrm{v}(\mathrm{v} \notin \mathbf{\mathring{a}}))$. The proof finishes by invoking Definitions 5.5.1.5 and 5.5.2.3. □

The autocombative theses are instrumental in showing, here, in §18, that the power set $\mathcal{P}(\mathrm{a})$ is paradoxical, lest a is universal and orthodox. This phenomenon concerning power sets is of importance for the librationist account of infinity, as it blocks the Cantorian conclusion that there are *uncountable* infinities.

Use standard notation, so that $\Vdash^{\mathrm{M}} \mathrm{a} \subset \mathrm{b} \leftrightarrow \forall \mathrm{x}(\mathrm{x} \in \mathrm{a} \rightarrow \mathrm{x} \in \mathrm{b})$, and posit

18.0.4. Definition The power set of a

$$\mathcal{P}(\mathrm{a}) \Longleftarrow\!\!\!= \{\mathrm{v}|\mathrm{v} \subset \mathrm{a}\}.$$

18.0.5. Definition The universal set

$$\mathbf{U} \Longleftarrow\!\!\!= \{\mathrm{v}|\mathrm{v} = \mathrm{v}\}$$

18.0.6. Theorem The quasi-universal paradoxicality of power sets

$$\mathcal{P}(\mathrm{a}) \text{ is paradoxical just if } \nVdash^{\mathrm{M}} \forall \mathrm{x}(\mathrm{x} \in \mathrm{a} \leftrightarrow \mathrm{x} \in \mathbf{U}).$$

*Proof.* We use a case distincion to provide a distinct proof for the case where $\Vdash^{\mathrm{M}} \exists \mathrm{v}(\mathrm{v} \notin \mathrm{a})$.

(1) If $\Vdash^{\mathrm{M}} \exists \mathrm{v}(\mathrm{v} \notin \mathrm{a})$, use the autocombatant $\mathbf{\mathring{a}}$, of Theorem 18.0.2, for which

$$\Vdash \forall \mathrm{v}(\mathrm{v} \in \mathbf{\mathring{a}}) \ \& \ \Vdash \forall \mathrm{v}(\mathrm{v} \notin \mathbf{\mathring{a}}).$$

In this case $\Vdash \mathbf{\mathring{a}} \notin \mathcal{P}(\mathrm{a})$ and $\Vdash \mathbf{\mathring{a}} \in \mathcal{P}(\mathrm{a})$, so $\mathcal{P}(\mathrm{a})$ is paradoxical.

(2) If $\Vdash \exists \mathrm{v}(\mathrm{v} \notin \mathrm{a})$ & $\Vdash \forall \mathrm{v}(\mathrm{v} \in \mathrm{a})$, $\Vdash \mathbf{U} \in \mathcal{P}(\mathrm{a})$ & $\Vdash \mathbf{U} \notin \mathcal{P}(\mathrm{a})$ and $\mathcal{P}(\mathrm{a})$ is paradoxical. □

# 19 Non-extensionality and nominal ursets

*"one might be tempted to call the axiom of extensionality "analytic," true by virtue of the meanings of the words contained in it, but not to consider the other axioms analytic."*

---

*George Boolos*, in (G. Boolos 1971, p. 501).

The principle of extensionality's failure in type free theories is well known, and many have contributed to the deposit of knowledge.

Let us first posit

19.0.1. Definition

The principle of extensionality holds for a theory just if it proves

$$a \stackrel{E}{=} b \rightarrow a = b.$$

A particularly easy semantical proof of the failure of the extensionality principle in **£** is obtained by observing that for any limit ordinal $\lambda$,

$$\VDash^{\lambda+1} \{v|v = v\} \stackrel{E}{=} \{v|v \notin v\} \wedge \{v|v = v\} \neq \{v|v \notin v\}.$$

As a consequence, there are sets a and b such that $\nVDash a \stackrel{E}{=} b \rightarrow a = b$. But it is a soundness requirement that $\Vdash^{M} a \stackrel{E}{=} b \rightarrow a = b \Rightarrow \VDash a \stackrel{E}{=} b \rightarrow a = b$ , so that $\nVdash^{M} a \stackrel{E}{=} b \rightarrow a = b$.

(Gilmore 1974), which appeared in (Jech 1974), showed that a partial set theory proves that there is an orthodox set a such that $a \stackrel{E}{=} \emptyset$ and $a \neq \emptyset$. (Bjørdal 2012, p. 345) relates Lev Gordeev's more concise proof of the same result as Gilmore's, in the context of Explicit Mathematics, and some on why it was published in (Beeson 1985), with acknowledgement.

Define Gordeev's set with the manifestation theorem so that one may posit

19.0.2. Definition (Via manifestation point)$\forall x(x \in \dot{g} \leftrightarrow \mathcal{T}\mathcal{T}(x = \emptyset \wedge x = \dot{g}))$.

19.0.3. Theorem[Gordeev] $\dot{g}$ is (i) orthodox, so $\Vdash^{M} x \in \dot{g} \leftrightarrow (x = \emptyset \wedge x = \dot{g})$, (ii) empty and (iii) distinct from $\emptyset$.

*Proof.* As the proof of Theorem 4 in (Bjørdal 2012, p. 345): (i) $\dot{g}$ is orthodox, on account of the theory of identity. (ii) As $\Vdash^{M} x \in \dot{g} \rightarrow (x = \emptyset \wedge x = \dot{g})$, $\Vdash^{M} x \in \dot{g} \rightarrow \dot{g} = \emptyset$, so $\dot{g}$ is empty. (iii) $\dot{g} \neq \emptyset$, for else $\dot{g} = \{\dot{g}\}$ on account of Theorem 19.0.3 (i), which contradicts (ii). □

19.0.4. Definition The Gordeev-progression

(1) $\dot{g}_0 \Longleftarrow\!\!= \emptyset$.

(2) $\dot{g}_{l+1} \Longleftarrow\!\!= \llbracket {}^{v}_{w}{}_{x}\, x = y \wedge \bigvee_{i=0}^{i=l} x = \dot{g}_i \rrbracket$ for any natural number l.

19.0.5. THEOREM (THE *URSETS*)

For all natural numbers $m$ and $n$,

(1) $\dot{g}_m$ is orthodox.

(2) $\forall x(x \in \dot{g}_{m+1} \leftrightarrow (x = \dot{g}_{m+1} \wedge \bigvee_{i=0}^{i=m} x = \dot{g}_i))$.

(3) $\dot{g}_m$ is empty.

(4) $\dot{g}_m$ and $\dot{g}_n$ are distinct if $m$ and $n$ are distinct.

*Proof.* (1) holds on account of the theory of identity. (2) is a consequence of (1) and Definition 19.0.2. (3) If $m = 0$ we are done, as $\dot{g}_0$ by Definition 19.0.4 (1) is the canonical empty set $\emptyset$. If $m = n + 1$, then, given (2), $b \in \dot{g}_{n+1}$ only if $b = \dot{g}_{n+1}$ and $(b = \dot{g}_0 \dots \vee \dots b = \dot{g}_n)$, where all of $\dot{g}_0, \dots, \dot{g}_n$ are orthodox and empty; so b is a member of an empty set if $b \in \dot{g}_{n+1}$, thus $\dot{g}_{n+1} = \dot{g}_m$ is empty. (4) It is sufficient to show that $\dot{g}_m$ and $\dot{g}_{n+1}$ are distinct if $m < n + 1$. As by (2), $\Vdash^M \dot{g}_{n+1} \in \dot{g}_{n+1} \leftrightarrow \dot{g}_{n+1} = \dot{g}_{n+1} \wedge (\dot{g}_{n+1} = \dot{g}_0 \dots \vee \dots \dot{g}_{n+1} = \dot{g}_n)$, and $\Vdash^M \dot{g}_{n+1} \notin \dot{g}_{n+1}$ on account of (3), it follows that $\Vdash^M (\dot{g}_{n+1} \neq \dot{g}_0 \dots \wedge \dots \dot{g}_{n+1} \neq \dot{g}_n)$ - and we are done. □

## 20 Transfinite urelements

*"Os números são as regras dos seres, e a matemática, é o regulamento do mundo."*

Francisco G. Teixeira in *Conferência sobre o poder e a beleza das Matemáticas,* (Gomes Teixeira 1925, p. 271).

In §19 the nominal ursets $\ulcorner n \urcorner$ for all $n \varepsilon \Omega$ were isolated, and the construction made use of the finite meta-ordinals $0, 1, 3, \ldots$ of $\Omega$.

The philosophically minded librationist set theorist may want urelements even beyond the nominal ursets of §19, so that $\ulcorner \alpha \urcorner$ for any transfinite meta-ordinal $\alpha$ with a standard notation is the $\alpha$'th empty set such that for all $\beta \prec \alpha$: $\ulcorner \alpha \urcorner \overset{E}{=} \ulcorner \beta \urcorner \ \& \ \ulcorner \alpha \urcorner \neq \ulcorner \beta \urcorner$.

We add the

20.0.1. Postulation of the transfinite urelement constants:

> The language of **£** has constant $c^{\alpha}$ just if transfinite meta-ordinal $\alpha$ has a standard notation, and any member $\beta$ of the ordinal $\alpha$ has an associated urelement $\ulcorner \beta \urcorner$.

20.0.2. Definition Given Definition 19.0.2.

$$\Omega^{\ulcorner\urcorner} \Longleftrightarrow \{\!\!\{x \vdots \mathbf{\Pi} y(\ulcorner 0 \urcorner \in y \land \mathbf{\Pi} n(n\varepsilon\Omega \rightarrow (\ulcorner n \urcorner \in y \rightarrow \ulcorner n+1 \urcorner \in y) \rightarrow x \in y)\}\!\!\}.$$

20.0.3. Theorem on manifestation set $\ulcorner \Omega \urcorner$

$$\Vdash^{M} \forall x(x \in \ulcorner \Omega \urcorner \leftrightarrow (x = \ulcorner \Omega \urcorner \land \exists y(y \in \Omega^{\ulcorner\urcorner} \land x = y))).$$

20.0.4. Remark

$\upsilon_{\omega}$ in Theorem 20.0.3 is orthodox as $\upsilon^{\omega}$ in Definition 20.0.2 is orthodox.

20.0.5. Theorem

(i) $\upsilon_{\omega}$ is empty. (ii) $\Vdash^{M} \forall y(y \in \upsilon^{\omega} \rightarrow \upsilon_{\omega} \neq y)\}$.

*Proof.* (i) If a were an element of $\upsilon_{\omega}$, $\Vdash^{M} a = \upsilon_{\omega} \land \exists y(y \in \upsilon^{\omega} \land a = y)$. Given Theorem 21.1.11, all members of $\upsilon^{\omega}$ are empty sets. Consequently, if a were an element of $\upsilon_{\omega}$ then $\upsilon_{\omega}$ would be an empty set. So $\upsilon_{\omega}$ is an empty set. (ii) A rendition of Theorem 20.0.3 is $\Vdash^{M} \forall x(x \notin \upsilon_{\omega} \leftrightarrow (x = \upsilon_{\omega} \rightarrow \forall y(y \in \upsilon^{\omega} \rightarrow x \neq y)))$, so, as a consequence,

$$\Vdash^{M} \upsilon_{\omega} \notin \upsilon_{\omega} \rightarrow (\upsilon_{\omega} = \upsilon_{\omega} \rightarrow \forall y(y \in \upsilon^{\omega} \rightarrow \upsilon_{\omega} \neq y)).$$

The proof finishes by invoking the maxim mode 7.2.7, as $\Vdash^{M} \upsilon_{\omega} \notin \upsilon_{\omega}$ on account of (i). □

20.0.6. Definition

(1) c is an urset just if $c = \ulcorner n \urcorner$, for some $n \in \Omega_{+}$.

(2) In accordance with Definition 4.5.3, $\ulcorner 2222221 \urcorner$, $\ulcorner 222221 \urcorner$, $\ulcorner 22221 \urcorner$, $\ulcorner 2221 \urcorner$, $\ulcorner 221 \urcorner$, $\ulcorner 21 \urcorner$ and $\ulcorner 1 \urcorner$ are the *symbolic* ursets: $\ulcorner * \urcorner$, $\ulcorner \ddot{c} \urcorner$, $\ulcorner ç \urcorner$, $\ulcorner \forall \urcorner$, $\ulcorner \downarrow \urcorner$, $\ulcorner \ddot{v} \urcorner$, $\ulcorner \bullet \urcorner$.

(3) The symbolic ursets are the *atomic* names, which denote the primitive symbols.

(4) Recall Definition 4.5.4.3 of $\ell(n_0) \mathrel{=\!=} \lfloor \log_2(n_0 + 1) \rfloor$, which uses $\log_2$ and the floor function $\lfloor\ \rfloor$, to define the length $\ell(n_0)$ of the bijective base-8 metanumeral needed for given metanumber $n_0$.

(5) So

$$
\begin{aligned}
\ell(\bullet) &= \ell(1) = 1\\
\ell(\ddot{v}) &= \ell(21) = 2\\
\ell(\downarrow) &= \ell(221) = 3\\
\ell(\forall) &= \ell(2221) = 4\\
\ell(\varsigma) &= \ell(22221) = 5\\
\ell(\ddot{c}) &= \ell(222221) = 6\\
\ell(*) &= \ell(2222221) = 7.
\end{aligned}
$$

(6) Given Definition 4.5.4, the *joining of names* is defined by positing

$$
\ulcorner n_0 \urcorner\ \ulcorner n_1 \urcorner \mathrel{=\!=} \ulcorner n_0 n_1 \urcorner \mathrel{=\!=} \ulcorner n_0 \frown n_1 \urcorner \mathrel{=\!=} \ulcorner n_0 \cdot 2^{\ell(n_1)} + n_1 \urcorner .
$$

(7) Given Definition 20.0.6.5, we may use Definition 20.0.6.6 to construe composite names gramatically correct by joining *names* whilst obeying the formation rules of § 4.5.

(8) For good $\in$ {symbol, symbol string, variable, formula, constant, term, sentence}, $\ulcorner \mathcal{N} \urcorner$ is a good name just if $\mathcal{N}$ is a good.

20.0.7. Caveat In formula $\forall v \mathcal{T} \ulcorner A \urcorner$, $\ulcorner\ \urcorner$ is a term operating formula forming operator, so the evaluation of $\mathcal{T} \ulcorner A \urcorner {}^{b}_{v}$ is comparable with $\Box A^{b}_{v}$, where $\Box$ is any formula operating formula forming operator. So, for example, $\mathcal{T} \ulcorner v = v \urcorner {}^{b}_{v}$ is $\mathcal{T} \ulcorner b = b \urcorner$. A subtle substitution function, e.g. as with (Smoryński 1977, 837 et passim) in the proof of Gödel's incompleteness theorem, is not needed, for there is no use of quantification into an opaque, or otherwise "intensional", context.

Notice that at this point Theorem 21.1.11 may straightforwardly be extended:

20.0.8. Theorem (Sets of urelements to any order)

$\upsilon^{\omega}$ in Definition 20.0.2 serves as the set of the expression names defined in § 19. Given Theorems 20.0.3 and 20.0.5 (i), $\upsilon_{\omega}$ is another empty set distinct from all members of $\upsilon^{\omega}$ But we may now define a new omega ordered set of ursets

$$
\upsilon^{\omega \cdot 2} \mathrel{=\!=} \{x | \forall y(\upsilon_{\omega} \in y \wedge \forall z(\upsilon_z \in y \rightarrow \upsilon_{z+1} \in y) \rightarrow x \in y)\}.
$$

$\upsilon^{\omega \cdot 2}$, and indeed $\upsilon^{\beta}$ for any ordinal $\beta$, may serve as sets of ursets, or *urelements*, for whatever purpose one may have in mind, including that of naming extramathematical things to equip **£** with domains useful for applied mathematics, including logic.

# 21 The infinitary non-extensionality of all orthodox sets

*"There are more things in Heaven and Earth, Horatio, than are dreamt of in your philosophy."*

Shakespeare, in *Hamlet*.

## 21.1 Extending Minari

(Cantini 1996)(74), relates a proof, suggested by P. Minari that we for *any* orthodox set a may find a *distinct* orthodox set b such that a and b are nevertheless co-extensional:

21.1.1. Definition The Minari set

$$\Vdash^{M} x \in b \leftrightarrow [(a = b \wedge a \notin a) \vee (a \neq b \wedge x \in a)]$$

21.1.2. Theorem

$$a \neq b.$$

*Proof.* If $b \in b$, $(a = b \wedge a \notin a) \vee (a \neq b \wedge b \in a)$; if the first disjunct, $a \neq b$ because $b \in b$ and $a \notin a$ and if the second disjunct $a \neq b$ as stated. If $b \notin b$, $(a \neq b \vee a \in a) \wedge (a = b \vee b \notin a)$ so that $a = b \rightarrow a \in a$ and hence $a = b \rightarrow b \in b$; consequently $a \neq b$ by modus ponens. □

21.1.3. Theorem

$$a \overset{E}{=} b.$$

*Proof.* Combine Definition 21.1.1 and Theorem 21.1.2 to see that $\Vdash^{M} c \in b \leftrightarrow c \in a$. □

21.1.4. Lemma (The infinitary generalization of Minari's non-extensionality)

For orthodox set $a^1$ and any natural number $n \geq 1$ there is an orthodox set $a^{n+1}$ such that

$$\Vdash^{M} x \in a^{n+1} \leftrightarrow [\bigvee_{i=1}^{i=n} (a^i = a^{n+1} \wedge a^i \notin a^i) \vee \bigwedge_{i=1}^{i=n} (a^i \neq a^{n+1} \wedge x \in a^i)].$$

*Proof.* We must show first that this definition is adequate given initial orthodox set $a^1$ and the definition by manifestation presupposed in Theorem 17.1.3 on account of the induced hereditary orthodoxy: In the base case $n = 1$ the biconditional in the lemma coincides with Minari's original example, and the definition is adequate on account of the orthodoxy induced by the logic of identity. If all $a^i$, for $1 \leq i \leq n$, are orthodox, $a^{n+1}$ is again orthodox on account of the logic of identity. □

21.1.5. Lemma

With notation as in Lemma 21.1.4,

$$\Vdash^{M} \bigwedge_{i=1}^{i=n} a^i \neq a^{n+1}.$$

*Proof.* Suppose $a^{n+1} \in a^{n+1}$. If the first disjunct on the right hand side of Lemma 21.1.4 were to hold, a positive $i \leq n$ would be such that $a^i = a^{n+1}$ and as well $a^i \notin a^i$, which is impossible; so the second disjunct with the target result would hold. Suppose on the other hand that $a^{n+1} \notin a^{n+1}$. It follows by means of of Lemma 21.1.4 that

$$\bigwedge_{i=1}^{i=n} (a^i \neq a^{n+1} \vee a^i \in a^i) \wedge \bigvee_{i=1}^{i=n} (a^{n+1} \in a^i \rightarrow a^i = a^{n+1}),$$

and the first conjunct is sufficient to finish the proof, as surely $a^i \in a^i$ entails that $a^i \neq a^{n+1}$ in the presence of $a^{n+1} \notin a^{n+1}$. □

21.1.6. THEOREM Given an orthodox initial set for the construction of $a^{n+1}$, the latter is orthodox and co-extensional with as well as distinct from any $a^i$, for positive $i < n + 1$:

$$\bigwedge_{i=1}^{i=n} (a^i \neq a^{n+1} \wedge a^i \stackrel{E}{=} a^{n+1}).$$

*Proof.* Combine Lemma 21.1.4 and Lemma 21.1.5 to get:

$$\Vdash^M \forall x (x \in a^{n+1} \leftrightarrow \bigwedge_{i=1}^{i=n} (x \in a^i \wedge a^i \neq a^{n+1})),$$

□

and keep track of the successive inequalities in the progression.

21.1.7. DEFINITION

Let $\llparenthesis a, n \rrparenthesis$, for orthodox set a and natural number $n > 0$, be the manifestation set of

$$\bigwedge_{i=1}^{i=n} (x \in \upsilon_i \wedge \upsilon_i \neq v_{n+1}) \vee (\bigwedge_{i=1}^{i=n} \upsilon_i \notin \upsilon_i \wedge \bigvee_{i=1}^{i=n} \upsilon_i = v_{n+1}).$$

21.1.8. LEMMA

$$\Vdash^M \forall x \Big( x \in \llparenthesis a, n \rrparenthesis \leftrightarrow (\bigwedge_{i=1}^{i=n} (x \in \upsilon_i \wedge \upsilon_i \neq \llparenthesis a, n \rrparenthesis) \vee (\bigwedge_{i=1}^{i=n} \upsilon_i \notin \upsilon_i \wedge \bigvee_{i=1}^{i=n} \upsilon_i = \llparenthesis a, n \rrparenthesis))\Big).$$

*Proof.* By the canonical argument for the existence of manifestation sets. □

21.1.9. LEMMA

$$\bigwedge_{i=1}^{i=n} \upsilon_i \neq \llparenthesis a, n \rrparenthesis.$$

*Proof.* If $\bigvee_{i=1}^{i=n}(\upsilon_i = \llparenthesis a,n \rrparenthesis)$, Lemma 21.1.8 reduces to

$$(x \in \llparenthesis a,n \rrparenthesis \leftrightarrow (\bigwedge_{i=1}^{i=n} \upsilon_i \notin \upsilon_i \wedge \bigvee_{i=1}^{i=n} \upsilon_i = \llparenthesis a,n \rrparenthesis)). \tag{21.1.10}$$

If we instantiate the left side of Equation 21.1.10 with $\llparenthesis a,n \rrparenthesis \in \llparenthesis a,n \rrparenthesis$, it follows that for some positive $i < n$: $\upsilon_i \notin \upsilon_i \wedge \upsilon_i = \llparenthesis a,n \rrparenthesis$. So $\llparenthesis a,n \rrparenthesis \in \llparenthesis a,n \rrparenthesis$ only if $\llparenthesis a,n \rrparenthesis \notin \llparenthesis a,n \rrparenthesis$. If we instantiate the left side with $\llparenthesis a,n \rrparenthesis \notin \llparenthesis a,n \rrparenthesis$, it follows that $(\bigvee_{i=1}^{i=n} \upsilon_i \in \upsilon_i \vee \bigwedge_{i=1}^{i=n} \upsilon_i \neq \llparenthesis a,n \rrparenthesis)$, and so consequently $\llparenthesis a,n \rrparenthesis \in \llparenthesis a,n \rrparenthesis$ only if $\llparenthesis a,n \rrparenthesis \notin \llparenthesis a,n \rrparenthesis$. So assumption $\bigvee_{i=1}^{i=n}(\upsilon_i = \llparenthesis a,n \rrparenthesis)$ has the consequence that $\llparenthesis a,n \rrparenthesis \in \llparenthesis a,n \rrparenthesis$ just if $\llparenthesis a,n \rrparenthesis \notin \llparenthesis a,n \rrparenthesis$, and we are done. □

21.1.11. Theorem

Given orthodox set a, the members of the infinite series $\llparenthesis a,1 \rrparenthesis$, $\llparenthesis a,2 \rrparenthesis, \llparenthesis a,3 \rrparenthesis, \ldots$ are pairwise distinct and co-extensional with a.

Theorem 5 (ii) in (Bjørdal 2012)(346), whose proof was left as an exercise, stated the result that Minari's construction can be generalized, as in Theorem 21.1.11. This is the most general non-extensionality result possible, as it were, and the mentioned exercise is solved, by the following verification of Theorem 21.1.11.

Proof:

LSRS - Suppose $\llparenthesis a,n \rrparenthesis \in \llparenthesis a,n \rrparenthesis$. It follows that e.g. $\upsilon_1 = \llparenthesis a,n \rrparenthesis$. Then, from the conjunct $\bigwedge_{i=1}^{i=n} \upsilon_i \notin \upsilon_i$, $\llparenthesis a,n \rrparenthesis \notin \llparenthesis a,n \rrparenthesis$.

RSLS - instantiating, and negating both sides:

$$(\llparenthesis a,n \rrparenthesis \notin \llparenthesis a,n \rrparenthesis \leftrightarrow (\bigvee_{i=1}^{i=n} \upsilon_i \in \upsilon_i \vee \bigwedge_{i=1}^{i=n} \upsilon_i \neq \llparenthesis a,n \rrparenthesis)); \tag{21.1.12}$$

but $a \in \llparenthesis a,n \rrparenthesis \leftrightarrow a \notin \llparenthesis a,n \rrparenthesis$ is impossible, so we have established that

$$\bigwedge_{i=1}^{i=n}(\upsilon_i \neq \llparenthesis a,n \rrparenthesis), \tag{21.1.13}$$

and so

$$\forall x(x \in \llparenthesis a,n \rrparenthesis \leftrightarrow (\bigwedge_{i=1}^{i=n} x \in \upsilon_i \wedge \upsilon_i \neq \llparenthesis a,n \rrparenthesis \wedge \llparenthesis a,n \rrparenthesis \stackrel{E}{=} a). \tag{21.1.14}$$

Consequently,

$$\forall x(x \in \llparenthesis a,n \rrparenthesis \leftrightarrow \bigwedge_{i=1}^{i=n} x \in \upsilon_i \wedge \llparenthesis a,n \rrparenthesis \stackrel{E}{=} a). \tag{21.1.15}$$

It follows by induction that

$$\bigwedge_{i=1}^{i=n} (\upsilon_i \stackrel{E}{=} \llbracket a, n \rrbracket). \tag{21.1.16}$$

and so it follows that

$$\bigwedge_{i=1}^{i=n} (\llbracket a, n \rrbracket \stackrel{E}{=} a). \tag{21.1.17}$$

Given Equations 21.1.13 and 21.1.17, the infinite omega ordered series

$$\sum_{n=1}^{\infty} \llbracket a, n \rrbracket$$

has just orthodox, pairwise distinct sets which are co-extensional with a.

# 22 The Holder, its heritors and the regulars

*“Совершенно согласен ... что смирение есть величайшая и необходимая добродетель. Как я всегда говорю, человек подобен дроби, в которой знаменатель определяет его мнение о самом себе. Самое лучшее, когда знаменатель этот ноль (полное смирение), а ужасно, когда знаменатель этот возрастает до бесконечности. В первом случае, каков бы ни был знаменатель, он имеет действительное значение, во втором же случае — никакого.”*

Leo Tolstoy, in (Tolstoy 1910).[a]

[a]“I perfectly agree with you that humility is the greatest and most needful virtue. As I always say, man is like a fraction, in which the denominator indicates his opinion about himself. It is best for this denominator to be zero (complete humility), and it is terrible when it is augmented to infinity. In the first case, man has a true significance, whatever the denominator; but in the second case — none.” – (Tolstoy 1910) translated to (Eves 1988, p. 81).

Heritors and regulars are defined, and their behavior regulated so as to support the development of the interpretation of NBG set theory of § 25.

22.0.1. Definition

$$\text{The } \textit{Holder} \text{ is } \mathcal{H} \Longleftarrow \{x|x = \{y|y \in x\}\}.$$

22.0.2. Definition

$$a \text{ is an } \textit{heritor} \text{ just if } \mathcal{H}(a).$$

22.0.3. Definition

$$\mathcal{H}(a) \Longleftarrow a = \{x|x \in a\}$$

22.0.4. Theorem

The Holder and its heritors are orthodox.

*Proof.* The Holder is orthodox by identity theory, and heritors by Postulate 6.3.2.6. □

22.0.5. Definition

$$a \text{ is an } \textit{hyposet} \text{ of set } b \text{ just if } a = \{x|x \in a \land x \in b\}.$$

22.0.6. Axiom

$$\Vdash^{M} \mathcal{H}(a) \land \mathcal{H}(b) \land a \subset b \rightarrow a = \{x|x \in a \land x \in b\}.$$

22.0.7. Theorem

$$\Vdash^{M} \mathcal{H}(a) \land \mathcal{H}(b) \land a \overset{E}{=} b \rightarrow a = b.$$

*Proof.* An instance of Axiom 22.0.6 is $\Vdash^{M} \mathcal{H}(b) \land \mathcal{H}(a) \land b \subset a \rightarrow b = \{x|x \in b \land x \in a\}$. $\{x|x \in a \land x \in b\} = \{x|x \in b \land x \in a\}$, given § 10, so just wed with the statement instance $\mathcal{H}(a) \land \mathcal{H}(b) \land a \subset b \rightarrow a = \{x|x \in a \land x \in b\}$ of Axiom 22.0.6. □

22.0.8. Axiom

$$\Vdash^{M} \mathcal{H}(a) \wedge \mathcal{H}(b) \wedge a \subset b \leftarrow a = \{x | x \in a \wedge x \in b\},$$

so if a is a hyposet of b, then a and b are heritors, and a is a *subheritor* of b.

22.0.9. Theorem

$$\Vdash^{M} \mathcal{H}(a) \wedge \mathcal{H}(b) \rightarrow (a \subset b \leftrightarrow a = \{x | x \in a \wedge x \in b\}).$$

*Proof.* Invoke Axioms 22.0.6 and 22.0.8. □

22.0.10. Axiom (Heritors are hereditarily heritors)

$$\Vdash^{M} \mathcal{H}(y) \rightarrow \forall x(x \in y \rightarrow \mathcal{H}(x)).$$

22.0.11. Observation

This section's axioms do not commit to the existence of heritors.

22.0.12. Definition (Regular sets)

$$\mathcal{R}(x) \Longleftarrow\!\!= \exists y(y \in x) \rightarrow \exists y(y \in x \wedge \forall z(z \notin x \vee z \notin y))$$

22.0.13. Exercise

Show that regular hereditarily orthodox sets are hereditarily regular.

…

Our attention below will be upon regular heritors.

# 23 ZF$^-$ without extensionality plus the subterset axioms

> *"Logic is the hygiene which mathematics practices to keep its ideas healthy and strong."*
>
> Attributed to Hermann Weyl in (Kline 1958, p. 24).

The heritor axiom: $\forall a(a \in \mathbf{V} \to a = \{x|x \in a\})$.

The commutation axiom:

$$\forall a \forall b(\{x|x \in a \wedge x \in b\} = \{x|x \in b \wedge x \in a\}).$$

The subterset definition:

$$\forall a \forall b(a \text{ is a subterset of } b \text{ just if } a = \{x|x \in a \wedge x \in b\}).$$

The subterset axiom: $\forall b \exists a(a = \{x|x \in a \wedge y \in b\}\})$.

It is a theorem from comprehension that

$$\forall a \forall b(a = \{x|x \in a \wedge y \in b\} \to \forall x(x \in a \to x \in b).$$

The other direction holds for pairs of elements of V, given the heritor axiom:

The axiom that subsets restrivted to V are subtersets:

$$\forall a \forall b(\langle a, b\rangle \in V^2 \to ((a \subset b) \to a = \{x|x \in a \wedge y \in b\})).$$

*The extensionality theorem*: If a and b are coextensional, it follows from the axiom that subsets are subtersets that

$$a = \{x|x \in a \wedge x \in b\} \,\&b = \{x|x \in b \wedge x \in a\}.$$

On account of the commutation axiom, $a = \{x|x \in a \wedge y \in b\} = b$, so that $a = b$ and the extensionality principle holds.

The potencyset definition: $♂(a) = \{x|x = \{y|y \in x \wedge y \in a\}\}$.

*The potencyset axiom*: $\forall a \exists x(x = ♂(a))$.

So, why does Cantor's theorem now fail for the potency set of $\mathbb{N}$?

Assume there is a function f from $\mathbb{N}$ onto $♂(\mathbb{N})$. We are asked to consider

$$S = \{x|x \in \mathbb{N} \wedge x \notin f(x)\}.$$

For any $n \in \mathbb{N}$ such that $f(n) = S$ we have $n \in S \leftrightarrow n \notin S$.

Thus *if* S is an heritor a problem arises. Have we proved that S is not an heritor, then? Indeed, that is precisely what the Cantorian proof just did - i.e. it proved that S is not an heritor, and it did not prove that $♂(\mathbb{N})$ has a larger cardinality than $\mathbb{N}$

# 24 Choice, power, potency and countability

*"The axiom of choice is obviously true, the well-ordering principle obviously false, and who can tell about Zorn's lemma?"*

According to (Simon 2015, p. 12), Jerry L. Bona related that he taught this in class.

We show that the librationist universe is countable. Theorem 18.0.6 is one of the important reasons why that is so. Theorem 24.4.5 establishes that there is an orthodox bijection from the set of natural numbers $\omega$ to the full universe **U**. § 24.6 spells out how it is that Cantor's arguments, linked to power sets, are circumvented in **£**, with recourse to the bijection € from $\omega$ to the universe, and the choice-function $\wp`x$ upon which it is based.

## 24.1 The denumerable wellordering

24.1.1. Definition

$$\mathbf{\Pi} a, b\Big[(\mathrm{Constant}(a) \wedge \mathrm{Constant}(b) \wedge \Vdash^{\alpha} a \trianglelefteq b)$$
$$\Updownarrow$$
$$(\mu x(x\eta\Omega \;\&\; x \preceq a \;\&\; \Vdash^{\alpha} x = a) \preceq \mu y(y\eta\Omega \;\&\; y \preceq b \;\&\; \Vdash^{\alpha} y = b))\Big]$$

24.1.2. Corollary

$$\Vdash^{\alpha} a = b \Leftrightarrow \Vdash^{\alpha} a \trianglelefteq b \;\&\; \Vdash^{\alpha} a \trianglerighteq b$$

24.1.3. Definition

$$\Vdash^{\alpha} a \lhd b \Leftrightarrow \Vdash^{\alpha} a \trianglelefteq b \;\&\; \Vdash^{\alpha} a \neq b$$

24.1.4. Axiom for the wellordering

$$\Vdash^{\alpha} \forall x, y(x \lhd y \vee x = y \vee x \rhd y)$$

24.1.5. Axiom for the orthodoxy of the wellordering

$\lhd$, and its cognate relations, are orthodox.

## 24.2 Function application notation

24.2.1. Definition

$$f{`}a \Bumpeq b := \forall x \forall y \forall z(((x, y) \in f \wedge (x, z) \in f) \rightarrow y = z) \wedge (a, b) \in f.$$

24.2.2. Definition

$$b \Bumpeq f{`}a := f{`}a \Bumpeq b$$

24.2.3. Definition

$$x \in f{`}a := \exists y(f{`}a \Bumpeq y \wedge x \in y)$$

24.2.4. Definition

$$f\text{‘}a \in x := \exists y(f\text{‘}a \Bumpeq y \land y \in x)$$

24.2.5. Remark

The notation $\Bumpeq$ is used instead of =, for there are paradoxical functions as e.g.

$$g = \{(x,y)|x = \{\emptyset\} \land ((r \in r \to y = \emptyset) \land (r \notin r \to y = \{\emptyset\})\},$$

for $r = \{x|x \notin x\}$. For g we do have that $\Vdash^{M} \forall x\forall y\forall z(((x,y) \in g \land (x,z) \in g) \to y = z)$. But, notice that $\VDash^{\alpha} (\{\emptyset\},\{\emptyset\}) \in g$ just if $\VDash^{\alpha+1} (\{\emptyset\},\emptyset) \in g$. So we cannnot write $g\text{‘}\{\emptyset\} = \{\emptyset\}$ in the former case, and $g\text{‘}\{\emptyset\} = \emptyset$ in the latter case. For identity is an orthodox equivalence relations. So we use $\Bumpeq$ to avoid problems with the theory of identity in exotic cases.

24.2.6. Remark

There certainly are sets of more orthodox functions so that a function h is an element in one of them only if $\Vdash^{M} \forall x\forall y(h\text{‘}x \Bumpeq y \to \mathcal{T}\ulcorner h\text{‘}x \Bumpeq y\urcorner)$.

24.2.7. Remark

The author introduced and discussed the notation $\Bumpeq$ in the article (Bjørdal 2008, pp. 55–66), whose title's English translation is «"2+2=4" is misleading», for such reasons which are adduced here.

## 24.3 The choice function

On account of ancient Greek διάλεξε, for *was selected*, we define $ϙ\text{‘}w$, the *atled* of w:

24.3.1. Definition The choice function

$$ϙ\text{‘}w \Bumpeq \{x|(x \in w \land \forall y(y \in w \to x \trianglelefteq y))\}.$$

24.3.2. Definition Iterated choices from b

$$a \Bumpeq ϙ^{m}b \Leftrightarrow ((m = 0 \land a \Bumpeq ϙ\text{‘}b) \lor \exists n(n \in \omega \land m = n + 1 \land a \Bumpeq ϙ\text{‘}(b \setminus \bigcup_{i=0}^{i=n} ϙ^{i}b))).$$

## 24.4 The enumerator

Given Axiom 24.1.5, the orthodoxy of $\omega$ and $\mathbf{U}$, and the fact that $ϙ^{n}w$ is orthodox if w is orthodox, we posit

24.4.1. Definition of the Enumerator

$$€ = \!= \{(n,x)|n \in \omega \land x \in ϙ^{n}\mathbf{U}\}$$

24.4.2. Theorem (€ is orthdox)

*Proof.* As $\wp^{n}\mathbf{U}$, for $n \in \omega$ is orthodox. □

24.4.3. Theorem (The functionality of €)

$$\Vdash^{M} \forall x \forall y \forall z(((x, y) \in \text{€} \wedge (x, z) \in \text{€}) \rightarrow y = z)$$

*Proof.* Obvious, given Definitions 24.3.2 and 24.4.1 and Theorem 24.4.2. □

24.4.4. Theorem

$$\Vdash^{M} \text{€}\text{`}n \Bumpeq x \leftrightarrow (n, x) \in \text{€} \leftrightarrow n \in \omega \wedge x \in \wp^{n}(\mathbf{U})$$

*Proof.* On account of Definitions 24.2.1 and 24.4.1, and Theorem 24.4.3. □

24.4.5. Theorem

€ is a bijection from $\omega$ to **U**.

*Proof.* Given § 24.1, as the orders of $\omega$ and $\Omega$ match, and for any constant a, $a\eta\Omega$, as all sets are finite positive von Neumann ordinals according to the meta language. □

## 24.5 The enumeration postulates

For any ordinal $\alpha$:

24.5.1. Postulate

$$\mathbf{\Pi} a \mathbf{\Pi} b\Big(\textsc{constant}(a) \;\&\; \textsc{constant}(b) \Rightarrow \\ \Vdash^{\alpha} \forall n(n \in \omega \rightarrow \big(\exists^{=n}x(x \lhd a) \wedge \exists^{=n}y(y \lhd b) \rightarrow a = b\big))\Big).$$

24.5.2. Postulate

$$\mathbf{\Pi} a, b, c, (\Vdash^{\alpha} \forall n\Big(n \in \omega \rightarrow (\big(\langle a, b\rangle \in \text{€} \;\&\; \exists^{=n}x(x \lhd b)\big) \leftrightarrow \\ \big(\langle\{v | v \in a \vee v = a\}, c\rangle \in \text{€} \;\&\; \exists^{=(n+1)}x(x \lhd c)\big))\Big).$$

24.5.3. Postulate

$$\Vdash^{\alpha} \forall n(n \in \omega \rightarrow \exists y(\langle n, y\rangle \in \text{€})).$$

24.5.4. Postulate

$$\Vdash^{\alpha} \forall y \exists n(n \in \omega \wedge \langle n, y\rangle \in \text{€}).$$

24.5.5. Postulate

$$\Vdash^{\alpha} \forall n \forall n' \forall y(\langle n, y\rangle \in \text{€} \wedge \langle n', y\rangle \in \text{€} \rightarrow n = n').$$

24.5.6. Postulate

$$\Vdash^{\alpha} \forall n \forall y \forall z(\langle n, y\rangle \in \text{€} \wedge \langle n, y\rangle \in \text{€} \rightarrow y = z).$$

**Some consequences of the enumeration postulates:**

24.5.7. Theorem

$$\Vdash^{M} \exists^{=0}x(x \lhd L\,).$$

24.5.8. Theorem

$$\Pi a \Pi b\Big(\textsc{constant}(a)\ \&\ \textsc{constant}(b) \Rightarrow \\ \Vdash^{M} \forall n(n \in \omega \rightarrow \big(\exists^{=n}x(x \lhd a) \wedge \exists^{=n}y(y \lhd b) \rightarrow a = b\big))\Big).$$

24.5.9. Theorem

$$\Pi b\Big[\textsc{constant}(b) \Rightarrow \Vdash^{M} \big(\langle \emptyset, b\rangle \in \text{€} \leftrightarrow \exists^{=0}x(x \lhd b)\big)\Big].$$

24.5.10. Theorem

$$\Pi a, b, c, (\Vdash^{M} \forall n\Big(n \in \omega \rightarrow (\big(\langle n, b\rangle \in \text{€}\ \&\ \exists^{=n}x(x \lhd b)\big) \leftrightarrow \\ \big(\langle \{v|v \in n \vee v = n\}, c\rangle \in \text{€}\ \&\ \exists^{=(n+1)}x(x \lhd c)\big))\Big)).$$

24.5.11. Theorem

$$\Vdash^{M} \forall y \exists n(n \in \omega \wedge \langle n, y\rangle \in \text{€}).$$

*Proof.* As all sets are finite von Neumann ordinals of the meta language, and $\omega$ has the same order as $\Omega$. □

## 24.6 Absolutely all sets are countable

If some set is uncountable, some set of subsets of $\omega$ is uncountable. We have earlier introduced the power set $\mathcal{P}(a) = \{x|x \subset a\}$, and will first consider its import on the question. Thereon we consider the *potency set* of a set a as given by

24.6.1. Definition $\text{♂}(a) = \{x|x = \{y|y \in x \wedge y \in a\}\}$.

The potency set construction is very important in § 25. Here the preoccupation is with showing that neither power sets nor potency sets generate uncountable sets.

### 24.6.1 € restricted to $\mathcal{P}(\{x|x \in \omega\})$

€ restricted to the power set of $\{x|x \in \omega\}$ is

$$\text{€}|_{\mathcal{P}(\{x|x\in\omega\})} = \{(x, y)|(x, y) \in \text{€} \wedge y \in \mathcal{P}(\{x|x \in \omega\})\}, \tag{24.6.2}$$

which has $\omega$ as domain and $\mathcal{P}(\{x|x \in \omega\})$ as range. Given Definitions 24.3.2 and 24.4.1, equation 24.6.2 may be equivalently stated as

$$\text{€}|_{\mathcal{P}(\{x|x\in\omega\})} = \{(x, y)|x \in \omega \wedge y \in \wp^{x}\mathbf{U} \wedge y \in \mathcal{P}(\{x|x \in \omega\})\}. \tag{24.6.3}$$

24.6.4. Theorem

$$\Vdash^{M} \forall u \forall v \forall w((u,v) \in \text{€}|_{\mathcal{P}(\{x|x \in \omega\})} \wedge (u,w) \in \text{€}|_{\mathcal{P}(\{x|x \in \omega\})} \rightarrow v = w).$$

*Proof.* Obvious, from the built up of $\text{€}|_{\mathcal{P}(\{x|x \in \omega\})}$ with orthodox function $\wp^{x}$. □

To attempt Cantor's proof by contradiction for uncountability, assume that $\text{€}|_{\mathcal{P}(\{x|x \in \omega\})}$ surjects from $\omega$ to $\mathcal{P}(\{x|x \in \omega\})$ and posit

24.6.5. Definition

$$S = \{x|x \in \omega \wedge x \notin \text{€}|^{\backprime}_{\mathcal{P}(\{x|x \in \omega\})} x\}.$$

24.6.6. Theorem For an $m \in \omega$,

$$\Vdash^{M} (m,S) \in \text{€}|_{\mathcal{P}(\{x|x \in \omega\})}.$$

*Proof.* A consequence of Equation 24.6.3 and alethic comprehension is

$$\Vdash^{M} (m,S) \in \text{€}|_{\mathcal{P}(\{x|x \in \omega\})} \leftrightarrow \mathcal{T} \ulcorner m \in \omega \wedge S \in \wp^{m}\mathbf{U} \wedge S \in \mathcal{P}(\{x|x \in \omega\}) \urcorner.$$

Let $m \in \omega$ be the natural number such that

$$\Vdash^{M} S \in \wp^{m}\mathbf{U}, \tag{24.6.7}$$

thus

$$\Vdash^{M} m \in \omega \wedge S \in {}^{m}\mathbf{U}. \tag{24.6.8}$$

But besides,

$$\Vdash^{M} S \in \mathcal{P}(\{x|x \in \omega\}), \tag{24.6.9}$$

as

$$\Vdash^{M} S \subset \{x|x \in \omega\}. \tag{24.6.10}$$

$$\Vdash^{M} S \subset \{x|x \in \omega\}.$$

So

$$\Vdash^{M} m \in \omega \wedge S \in \wp^{m}\mathbf{U} \wedge S \in \mathcal{P}(\{x|x \in \omega\}). \tag{24.6.11}$$

Thus, on account of inference mode 7.1.2.1,

$$\Vdash^{M} \mathcal{T}\ulcorner m \in \omega \wedge S \in \wp^{m}\mathbf{U} \wedge S \in \mathcal{P}(\{x|x \in \omega\})\urcorner. \qquad (24.6.12)$$

Finish by using the maxim mode 7.2.7. □

24.6.13. THEOREM

There is an $m \in \omega$ such that $\Vdash^{M} \text{€}|^{\backprime}_{\mathcal{P}(\{x|x\in\omega\})}m \Bumpeq S$.

*Proof.* Invoke Theorems 24.6.4 and 24.6.6, and Definition 24.2.1. □

24.6.14. THEOREM

S is paradoxical.

*Proof.* From Definition 24.6.5 and alethic comprehension,

$$\Vdash^{M} m \in S \leftrightarrow \mathcal{T}\ulcorner m \in \omega \wedge m \notin \text{€}|^{\backprime}_{\mathcal{P}(\{x|x\in\omega\})}m\urcorner. \qquad (24.6.15)$$

Given Definition 24.2.3,

$$\Vdash^{M} m \in S \leftrightarrow \mathcal{T}\ulcorner m \in \omega \wedge m \notin \text{€}|^{\backprime}_{\mathcal{P}(\{x|x\in\omega\})}m\urcorner. \qquad (24.6.16)$$

Given Definition 24.2.1, it follows that

$$\Vdash^{M} m \in S \leftrightarrow \mathcal{T}\ulcorner m \in \omega \wedge m \notin \text{€}|^{\backprime}_{\mathcal{P}(\{x|x\in\omega\})}m\urcorner. \qquad (24.6.17)$$

So that

$$\Vdash^{M} m \in S \leftrightarrow \mathcal{T}\ulcorner m \in \omega \wedge \forall y(\text{€}|^{\backprime}_{\mathcal{P}(\{x|x\in\omega\})}m \Bumpeq y \rightarrow m \notin y)\urcorner. \qquad (24.6.18)$$

Given Theorems 24.6.4 and 24.6.13, and the fact that there is only one $m \in \omega$ such that $S \in \wp^{m}\mathbf{U}$, for the appropriate m, $\forall y(\text{€}|^{\backprime}_{\mathcal{P}(\{x|x\in\omega\})}m \Bumpeq y \leftrightarrow y = S)$. So that

$$\Vdash^{M} m \in S \leftrightarrow \mathcal{T}\ulcorner m \in \omega \wedge m \notin S)\urcorner. \qquad (24.6.19)$$

But it was assumed that $m \in \omega$, which is an orthodox statement, so that

$$\Vdash^{M} m \in S \leftrightarrow \mathcal{T}\ulcorner m \notin S)\urcorner. \qquad (24.6.20)$$

As

$$\Vdash \mathcal{T}\ulcorner m \notin S)\urcorner \rightarrow m \notin S \qquad (24.6.21)$$

and

$$\Vdash \mathrm{m} \in \mathrm{S} \rightarrow \mathcal{T}\ulcorner \mathrm{m} \in \mathrm{S})\urcorner, \tag{24.6.22}$$

it follows that

$$\Vdash \mathrm{m} \notin \mathrm{S} \tag{24.6.23}$$

and

$$\Vdash \mathcal{T}\ulcorner \mathrm{m} \notin \mathrm{S}\urcorner \rightarrow \mathcal{T}\ulcorner \mathrm{m} \in \mathrm{S}\urcorner. \tag{24.6.24}$$

But

$$\Vdash^{\mathrm{M}} \mathcal{T}\ulcorner \mathrm{m} \in \mathrm{S}\urcorner \rightarrow \neg\mathcal{T}\ulcorner \mathrm{m} \notin \mathrm{S}\urcorner, \tag{24.6.25}$$

so that

$$\Vdash \mathcal{T}\ulcorner \mathrm{m} \notin \mathrm{S}\urcorner \rightarrow \neg\mathcal{T}\ulcorner \mathrm{m} \notin \mathrm{S}\urcorner, \tag{24.6.26}$$

and consequently

$$\Vdash \neg\mathcal{T}\ulcorner \mathrm{m} \notin \mathrm{S}\urcorner. \tag{24.6.27}$$

On account of Postulate 7.1.1,

$$\Vdash \neg\mathcal{T}\ulcorner \mathrm{m} \notin \mathrm{S}\urcorner \Rightarrow \Vdash \mathrm{m} \in \mathrm{S}, \tag{24.6.28}$$

so that

$$\Vdash \mathrm{m} \in \mathrm{S}. \tag{24.6.29}$$

A joining of equations 24.6.23 and 24.6.29 results in

$$\Vdash \mathrm{m} \in \mathrm{S} \ \& \Vdash \mathrm{m} \notin \mathrm{S}. \tag{24.6.30}$$

□

Notice that the proof of Equation 24.6.30 merely amounts to an argumentum ad paradoxo, and it has not been proven that $€|_{\mathcal{P}(\{\mathrm{x}|\mathrm{x}\in\omega\})}$ is not a function with domain $\omega$ which is onto its range $\mathcal{P}(\{\mathrm{x}|\mathrm{x} \in \omega\})$.

### 24.6.2 € restricted to $♂(\omega)$

The potency set of $\omega$ is

$$♂(\omega) = \{x|x = \{y|y \in x \wedge y \in \omega\}\}. \tag{24.6.31}$$

€ restricted to the potency set of $\omega$ is

$$€|_{♂(\omega)} = \{(x, y)|(x, y) \in € \wedge y \in ♂(\omega)\}, \tag{24.6.32}$$

which has $\omega$ as domain and $♂(\omega)$ as its range. Given Definition 24.4.1, equation 24.6.32 may be equivalently stated as

$$€|_{♂(\omega)} = \{(x, y)|x \in \omega \wedge y \in \wp^{x}(\mathbf{U}) \wedge y \in ♂(\omega)\}, \tag{24.6.33}$$

24.6.34. Fact

$$€, ♂(\omega) \text{ and } €|_{♂(\omega)} \text{ are orthodox.}$$

*Proof.* € is orthodox given Theorem 24.4.2, $♂(\omega)$ on account of the theory of identity, and $€|_{♂(\omega)}$ is orthodox because € and $♂(\omega)$ are orthodox. □

24.6.35. Fact

$$\Vdash^{M} \forall x \forall y \forall z((x, y) \in €|_{♂(\omega)} \wedge (x, z) \in €|_{♂(\omega)} \rightarrow y = z).$$

*Proof.* As € is functional. □

24.6.36. Assumption

Orthodox function $€|_{♂(\omega)}$ surjects from $\omega$ to $♂(\omega)$:

$$\forall w(w \in ♂(\omega) \rightarrow \exists v(v \in \omega \wedge €|^{`}_{♂(\omega)} v \Bumpeq w)).$$

24.6.37. Definition

$$S = \{x|x \in \omega \wedge x \notin €|^{`}_{♂(\omega)} x\}.$$

24.6.38. Assumption

$$S = \{y|y \in \omega \wedge y \in S\}.$$

24.6.39. Consequence of Assumption 24.6.36

$$S \in ♂(\omega).$$

24.6.40. Assumption S is orthodox.24.6.38.

*Proof.* From Assumption 24.6.38, Axiom 22.0.8 and Theorem 22.0.4. □

24.6.41. Assumption

$$\text{An } m \in \omega \text{ is such that } €|‘_{♂(\omega)}\, m \Bumpeq S.$$

*Proof.* From Assumption 24.6.36. □

24.6.42. Assumption

$\Vdash^{M} \forall x(x \in S \leftrightarrow x \in \omega \land x \notin €|‘_{♂(\omega)} x).$

*Proof.* Given Definition 24.6.37 and the fact that S is orthodox. □

24.6.43. Assumption

$\Vdash^{M} \forall x(x \in S \leftrightarrow x \in \omega \land \forall y(€|‘_{♂(\omega)} x \Bumpeq y \rightarrow x \notin y).$

*Proof.* On account of Definition 24.2.3 and Assumption 24.6.42. □

24.6.44. Assumption

$\Vdash^{M} (m \in S \rightarrow m \notin S).$

*Proof.* It was agreed in Assumption 24.6.41 that for an $m \in \omega$, $€|‘_{♂(\omega)} m \Bumpeq S$. □

24.6.45. Assumption

$m \notin S \rightarrow \exists y(€|‘_{♂(\omega)} m \Bumpeq y \land m \in y).$

*Proof.* From Assumption 24.6.43, the agreement of Assumption 24.6.41. □

24.6.46. Theorem For functional f:

$$\text{if } \Vdash^{M} \exists y(f‘a \Bumpeq y \land a \in y) \text{ and } \Vdash^{M} f‘a \Bumpeq c, \text{then } \Vdash^{M} a \in c.$$

*Proof.* Because $\Vdash^{M} [(a, y) \in f \land (a, c) \in f] \rightarrow y = c$, as f is functional, and because $\Vdash^{M} ((d, e) \in f \leftrightarrow f‘d \Bumpeq e)$ if $\Vdash^{M}$ (f is functional). □

24.6.47. Assumption$\Vdash^{M} (m \notin S \rightarrow m \in S).$

*Proof.* Appeal to Assumption 24.6.45 and Theorem 24.6.46. □

24.6.48. Assumption$\Vdash^{M} m \in S \land m \notin S$

*Proof.* From Assumptions 24.6.44 and 24.6.47. □

The contradiction in the maximal context of Assumption 24.6.48 is false, so it follows that a previous assumption is to be discarded. We do that by stating the following

24.6.49. THEOREMAssumption 24.6.38 is false, and so $\Vdash^{M} S \neq \{x|x \in S \ \wedge \ x \in \omega\}$.

*Proof.* The discussion in § 24.6.2. □

# 25 The theory $\mathfrak{th}\mathcal{H}\mathcal{R}(\mathbf{G})$ of vonsets

> *"Weil die Zermelo'schen Axiome den Bereich B nicht eindeutig bestimmen, ist es sehr unwahrscheinlich, dass mit hilfe dieser Axiome alle Mächtigkeitsprobleme einscheidbar sein sollten. Es ist z. B. sehr wohl möglich, dass das sogenannte Kontinuumproblem, nähmlich ob* $2^{\mathrm{Alef}_0} >$ *oder* $= \mathrm{Alef}_1$ *ist, auf dieser Grundlage überhaupt nicht lösbar ist; es braucht eben nichts darüber entschieden zu sein. Der Sachverhalt kann genau derselbe sein wie im folgenden Falle: Es ist ein unbestimmter Rationalitetsbereich gegeben, und man fragt, ob in diesem Bereiche eine Grösse x vorhanden ist, sodass* $x^2 = 2$ *ist. Dies ist eben nicht bestimmt, wegen der Mehrdeutigkeit des Bereiches."*

---

`Thoralf Albert Skolem in (Skolem 1922, p. 149).`

Recall Definitions 22.0.3 and 22.0.12.

25.0.1. Definition
Set theory $\mathfrak{th}$ is **£** plus Axioms 22.0.6, 22.0.8 and 22.0.10.

25.0.2. Definition
$\mathfrak{th}\mathcal{H}\mathcal{R}(\mathbf{G})$ is $\mathfrak{th} + \mathcal{H}(\mathbf{G}) + \mathcal{R}(\mathbf{G})$, with **G** as in Definition 25.4.6.

Let **NBGC + TA** be Neumann-Bernays-Gödel set theory with Global Choice and Tarski's Axiom. An interpretation of **NBGC + TA** is developed in $\mathfrak{th}\mathcal{H}\mathcal{R}(\mathbf{G})$ below.

Natural weakenings and extensions of **NBGC + TA** are as well taken to be theories of *vonsets*. Needless to say, but all vonsets are sets, though some sets are not vonsets.

The term "natural" in the previous paragraph is left undefined, as investigations should not be restrained. So we here disregard philosophical quandaries s that the term "vonset" may have different natural extensions of **NBG** which are not consistent with each other.

## 25.1 The potency vonset

We saw in § 18 that power sets are mathematically useless, as they are paradoxical lest of a non-paradoxical universal set.

The notion of potency set was introduced in Definition 24.6.1, which states:

$$♂(a) == \{x|x = \{y|y \in x \wedge y \in a\}\}.$$

Potency vonsets are potency sets, as all vonsets are sets.

25.1.1. Theorem

The *Potency vonset* of vonset a contains precisely a's *hypovonsets*, as per Definition 22.0.5.

*Proof.* Use Axioms 22.0.10 and 25.4.10 and Theorem 22.0.9.entail that vonsets are heritors, and from Axioms 22.0.6 and 22.0.8. □

25.1.2. Theorem

$$\male(a) \text{ is orthodox, and all of its members are hereditarily heritors.}$$

*Proof.* $\male(a)$ is orthodox by the logic of identity. Its members, if any, are heritors on account of Axiom 22.0.8, and are hereditarily heritors given Axiom 22.0.10. □

25.1.3. Theorem

$$\male(a) \text{ is empty if a is not an heritor.}$$

*Proof.* Appeal to Axiom 22.0.8. □

25.1.4. Theorem

$$\Vdash^{M} \forall x(x \in \male(a) \leftrightarrow \mathcal{H}(x) \wedge \mathcal{H}(a) \wedge x \subset a).$$

*Proof.* Appeal to Theorem 22.0.9 and Definition 24.6.1. □

## 25.2 The Grothendieck vonset of w relative to v

25.2.1. Definition

$$\begin{aligned}
G(v, w, v_0, v_1,) == \forall y\big(& w \in y \wedge \forall z[z \in y \rightarrow (z \in \male(v_1) \wedge \male(z) \in \male(v_1) \wedge \male(z) \in y)] \wedge \\
& \forall z(z \in \male(y) \wedge z \notin y \rightarrow \exists f[f \in v \wedge \text{Bijection}(f) \wedge \\
& (\forall x_0)(x_0 \in y \rightarrow \exists x_1(x_1 \in z \wedge (x_0, x_1) \in f)) \\
& (\forall x_1)(x_1 \in z \rightarrow \exists x_0(x_0 \in y \wedge (x_0, x_1) \in f))]) \rightarrow v_0 \in y\big).
\end{aligned}$$

Use Theorem 17.1.4 to obtain the manifestation set with parameters $\mathcal{G}(v, w)$,

25.2.2. Definition of the Grothendieck of w relative to v

$$\begin{aligned}
\Vdash^{M} \forall u \forall w(u \in \mathcal{G}(v, w) \leftrightarrow \mathcal{T}\mathcal{T} \forall y\big(& w \in y \wedge \forall z[z \in y \rightarrow \\
& (z \in \male(\mathcal{G}(v, w)) \wedge \male(z) \in \male(\mathcal{G}(v, w)) \wedge \male(z) \in y)] \wedge \\
& \forall z(z \in \male(y) \wedge z \notin y \rightarrow \exists f[f \in v \wedge \text{Bijection}(f) \wedge \\
& (\forall x_0)(x_0 \in y \rightarrow \exists x_1(x_1 \in z \wedge (x_0, x_1) \in f)) \\
& (\forall x_1)(x_1 \in z \rightarrow \exists x_0(x_0 \in y \wedge (x_0, x_1) \in f))]) \rightarrow u \in y\big))
\end{aligned}$$

25.2.3. Theorem

$\mathcal{G}(v, w)$ is orthodox for orthodox v and w, so that

$$\begin{aligned}
\Vdash^{M} \forall u \forall w(u \in \mathcal{G}(v, w) \leftrightarrow \forall y\big(& w \in y \wedge \forall z[z \in y \rightarrow \\
& (z \in \male(\mathcal{G}(v, w)) \wedge \male(z) \in \male(\mathcal{G}(v, w)) \wedge \male(z) \in y)] \wedge \\
& \forall z(z \in \male(y) \wedge z \notin y \rightarrow \exists f[f \in v \wedge \text{Bijection}(f) \wedge \\
& (\forall x_0)(x_0 \in y \rightarrow \exists x_1(x_1 \in z \wedge (x_0, x_1) \in f)) \\
& (\forall x_1)(x_1 \in z \rightarrow \exists x_0(x_0 \in y \wedge (x_0, x_1) \in f))]) \rightarrow u \in y\big))
\end{aligned}$$

*Proof.* As in the proof that $\omega$ is orthodox, of Theorem 11.0.2.3 on page 47, noting that ♂$(\mathcal{G}(v, w))$ is an orthodox heritor by cause of Theorem 25.1.2. □

25.2.4. Remark For appropriate v and w, Theorem 25.2.3 amounts to *Tarski's axiom*, which states that all sets are members of a Grothendieck-universe. Tarski-Grothendieck set theory is usually presented as **ZFC** + Tarski's axiom.

## 25.3 Capture

In this section we presuppose that the sets and conditions invoked are orthodox.

25.3.1. Definition Capture with B from w

$$\mathcal{C}(B, w) =\!\!=\!\!= \{x | \exists y(y \in w \land \forall z((x, y)^z_y \in B \leftrightarrow y = z))\}$$

25.3.2. Theorem

Capture is equivalent with replacement.

*Proof.* i) If a vonset is obtained from capture with B from w, it can be obtained from replacement by using the functional condition $\forall z((x, y)^z_y \in B \leftrightarrow y = z)$. ii) If a vonset is obtained from replacement by functional B so that $\forall x \forall y \forall z((x, y) \in B \land (x, z) \in B \rightarrow y = z)$, it can be obtained from capture by using the condition as in Definition 25.3.1. □

25.3.3. Theorem

Capture, as replacement, entails specification.

*Proof.* Use the functional $B' =\!\!=\!\!= \{(x, y) | x \in B \land x = y)\}$ as capture vonset relative to a vonset a, and observe that the existence of the vonset $\{x | x \in a \land B(x)\}$ is justified by capture and extensionality, which holds for **V** and **G** below, as per Theorem 25.4.13. □

## 25.4 V and G

25.4.1. Definition of the generation of u

$$\begin{aligned}
\mathfrak{G}(u) = \Big\{ w | w \in u \lor \forall v([&u \in v \land E = \{(x_i, x_j) | (x_i, x_j) \in u^2 \land x_i \in x_j\} \in v \land \\
&\forall x_i(x_i \in v \rightarrow \{y | y \in u \land y \notin x_i\} \in v) \land \forall x_i \forall x_j(x_i \in v \land x_j \in v \rightarrow x_i \cap x_j \in v) \land \\
&\forall x_i(x_i \in v \rightarrow \mathrm{dom}(x_i) = \{y | \exists x((y, x) \in x_i)\} \in v) \land \\
&\forall x_i(x_i \in v \rightarrow \{y | \exists x_j, x_k(y = (x_j, x_k) \land x_j \in x_i \land v_k \in v_1)\} \in v) \land \\
&\forall x_i(x_i \in v \rightarrow \{y | \exists x_j \exists x_k \exists x_l(y = (x_j, x_k, x_l) \land (x_k, x_l, x_j) \in x_i\} \in v) \land \\
&\forall x_i(x_i \in v \rightarrow \{y | \exists x_j \exists x_k \exists x_l(y = (x_j, x_k, x_l) \land (x_j, x_l, x_k) \in x_i\} \in v)] \rightarrow w \in v) \Big\}
\end{aligned}$$

25.4.2. Definition of $V(v_0, v_1)$

$$
\begin{aligned}
V(v_0, v_1) = \forall v(&(\omega \in v \wedge \forall w \in v \forall x \in v : \{w, x\} \in v \wedge \forall w \in v : \bigcup w \in v \wedge \\
&\forall w \in v : \mars(w) = \{x | x = \{y | y \in x \wedge y \in w\}\} \in v \wedge \\
&\forall w \in v : \wp(w) = \{x | (x \in w \wedge \forall y (y \in w \rightarrow x \trianglelefteq y)\} \in v \wedge \\
&\forall w \in v : \mathcal{G}(\mathfrak{D}(v_1), w) \in v \wedge \\
&\forall w \in v \forall B \in \mathfrak{D}(v_1) : \mathcal{C}(B, w) = \{x | \exists y (y \in w \wedge \forall z ((x, y)^z_y \in B \leftrightarrow y = z))\} \in v) \\
&\rightarrow v_0 \in v)
\end{aligned}
$$

25.4.3. Definition of $\mathbf{V}$ via manifestation from Definition 25.4.2

$$
\begin{aligned}
\Vdash^M \forall u\Big[u \in \mathbf{V} \leftrightarrow \mathcal{T}\mathcal{T} \forall v[&[\omega \in v \wedge \forall w \in v \forall x \in v : \{w, x\} \in v \wedge \forall w \in v : \bigcup w \in v \wedge \\
&\forall w \in v : \mars(w) = \{x | x = \{y | y \in x \wedge y \in w\}\} \in v \wedge \\
&\forall w \in v : \wp(w) = \{x | (x \in w \wedge \forall y (y \in w \rightarrow x \trianglelefteq y)\} \in v \wedge \\
&\forall w \in v : \mathcal{G}(\mathfrak{D}(\mathbf{V}), w) \in v \wedge \\
&\forall w \in v \forall B \in \mathfrak{D}(\mathbf{V}) : \mathcal{C}(B, w) = \{x | \exists y (y \in w \wedge \forall z ((x, y)^z_y \in B \leftrightarrow y = z))\} \in v] \\
&\rightarrow u \in v]\Big]
\end{aligned}
$$

As $\mathbf{V}$ is orthodox on account of Theorems 22.0.4 and 25.4.11,

25.4.4. Theorem

$$
\begin{aligned}
\Vdash^M \forall u\Big[u \in \mathbf{V} \leftrightarrow \forall v[&[\omega \in v \wedge \forall w \in v \forall x \in v : \{w, x\} \in v \wedge \forall w \in v : \bigcup w \in v \wedge \\
&\forall w \in v : \mars(w) = \{x | x = \{y | y \in x \wedge y \in w\}\} \in v \wedge \\
&\forall w \in v : \wp(w) = \{x | (x \in w \wedge \forall y (y \in w \rightarrow x \trianglelefteq y)\} \in v \wedge \\
&\forall w \in v : \mathcal{G}(\mathfrak{D}(\mathbf{V}), w) \in v \wedge \\
&\forall w \in v \forall B \in \mathfrak{D}(\mathbf{V}) : \mathcal{C}(B, w) = \{x | \exists y (y \in w \wedge \forall z ((x, y)^z_y \in B \leftrightarrow y = z))\} \in v] \\
&\rightarrow u \in v]\Big]
\end{aligned}
$$

25.4.5. The genus equation:

$$\mathbf{G} = \mathfrak{G}(\mathbf{V}).$$

The Genus equation dictates that Genus $\mathbf{G}$ is the *generation* of $\mathbf{V}$ in the sense of Definition 25.4.1, so we can express the following more accurate definition.

25.4.6. Definition of the Genus, by means of the Genus equation and Definition 25.4.1:

$$
\begin{aligned}
\mathbf{G} = \Big\{ & w | w \in \mathbf{V} \vee \forall v([\mathbf{V} \in v \wedge E = \{(x_i, x_j) | (x_i, x_j) \in \mathbf{V}^2 \wedge x_i \in x_j\} \in v \wedge \\
& \forall x_i (x_i \in v \to \{u | u \in \mathbf{V} \wedge u \notin x_i\} \in v) \wedge \forall x_i \forall x_j (x_i \in v \wedge x_j \in v \to x_i \cap x_j \in v) \wedge \\
& \forall x_i (x_i \in v \to \mathrm{dom}(x_i) = \{w | \exists x((w, x) \in x_i)\} \in v) \wedge \\
& \forall x_i (x_i \in v \to \{u | \exists x_j, x_k (u = (x_j, x_k) \wedge x_j \in x_i \wedge v_k \in \mathbf{V})\} \in v) \wedge \\
& \forall x_i (x_i \in v \to \{u | \exists x_j \exists x_k \exists x_l (u = (x_j, x_k, x_l) \wedge (x_k, x_l, x_j) \in x_i\} \in v) \wedge \\
& \forall x_i (x_i \in v \to \{u | \exists x_j \exists x_k \exists x_l (u = (x_j, x_k, x_l) \wedge (x_j, x_l, x_k) \in x_i\} \in v)] \to w \in v) \Big\}
\end{aligned}
$$

25.4.7. Fact

$$\Vdash^{M} \mathbf{V} \subset \mathbf{G}.$$

25.4.8. Definition of **V**, with recourse to **G**:

$$
\begin{aligned}
\mathbf{V} \Longleftrightarrow \Big\{ & u | \forall v([\omega \in v \wedge \forall w \in v \forall x \in v : \{w, x\} \in v \wedge \forall w \in v : \bigcup w \in v \wedge \\
& \forall w \in v : \male(w) = \{x | x = \{y | y \in x \wedge y \in w\}\} \in v \wedge \\
& \forall w \in v : \female(w) = \{x | (x \in w \wedge \forall y (y \in w \to x \trianglelefteq y)\} \in v \wedge \\
& \forall w \in v : \mathcal{G}(\mathbf{G}, w) \in v \wedge \\
& \forall w \in v \forall B \in \mathbf{G} : \mathcal{C}(B, w) = \{x | \exists y (y \in w \wedge \forall z ((x, y)_y^z \in B \leftrightarrow y = z))\} \in v] \\
& \to u \in v) \Big\}.
\end{aligned}
$$

25.4.9. Definition

**V** is the class of all vonsets.

25.4.10. Axiom

$$\mathcal{H}(\mathbf{G}).$$

25.4.11. Theorem

$$\mathcal{H}(\mathbf{V}).$$

*Proof.* On account of Axioms 22.0.10 and 25.4.10. □

25.4.12. Corollary

**V** and **G** are orthodox.

*Proof.* Use Axiom 25.4.10, Theorem 25.4.11 and Theorem 22.0.4. □

25.4.13. Theorem

Co-extensional members of **G** are identical.

*Proof.* Use Axioms 22.0.10 and 25.4.10, Theorem 25.4.11 and, finally, Theorem 22.0.7. □

25.4.14. Theorem (**V** is a hyposet of **G** in the sense of Definition 22.0.5)

$$\mathbf{V} = \{x|x \in \mathbf{V} \wedge x \in \mathbf{G}\}.$$

*Proof.* As $\mathcal{H}(\mathbf{V})$ and $\mathcal{H}(\mathbf{G})$ on account of Axiom 25.4.10 and Theorem 25.4.11, appeals to Theorems 22.0.9 and Fact 25.4.7 suffice to finish the proof. □

25.4.15. Axiom (The Genus is wellfounded)

$$\mathcal{R}(\mathbf{G})$$

25.4.16. Theorem

All classes are wellfounded.

*Proof.* Invoke the result of Exercise 22.0.13. □

25.4.17. Theorem

All vonsets are wellfounded

*Proof.* Given Fact 25.4.7, a vonset in **V** is as well a class member of **G**. So the vonset is wellfounded on account of Theorem 25.4.16. □

25.4.18. Theorem

**G** is not a class.

*Proof.* If **G** were a class, it would, on account of Definition 25.4.6, follow that $\mathbf{G} \in \mathbf{G}$, which contradicts Axiom 25.4.15. □

25.4.19. Remark

Instead of postulating Axiom 25.4.15, one may obtain a suitable regular class **V*** of all regular vonsets by taking it to be the class of all elements of a potency set of an ordinal in **V**. That invokes the consistency proof of **ZFC** with regularity given the consistency of $\mathbf{ZFC}^{-} = \mathbf{ZFC}$ without regularity, by (Kunen 1980, chapter 3), or a similar relative consistency proof. Given Kunen's result, however, and the relative consistency results obtained earlier by (Skolem 1922) and (Neumann 1929), we know that we can safely posit Axiom 25.4.15.

25.4.20. Theorem

**V** is not a vonset.

*Proof.* Appeal to Definition 25.4.9 and Theorem 25.4.15. □

## 25.5 Primitive theorems for classes

We leave is as an exercise to prove the following from Definition 25.4.6.

25.5.1. THEOREM (**V** IS A CLASS)

$$\mathbf{V} \in \mathbf{G}.$$

25.5.2. THEOREM (MEMBERSHIP CLASS)

$$E = \{(x, y) | x \in \mathbf{V} \land y \in \mathbf{V} \land x \in y\} \in \mathbf{G}.$$

25.5.3. THEOREM (INTERSECTION CLASS)

$$\forall A \in \mathbf{G} \forall B \in \mathbf{G} \exists C \in \mathbf{G} \forall x (x \in C \leftrightarrow x \in A \land x \in B).$$

25.5.4. THEOREM (COMPLEMENT CLASS)

$$\forall A \in \mathbf{G} \exists B \in \mathbf{G} \forall x (x \in B \leftrightarrow x \notin A).$$

25.5.5. THEOREM (DOMAIN CLASS)

$$\forall A \in \mathbf{G} \exists B \in \mathbf{G} \forall x (x \in B \leftrightarrow \exists y ((x, y) \in A)).$$

25.5.6. THEOREM (PRODUCT BY V CLASS)

$$\forall A \in \mathbf{G} \exists B \in \mathbf{G} \forall x (x \in B \leftrightarrow \exists y \exists z (x = (y, z) \land y \in A \ \land z \in V)).$$

25.5.7. THEOREM (CIRCULAR PERMUTATION CLASS)

$$\forall A \in \mathbf{G} \exists B \in \mathbf{G} \forall x \forall y \forall z ((x, y, z) \in B \leftrightarrow (y, z, x) \in A).$$

25.5.8. THEOREM (TRANSPOSITION CLASS)

$$\forall A \in \mathbf{G} \exists B \in \mathbf{G} \forall x \forall y \forall z ((x, y, z) \in B \leftrightarrow (x, z, y) \in A).$$

## 25.6 The Tuple-lemmas

25.6.1. LEMMA

$$\forall A \in \mathbf{G} \exists B_1 \in \mathbf{G} \forall x \forall y \forall z ((x, y, z) \in B_1 \leftrightarrow (x, y) \in A \land z \in V).$$

25.6.2. LEMMA

$$\forall A \in \mathbf{G} \exists B_2 \in \mathbf{G} \forall x \forall y \forall z ((x, z, y) \in B_2 \leftrightarrow (x, y) \in A \land z \in V).$$

25.6.3. LEMMA

$$\forall A \in \mathbf{G} \exists B_3 \in \mathbf{G} \forall x \forall y \forall z ((z, x, y) \in B_3 \leftrightarrow (x, y) \in A \land z \in V).$$

25.6.4. LEMMA

$$\forall A \in \mathbf{G} \exists B_4 \in \mathbf{G} \forall x \forall y ((y, x) \in B_4 \leftrightarrow (x, y) \in A).$$

*Proof.* Use Theorem 25.5.6 to get $B_1$, Theorem 25.5.8 on $B_1$ to get $B_2$, Theorem 25.5.7 on $B_1$ to get $B_3$, and use Theorem 25.5.7 on $B_2$, plus Theorem 25.5.5, to get $B_4$. □

## 25.7 The class existence theorem

## 25.8 The expansion lemma

## 25.9 Proof that V is orthodox

## 25.10 Proof that all members of V are orthodox

As $\mathcal{H}$ by Axiom ... This is done already.

## 25.11 Global well ordering

Useful explanation of

Global well ordering given global choice.

# 26 The theory ᵺ$\mathcal{H}\mathcal{R}$(W) of vansets

> *"The analogy between the myth of mathematics and the myth of physics is, in some additional and perhaps fortuitous ways, strikingly close. Consider, for example, the crisis which was precipitated in the foundations of mathematics, at the turn of the century, by the discovery of Russell's paradox and other antinomies of set theory. These contradictions had to be obviated by unintuitive, ad hoc devices; our mathematical myth-making became deliberate and evident to all. But, what, of physics? An antinomy arose between the undular and the corpuscular accounts of light; and if this was not as out-and-out a contradiction as Russell's paradox, I suspect that the reason is that physics is not as out-and-out as mathematics."*
>
> `Willard van Orman Quine`, in (Quine 1961, pp. 18–19)

We give an account of Willard van Quine's set theory NF, of (Quine 1937), via the axiomatization offered by (Hailperin 1944, p. 10), which is adapted here:

$P_0 : \exists\beta\forall x(x \in \beta \leftrightarrow \exists y(x \in y \wedge x \notin y)$

$P_1 : \forall u\forall v\exists\beta\forall x(x \in \beta \leftrightarrow (x \notin u \vee x \notin v))$

$P_2 : \forall\alpha\exists\beta\forall x\forall y((\{x\}, \{y\}) \in \beta \leftrightarrow (x, y) \in \alpha)$

$P_3 : \forall\alpha\exists\beta\forall x\forall y\forall z((x, y, z) \in \beta \leftrightarrow (x, y) \in \alpha)$

$P_4 : \forall\alpha\exists\beta\forall x\forall y\forall z((x, z, y) \in \beta \leftrightarrow (x, y) \in \alpha)$

$P_5 : \forall\alpha\exists\beta\forall x\forall y((y, x) \in \beta \leftrightarrow x \in \alpha)$

$P_6 : \forall\alpha\exists\beta(x \in \beta \leftrightarrow \forall u((u, \{x\}) \in \alpha))$

$P_7 : \forall\alpha\exists\beta\forall x\forall y((y, x) \in \beta \leftrightarrow (x, y) \in \alpha)$

$P_8 : \exists\beta\forall x(x \in \beta \leftrightarrow \exists y(x = \{y\}))$

$P_9 : \exists\beta\forall x\forall y((\{x\}, y) \in \beta \leftrightarrow x \in y)$

Notice that $P_0$ was not included in (Hailperin 1944, p. 10).

**U** was reserved for the full universal set $\{x|x = x\}$ of **£** in Definition 18.0.5.

**V** was in §25 reserved for the class of all vonsets, as defined via manifestation there.

**W**, with associated mnemonic device *die* **W***elt*, is reserved for the Quinean New-Foundations-style vanset of all vansets, as defined via manifestation below in this section.

26.0.1. Definition

$$
\begin{aligned}
\mathrm{W}(v_0, v_1) \mathrel{=\!=} \forall v\Big(&[\{x|\exists y(x \in y \wedge x \notin y)\} \in v \wedge \\
&\forall w \forall x(w \in v \wedge x \in v \to \{y \in v_1|\,(y \notin w \vee y \notin x)\} \in v) \wedge \\
&\forall w(w \in v \to \{(\{x\}, \{y\} \in v_1|(x, y) \in w\} \in v) \wedge \\
&\forall w(w \in v \to \{(x, y, z) \in v_1|(x, y) \in w\} \in v) \wedge \\
&\forall w(w \in v \to \{(x, z, y) \in v_1|(x, y) \in w\} \in v) \wedge \\
&\forall w(w \in v \to \{(y, x) \in v_1|(x, y) \in w\} \in v) \wedge \\
&\forall w(w \in v \to \{x \in v_1|\forall y(y \in v_1 \to (y, \{x\}) \in w)\} \in v) \wedge \\
&\forall w(w \in v \to \{(y, x) \in v_1|x \in w\} \in v) \wedge \\
&\forall w(w \in v \to \{x \in v_1|\exists y(y \in v_1 \wedge x = \{y\})\} \in v) \wedge \\
&\forall w(w \in v \to \{(\{x\}, y) \in v_1|x \in y\} \in v)] \\
&\to v_0 \in v\Big)
\end{aligned}
$$

Use Definitions 26.0.1 and 17.1.2 to obtain

26.0.2. Theorem

$$
\begin{aligned}
\forall u(u \in \mathbf{W} \leftrightarrow \mathcal{T}\mathcal{T}\forall v\Big(&[\{x|\exists y(x \in y \wedge x \notin y)\} \in v \wedge \\
&\forall w \forall x(w \in v \wedge x \in v \to \{y \in \mathbf{W}|(y \notin w \vee y \notin x)\} \in v) \wedge \\
&\forall w(w \in v \to \{(\{x\}, \{y\} \in \mathbf{W}|(x, y) \in w\} \in v) \wedge \\
&\forall w(w \in v \to \{(x, y, z) \in \mathbf{W}|(x, y) \in w\} \in v) \wedge \\
&\forall w(w \in v \to \{(x, z, y) \in \mathbf{W}|(x, y) \in w\} \in v) \wedge \\
&\forall w(w \in v \to \{(y, x) \in \mathbf{W}|(x, y) \in w\} \in v) \wedge \\
&\forall w(w \in v \to \{x \in \mathbf{W}|\forall y(y \in \mathbf{W} \to (y, \{x\}) \in w)\} \in v) \wedge \\
&\forall w(w \in v \to \{(y, x) \in \mathbf{W}|x \in w\} \in v) \wedge \\
&\forall w(w \in v \to \{x \in \mathbf{W}|\exists y(y \in \mathbf{W} \wedge x = \{y\})\} \in v) \wedge \\
&\forall w(w \in v \to \{(\{x\}, y) \in \mathbf{W}|x \in y\} \in v)] \\
&\to u \in v\Big))
\end{aligned}
$$

26.0.3. Theorem

**W** is orthodox.

*Proof.* Adapt the the proof of Theorem 11.0.2.3. □

26.0.4. Theorem

$$
\begin{aligned}
\forall u(u \in \mathbf{W} \leftrightarrow \forall v\Big([\{x|\exists y(x \in y \wedge x \notin y)\} \in v \wedge \\
&\forall w\forall x(w \in v \wedge x \in v \to \{y \in \mathbf{W}|(y \notin w \vee y \notin x)\} \in v) \wedge \\
&\forall w(w \in v \to \{(\{x\},\{y\} \in \mathbf{W}|(x,y) \in w\} \in v) \wedge \\
&\forall w(w \in v \to \{(x,y,z) \in \mathbf{W}|(x,y) \in w\} \in v) \wedge \\
&\forall w(w \in v \to \{(x,z,y) \in \mathbf{W}|(x,y) \in w\} \in v) \wedge \\
&\forall w(w \in v \to \{(y,x) \in \mathbf{W}|(x,y) \in w\} \in v) \wedge \\
&\forall w(w \in v \to \{x \in \mathbf{W}|\forall y(y \in \mathbf{W} \to (y,\{x\}) \in w)\} \in v) \wedge \\
&\forall w(w \in v \to \{(y,x) \in \mathbf{W}|x \in w\} \in v) \wedge \\
&\forall w(w \in v \to \{x \in \mathbf{W}|\exists y(y \in \mathbf{W} \wedge x = \{y\})\} \in v) \wedge \\
&\forall w(w \in v \to \{(\{x\},y) \in \mathbf{W}|x \in y\} \in v)] \\
&\to u \in v\Big))
\end{aligned}
$$

*Proof.* A consequence of Theorem 26.0.2 as **W** is orthodox, given Theorem 26.0.3. □

26.0.5. Axiom

$$\mathcal{H}(\mathbf{W})$$

26.0.6. Theorem

Co-extensional sets in **W** are identical.

*Proof.* Use Axiom 22.0.10 and Theorem 22.0.7. □

The identity for **W** is of course given by

26.0.7. Definition

$$a \overset{\mathbf{w}}{=} b \Longleftrightarrow \forall v(v \in \mathbf{W} \to (a \in v \to b \in v)).$$

By Axiom 26.0.5, Theorem 26.0.6 and Theorem 26.0.4 combined with the results of (Hailperin 1944), it follows that ŧƀ$\mathcal{H}$**W** accounts for Quine's set theory **NF**.

# 27 Whither librationist category theory?

> *"La filosofia è scritta in questo grandissimo libro, che continuamente ci sta aperto innanzi agli occhi (io dico l'Universo), ma non si può intendere, se prima non il sapere a intender la lingua, e conoscer i caratteri ne quali è scritto. Egli è scritto in lingua matematica, e i caratteri son triangoli, cerchi ed altre figure geometriche, senza i quali mezzi è impossibile intenderne umanamente parola; senza questi è un aggirarsi vanamente per un oscuro labirinto."*
>
> Galileo Galilei in (Galilei 1623, pp. 16–17).

Some have claimed that set theories as Grothendieck set theory **NBGC + TA** are ideal for category theory, and so it seems worthwhile to investigate how category theory can best be done in the librationist framework set up for mentioned set theories in § 25.

Cambridge University Press.

# Index of names